\documentclass[11pt,english]{article}
\usepackage[sc]{mathpazo}
\usepackage{geometry}
\usepackage{color}
\usepackage{array}
\usepackage{bbding}
\usepackage{multirow}
\usepackage[fleqn]{amsmath}
\usepackage{amssymb}
\usepackage{amsthm}
\usepackage{graphicx}

\usepackage{natbib}
\usepackage{lscape}
\usepackage[colorlinks,
linkcolor=blue, 
anchorcolor=blue, 
citecolor=blue,  
]{hyperref}
\usepackage{comment}
\usepackage{bm}
\usepackage[title]{appendix}
\usepackage{longtable}
\usepackage{mathtools}
\usepackage{chngcntr}
\usepackage{authblk}
\usepackage{enumerate}

\makeatletter
\allowdisplaybreaks

\@ifundefined{showcaptionsetup}{}{%
 \PassOptionsToPackage{caption=false}{subfig}}
\usepackage{subfig}
\makeatother

\usepackage{babel}
\usepackage{enumitem}
\usepackage{enumerate}

\begin{document}\sloppy

\title{Network-level rhythmic control of heterogeneous automated traffic with buses}

\author{Xiangdong Chen\textsuperscript{a}\hspace{2em} Xi Lin\textsuperscript{a}\hspace{2em} Meng Li\textsuperscript{a}\footnote{Corresponding author. E-mail address: \textcolor{blue}{mengli@.tsinghua.edu.cn}}\hspace{2em} Fang He\textsuperscript{b}} 
\affil{\small\emph{\textsuperscript{a}Department of Civil Engineering, Tsinghua University, Beijing 100084, P.R. China}\normalsize}
\affil{\small\emph{\textsuperscript{b}Department of Industrial Engineering, Tsinghua University, Beijing 100084, P.R. China}\normalsize}

\date{\today}
\maketitle

\begin{abstract}
\noindent Guaranteeing the quality of transit service is of great importance to promote the attractiveness of buses and alleviate urban traffic issues such as congestion and pollution. Emerging technologies of automated driving and V2X communication have the potential to enable the accurate control of vehicles and the efficient organization of traffic to enhance both the schedule adherence of buses and the overall network mobility. This study proposes an innovative network-level control scheme for heterogeneous automated traffic composed of buses and private cars under a full connected and automated environment. Inheriting the idea of network-level rhythmic control proposed by  Lin et al. (2020), an augmented rhythmic control scheme for heterogeneous traffic, i.e., RC-H, is established to organize the mixed traffic in a rhythmic manner. Realized virtual platoons are designed for accommodating vehicles to pass through the network, including dedicated virtual platoons for buses to provide exclusive right-of-ways (ROWs) on their trips and regular virtual platoons for private cars along with an optimal assignment plan to minimize the total travel cost. A mixed-integer linear program (MILP) is formulated to optimize the RC-H scheme and a bilevel heuristic solution method is designed to relieve the computational burden of MILP. Numerical examples and simulation experiments are conducted to evaluate the performance of the RC-H scheme under different scenarios. The results show that the bus operation can be guaranteed and the travel delay can be minimized under various demand levels with transit priority. Moreover, compared with traffic signal control strategies, the RC-H scheme has significant advantages in handling massive traffic demand, in terms of both vehicle delay and network throughput.
\par
\hfill\break%
\noindent\textit{Keywords}: rhythmic control, heterogeneous traffic, connected and automated vehicles, bus
\end{abstract}

\section{Introduction} \label{intro_sec}

The rapid growth of metropolitan population and car ownership has aggravated the problem of traffic congestion and air pollution. With the ability of serving mass travel demands with a relatively low cost, bus transit is regarded as an effective approach to relieve the negative impacts of the enormous private car usage  (Levinson et al., 2003; Deng and Nelson, 2011). To improve the attractiveness and encourage the preference choice of buses compared to private cars, a variety of strategies for granting the bus system's quality of service have been proposed and implemented. The crux for guaranteeing the service quality is to alleviate the travel delay of buses, particularly at intersection, thereby enhancing their on-time performance and schedule adherence.\par

One prevailing strategy is the transit signal priority (TSP), which provides a priority right-of-way (ROW) for buses to pass through a signalized intersection (Smith et al., 2005). TSP can be roughly divided into three main categories, i.e., passive, active, and adaptive (Baker et al., 2002). Passive TSP reduces the intersection delay of buses by cycle length reduction or phase split strategy, and operates whether or not a bus is present; thus, it inevitably induces continuous waste of road capacity along the arterial road and is more suitable for application to the scenario where the buses operate frequently and predictably (Urbanik et al., 1977).  Active TSP, such as green extension and red truncation, is triggered only when a signal priority request is generated (Baker et al., 2002); it shows promising benefits for the transit service when the traffic demand is not high, but it potentially causes significant disruptions to the regular traffic; thus, the absolute personal travel time savings of buses may be compromised or even 
neutralized by the extra delays of regular vehicles (Garrow and Machemehl, 1999; Sunkari et al., 1995; Balke et al., 2000). By contrast, adaptive TSP considers both the service quality of buses and the impact on the regular traffic  (Dion  and  Hellinga,  2002;  Liao  and  Davis,  2007;  Hu  et  al.,  2014;  Chang  et  al.,1996;  Li  et  al.,  2011). It requires the acquisition of real-time information of the approaching buses through automatic vehicle location systems or connected vehicle technology, followed by concise prediction and continuous updating of their arrival times. Personal delay is proposed to normalize the bus delay into passengers' delay and is applied as a performance index to balance the mobility enhancement of buses and regular vehicles, minimizing the efficiency damage to the regular traffic under the premise of bus priority fulfillment. However, these studies have mainly focused on the design of adaptive TSP at an isolated intersection, where the advanced green time of buses at the upstream intersection may be compensated by the extra delay at the downstream intersection. To guarantee the effectiveness of the priority treatments at intersections, coordinated TSP strategy is then proposed to enable the cooperation of the signal controls between adjacent intersections  (He et al., 2012; He et al., 2014; Christofa et al., 2016; Ma et al., 2013; Hu et al., 2015). \par

To further enhance the performance of bus service, the TSP strategy is always bound to the implementation of dedicated bus lanes (DBLs). DBL provides an exclusive ROW for buses along their traveling, and the interference from regular traffic and the mobility damage caused by the traffic flow, especially in congestion situations, can be alleviated ated (Levinson et al., 2003; Deng andNelson, 2011;  Currie, 2006). However, exclusive lanes for bus transit will inevitably invade the road resources for regular traffic and substantially reduce the road capacity, potentially degrading the overall traffic efficiency. Intermittent bus lane (IBL) is further proposed to relieve this issue such that it is transferred to a DBL upon the appearance of a bus and is switched back to serve regular traffic upon the bus leaving the lane (Viegas and Lu,1997) . Furthermore, combining the implementation of IBL with the design of TSP strategy advances the improvement of both the service quality of buses and the efficiency of whole traffic, compared to the pure TSP strategy or DBL implementation at an isolated intersection ion (Eichler and Daganzo, 2006; Ma et al., 2014), and the benefit becomes more significant under an integrated control scheme on the network scale  (Viegas and Lu, 2004; Mesbah et al., 2008). \par

The benefit of implementing IBL implies that a permanent DBL is not indispensable for guaranteeing the operational performance of buses. Instead, a series of dedicated time-space slots are sufficient in the view of the infrequent presence of buses, and are more embraced for the whole-network mobility. The achievement of IBL benefits from the technology development of automatic vehicle location or wireless communication, which enables the real-time information of buses to be available, making their arrival times predictable. Furthermore, with the emergence of automated driving technology, vehicles can be controlled more accurately and organized centrally (Administration et al., 2013), enabling the operation of not only transit buses but also the whole traffic to be more predictable and reliable. Thus, leveraging the technologies of connected and automated vehicles (CAVs) has great potential to fulfill a more efficient and credible priority treatment of transit buses. Different from the approach of reserving DBLs for buses to guarantee the service quality, we propose an innovative idea of assigning “dedicated time-space slots” for buses under a CAV environment, which provides exclusive ROWs for buses to eliminate the potential delays along their traveling in the network. \par

Moreover, although the coordinated control strategies in an arterial or a network scale handle the problem of the transit priority effectiveness between adjacent intersections, most of these are only effective in the under-saturated conditions, while over-saturated traffic is still an intractable problem. Specifically, when the traffic becomes over-saturated, vehicle queues will form and overflow the road capacity, and then block the upstream intersection and even cause gridlocks; in this situations, green time is no longer effective as no traffic can be discharged, i.e., descending to de facto red (Abu-Lebdeh and Benekohal, 1997; Daganzo, 2007). To tackle the over-saturation problem, access control of the perimeter is required and should be combined with the coordinated control within the network, involving a more large-scale centralized control scheme. Chen et al. (2020b) proposed an innovative mixed-use BRT/AV lane model where the DBL is available for both BRT vehicles and automated private cars; to simultaneously avoid the interference from private cars on the operation of BRT vehicles and enhance the mobility of private cars, the access control and trajectory planning of private cars are incorporated; the results show that the proposed method can improve the overall efficiency of both the BRT/AV lane and regular lanes, under the premise of guaranteeing the operation of BRT system. The study of Chen et al. (2020b) illustrates the benefit of centralized organization of CAVs for the mobility of the heterogeneous traffic, even though it only focuses on a corridor with a relatively low market penetration ratio of CAVs. With the ongoing development of CAV technology and improvement of the customer acceptance of CAVs, more CAVs will be available on roads in the future, implying great potential to promote the overall mobility of a network. Therefore, a centralized organization and integrated control framework of heterogeneous traffic under a CAV environment is promising for enhancing both the transit service quality and the network mobility.\par

\begin{figure}[!ht]
	\centering
	\subfloat[][Virtual platoon]{\includegraphics[width=0.4\textwidth]{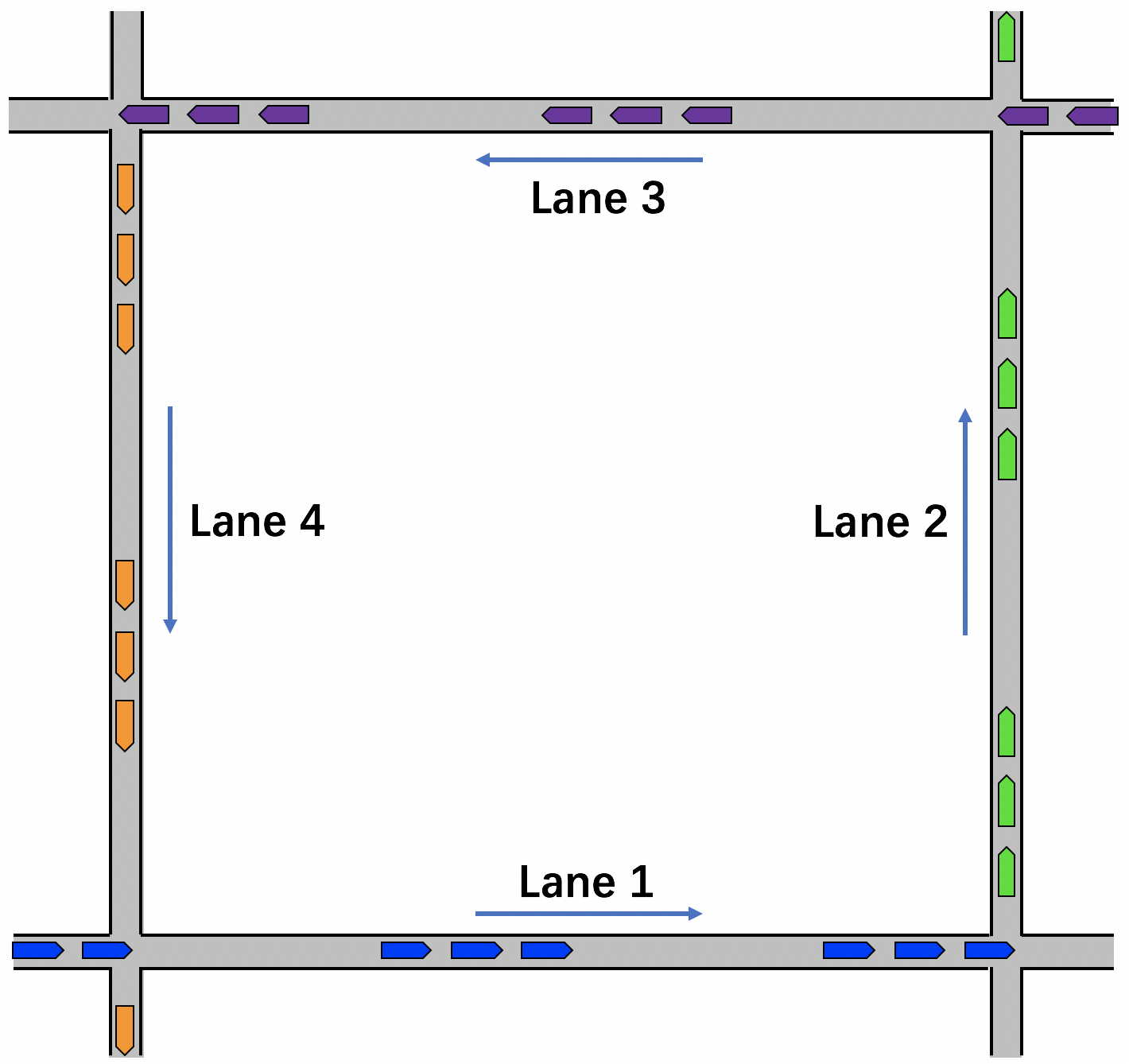}}\hspace{1em}
	\subfloat[][Collision-free trajectory]{\includegraphics[width=0.45\textwidth]{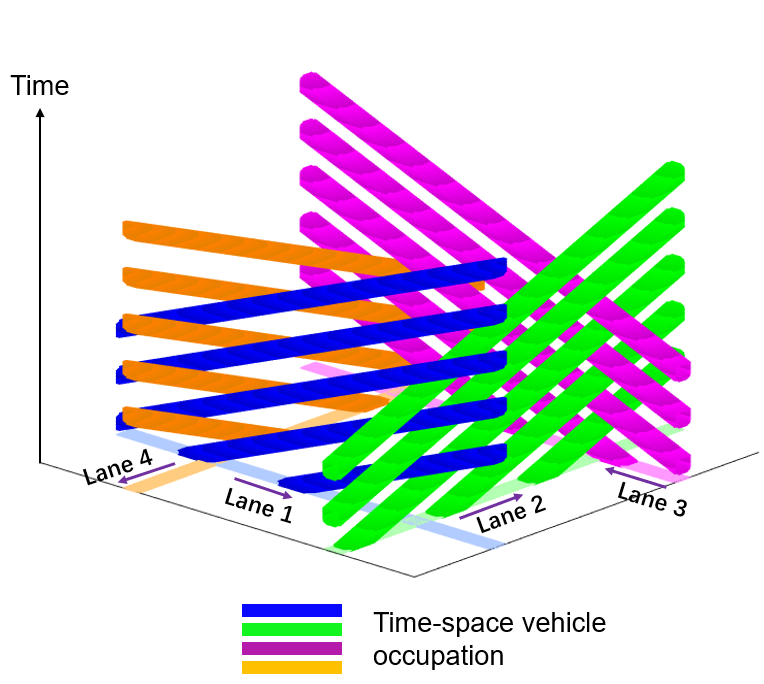}}\\	
	\caption[]{The concept of RC scheme (Lin et al., 2020)}
	\label{Fig_RCreview}
\end{figure}

Rhythmic control (RC), recently proposed by Chen et al. (2020a) and Lin et al. (2020), is a novel control framework to improve the automated traffic at intersections and on networks. The key idea of this approach is to organize the traffic into a regularly recurring and quickly-changing manner, which is analogous to the rhythm in music. To embody the rhythmic organization, the concept of "virtual platoon", abbreviated as VP, is introduced to represent a time-space slot that is available for being occupied by vehicles. Under the logic of RC, vehicles can follow the pace of VPs to pass through a traffic facility without any stops or collisions. As the name implies, the VPs just provide potential time-space resources for vehicles to travel with and the actual vehicle occupations are determined by the traffic flow. To resolve the complicated conflict relations at intersections, RC adopts a one-way network, which possesses the potential to improve traffic throughput and alleviate computational complexity in a real-time control implementation. Fig.\ref{Fig_RCreview} shows the concept of VPs and the collision-free trajectories.\par

Since RC is a highly centralized organization and control method, it is capable of meeting the requirement of the bus operation with regard to reliability. Moreover, as buses are generally operated in a cyclic manner, opportunely sharing the characteristic of rhythm with RC, there is great potential to incorporate bus travel into the RC scheme. However, it can be expected that the problem will become more complicated to design an RC scheme under the heterogeneous traffic scenario with both private cars and buses, due to the incompatibility of the rhythm for vehicles with distinctive travel features. Specifically, buses usually travel at a relatively lower speed than that of private cars, along with extra dwelling time at bus stations; thus, it is difficult for buses to catch up with the RC rhythm designed for the traffic of private cars, while accommodating private cars under the rhythm of buses will slow down the pace of private cars and incur enormous delays. \par

For the above research goal, in this study we propose an improved version of RC for heterogeneous CAVs, abbreviated as RC-H, to coherently organize and jointly control the traffic of both private cars and buses in a network. Inheriting the concept of “virtual platoons” proposed by Lin et al. (2020), “dedicated virtual platoons” are designed for buses, which are essentially space-time slots and provide exclusive ROWs for buses on their routes as priority treatment expected. Traveling within the dedicated VPs, buses will fulfill their trips with no collision and limited delays within the network, meanwhile on-time performance and schedule adherence can be guaranteed. In the meantime, “regular virtual platoons” are also designed for private cars, with whose paces the private cars can pass through the network without any stops and collisions. Both the dedicated VPs and regular VPs are unified as "realized VPs", which are extended from the notion of VP of the RC scheme and are regarded as the carriers of the background rhythm. Along with the design of realized VPs, the planning of bus itinerary, defined as a series of space-time trajectories along the travel path, and the traffic assignment of private cars are jointly optimized to acquire an optimal RC-H scheme. The optimization is aimed at minimizing the total personal travel cost and is formulated as a mixed-integer linear program (MILP). A bilevel solution method with a heuristic searching algorithm and tractable approximations is proposed to solve the MILP efficiently. To validate the performances of the RC-H scheme and the effectiveness of the heuristic algorithm, numerical experiments are conducted on both a toy example and a real-world network, while simulation experiments are carried out to compare the results of the RC-H scheme to those of the signal control strategies with and without DBLs.\par

For the remainder of the paper, Section \ref{RC-H_sec} first briefly reviews the RC scheme, and then introduces the concept of RC-H scheme. Section \ref{Optimal_design_sec} develops a joint optimization problem of RC-H scheme composed of realized VP design, bus itinerary planning and traffic assignment of private cars. A bilevel heuristic solution method is proposed to solve the optimization in Section \ref{Two_level_sec}. Section \ref{Numerical_sec} reports a series of numerical examples and simulation experiments for the performance validation of RC-H scheme, and a comparison between RC-H scheme and traffic signal control strategies. Finally, Section \ref{Concluding_sec} summarizes the paper and discusses some future research directions.

\section{Rhythmic control of heterogeneous automated traffic} \label{RC-H_sec}

This section introduces the concept of rhythmic control scheme of heterogeneous traffic of both buses and private cars, i.e., RC-H scheme. In the following, we first briefly review the RC scheme proposed in Lin et al. (2020), and then qualitatively describe how to coherently organize the traffic of both buses and private cars under the RC scheme.\par

\subsection{A brief review of RC} \label{RCreview_subsec}

Rhythmic control (RC) is a novel traffic control framework for vehicles on urban road networks. The basic idea of RC is to organize the whole network traffic into a regularly recurring and quickly-changing manner, where low delays and high traffic throughput can be achieved under the premise of collision avoidance. The concept of “virtual platoon”, i.e., VP, is proposed to represent a time-space slot that is available for being occupied by vehicles and is generated with a preset rhythm. Under the logic of RC, all of the vehicles are required to follow the paces of VPs along the vehicles’ paths to finish their trips and they will not encounter any stops or collisions after entering the network. 

Furthermore, Lin et al. (2021) proposes a rhythm design method to optimally determine the RC scheme composed of the cycle length, the VPs' sizes and their arrival times at conflicting nodes (i.e., the pace). The cycle length is the unified rhythm length of the RC scheme, which holds a trade-off between the network capacity and vehicle delay, as a longer rhythm can accommodate higher traffic demands but incurs long average waiting times at entrances or junctions. The sizes of VPs embody the traffic capacities of the links, that is, the maximum numbers of vehicles occupying on the platoons, violating which could potentially result in collisions at conflict points. As discussed above, the platoon sizes are determined according to the traffic demand pattern. Relative arrival time is proposed to denote the arrival time of VPs at the conflicting nodes in each cycle, based on which the travel time on links can be determined and the collision avoidance condition can be guaranteed. The design problem is formulated as a mixed-integer linear program (MILP), and by solving the MILP, an optimal RC scheme for minimizing the total delay can be obtained.\par

The readers can refer to Lin et al. (2020) and Lin et al. (2021) for detailed technical descriptions of the RC scheme. Here, we list some results that lay the foundation of our study of the RC-H scheme described below.
\begin{itemize}
    \item A VP represents a time-space slot that is available for being occupied by vehicles and is generated in a rhythmic manner.
    \item All vehicles are required to follow the VPs employed on the links and transfer between different platoons along the paths to finish their trips.
    \item Under the logic of RC, vehicles will not encounter any stops or collisions after entering the network.\par 
    \item The average vehicle delay through the entire network is almost negligible when the travel demand is not too high, and massive traffic demands can be well accommodated.
    \item RC provably maximizes network traffic throughput under some mild conditions.
\end{itemize}

\begin{figure}[!ht]
	\centering
	\subfloat[][]{\includegraphics[width=0.45\textwidth]{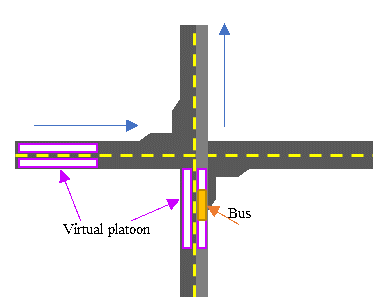}}\hspace{1em}
	\subfloat[][]{\includegraphics[width=0.45\textwidth]{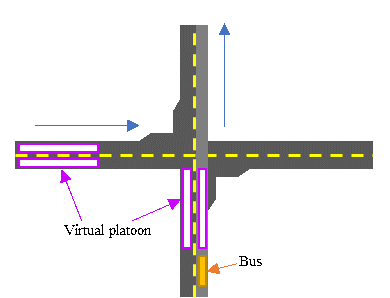}}\\
	\caption[]{Conflict relation of buses at an intersection}
	\label{Fig_Businters}
\end{figure}

\subsection{Concept of RC-H } \label{concept_RCH_subsec}
In a network, conflict points are the main crux of vehicle delays and traffic collisions, which should be given great focus in carrying out the centralized control of the network. Under the logic of the RC scheme, VPs are preset in each link, and all of the vehicles should travel with the pace of the VPs to avoid collisions at conflict points. The same way for buses, Fig.\ref{Fig_Businters} illustrates an example of the conflict relation of buses with VPs at a conflict point composed of two intersecting one-way roads. In Fig.\ref{Fig_Businters}(a), the bus catches up with a VP when passing through the intersection, while in Fig.\ref{Fig_Businters}(b), as the bus falls out of the VP, collision may occur between the bus and the intersecting vehicles at the intersection.\par

Unfortunately, as discussed above, due to their lower speed and dwelling process, it is difficult for buses to catch up with the rhythm designed for accommodating private cars. In the meantime, due to the interference from buses, which can be regarded as “moving bottlenecks”, the travel of private cars in the mixed traffic may also be blocked. Therefore, the original RC scheme is incompatible in the heterogeneous traffic, and an improved RC scheme is required for accommodating both private cars and buses in a harmonized rhythmic manner. Motivated by Fig.\ref{Fig_Businters}(a), in the heterogeneous traffic, we can still coordinate the arrivals of private cars and buses at intersections to form "common paces". It comes to an idea that by re-designing the VPs in the mixed traffic with some customized adjustments, the heterogeneous traffic of buses and private cars can be accommodated simultaneously under the RC framework. \par 

\begin{figure}[!ht]
	\centering
	{\includegraphics[width=0.7\textwidth]{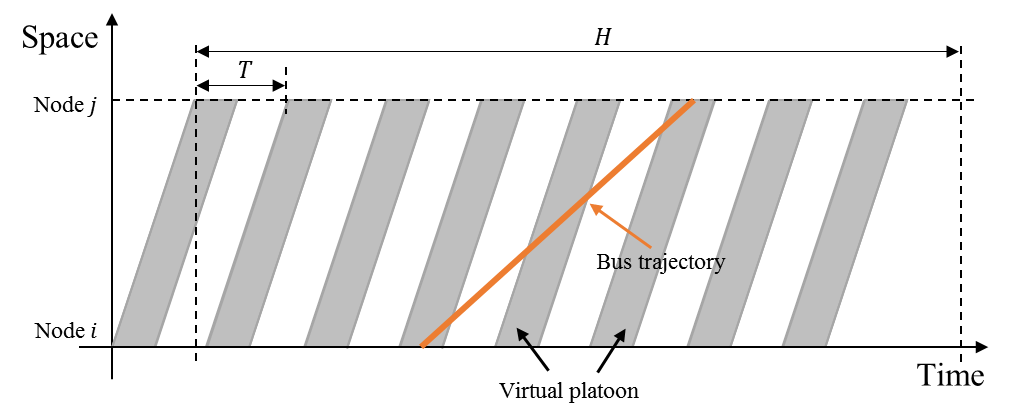}}
	\caption[]{Illustration of a valid bus trajectory}
	\label{Fig_RCH_bustrajetory}
\end{figure}

\begin{figure}[!ht]
	\centering
	{\includegraphics[width=0.7\textwidth]{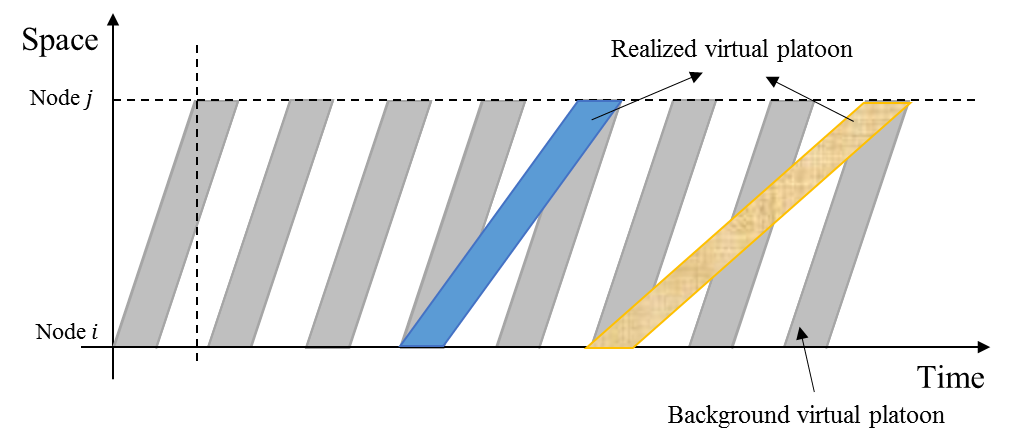}}
	\caption[]{Trajectories for different types of virtual platoons}
	\label{Fig_RCH_concept}
\end{figure}

Specifically, taking an RC scheme obtained from the study of  Lin et al. (2020)  as the background, which is designed for private cars with the preset rhythm, platoon size and relative node arrival time, we draw the trajectories of VPs in Fig.\ref{Fig_RCH_bustrajetory}, where the grey bands denote the VPs on link $(i,j)$. For a vehicle intended to pass through the link, as long as it could catch up with some VP at the entrance (node $i$) and the exit (node $j$) of the link, collision avoidance can be guaranteed; Fig.\ref{Fig_RCH_bustrajetory} shows an example of a valid bus trajectory. Then, by combining different VPs at the entrance and exit, a set of new “distorted” VPs can be generated, as the blue and yellow bands in Fig.\ref{Fig_RCH_concept} shows. To avoid confusion, we define the newly generated VP as “realized VP”, and the preset VP as “background VP”. The implication is that the realized VP will be assigned with real traffic while the background VP only provides the background rhythm without actual traffic assignment.\par
From Fig.\ref{Fig_RCH_concept}, the relationship between the realized VP and the background VP can be described as follows:
\begin{itemize}
    \item The background VPs are predetermined and provide a fundamental rhythm for the realized VPs, while the realized VPs can be generated as different for the same link in order to accommodate vehicles with different speeds.
    \item Under the mixed traffic, the realized VP can be regarded as a “stretched” background VP as the travel time of the former should be no smaller than that of the latter in the presence of buses. However, there is still a chance that the realized VP coincides with the background platoon, such as when the interference of buses can be dodged for private cars.
    \item Under the homogeneous traffic of private cars, the realized VPs overlap with the background platoons, as the background VPs are designed for the traffic of private cars, and they can travel freely without the interference from buses.
\end{itemize}\par

Moreover, in general, the system of buses operates in a cyclic manner with a fixed cycle denoted as $H$, while the background VP is generated with a cycle of $T$ under the RC scheme as shown in Fig.\ref{Fig_RCH_bustrajetory}. Therefore, the proposed RC-H scheme also possesses periodicity with a cycle length of the lowest common multiple of $T$ and $H$, denoted as $C$, which characterizes the rhythm of the realized VPs. To simplify the problem, we assume that $H=QT,Q\in\mathbb{N}^+$, where $Q$ is the number of the RC cycles, equivalently the number of background VPs, for a single RC-H cycle\footnote{In fact, $H=QT$ can be easily realized as on the one hand, $T$ is always far smaller than $H$ ($T$ is no greater than 20 seconds while $H$ is generally on the order of hundreds of seconds) and the setting of $T$ can be adjusted as expected; on the other hand, the cycle of the bus system can be regarded as $C$, where the repeated buses for the same line are generated.}. To differentiate the two cycles, "RC cycle" refers to the cycle of the RC scheme (i.e.,$T$) while "RC-H cycle" refers to that of the RC-H scheme (i.e.,$H$). As the movements of vehicles are the same from one RC-H cycle to another, we can model the movements of vehicles in a single RC-H cycle only.  \par

To provide exclusive ROWs of buses and avoid the interference from private cars, the VPs assigned for buses are separated from those for private cars, i.e., "dedicated VPs" for buses and "regular VPs" for private cars, both of which are realized VPs. With the introduction of realized VPs, the heterogeneous traffic of buses and private cars can be coherently organized into a harmonized rhythmic manner. Along with the design of realized VPs, the itinerary planning of buses and traffic assignment of private cars can be further jointly considered, which is expected to achieve the optimal design of the RC-H scheme, as discussed in the next section. 

\section{Optimal design of RC-H} \label{Optimal_design_sec}
The RC-H scheme proposed in Section \ref{RC-H_sec} guarantees that by following the realized VP, vehicles would be collision free at all conflict points in the network, but may incur delays due to the "stretched deformation" of the background VPs. Therefore, in this section, we present a framework to minimize the total travel cost by designing the RC-H parameters, when given bus line and traffic demand of private cars information.\par 

Since previous studies have addressed the design problem of RC scheme, the related parameters including the cycle length, platoon sizes and relative node arrival times of the RC scheme are not decision variables and are already provided prior to solving the RC-H design problem. Therefore, the RC-H design can be modeled as a combined optimization problem with three components; (i) realized VP design is carried out to determine the dedicated VPs for buses and the regular VPs for private cars; (ii) based on the dedicated platoons, the itinerary of the buses can be planned given the bus routes and time schedule; (iii) meanwhile, the itinerary of private cars can also be determined, along with the traffic assignment on the different paths and different types of lanes. Then, by jointly optimizing the three components, the operation of buses and the overall traffic efficiency can be enhanced concurrently.\par

\subsection{Basic settings} \label{basic_RC-H_subsec}
Consider a network composed of links with multiple lanes, where buses are only allowed on a bus lane although the bus lane is shared with other traffic \footnote{Note that here we reserve only one lane for buses as the number of buses is generally limited compared to the private cars, and this approach can be easily generalized to the scenario where multiple bus lanes are required.}. Note that the lane served for buses is not a dedicated lane, but some time-space slots on the lane are dedicated, following the concept of dedicated VPs. Thus, two types of lanes are defined as mixed-traffic lanes and normal lanes. The network can be represented by a directed graph, i.e., $\mathcal{G}=(\mathcal{A}, \mathcal{V})$, where $\mathcal{A}$ denotes the link set and $\mathcal{V}$ denotes the node set. We can also use two nodes to denote a link, i.e., $(i,j)\in\mathcal{A}$, where $i,j\in\mathcal{V}$. The node set can be further divided into origin node set $\mathcal{V}_o$, destination node set $\mathcal{V}_d$, bus station node set $\mathcal{V}_s$ and intersection node set $\mathcal{V}_{in}$.
The notations used for this study are listed in Appendix A.

\subsection{Realized virtual platoon design} \label{virtualpla_subsec}

\begin{figure}[!ht]
	\centering
	{\includegraphics[width=0.7\textwidth]{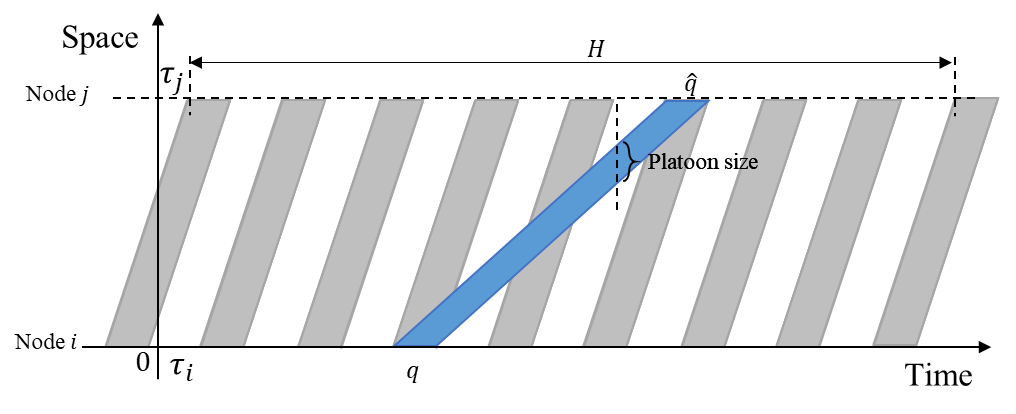}}\\	
	\caption[]{Trajectory of a realized virtual platoon}
	\label{Fig_realized_virtualpla}
\end{figure}

We introduce a binary variable $\Theta_{q,\hat{q}}^{ij}$ to denote the VP on link $(i,j)$ arriving with platoon $q$ and leaving with platoon $\hat{q}$, where $q,\hat{q}\in\mathcal{Q}$, defined as $q\rightarrow\hat{q}$. Fig.\ref{Fig_realized_virtualpla} illustrates the trajectory of the realized VP. If the VP is realized, $\Theta_{q,\hat{q}}^{ij}=1$, otherwise, $\Theta_{q,\hat{q}}^{ij}=0$. Note that the size of the realized VP is identical to the background platoon and any arrival platoon $q$ or departure platoon $\hat{q}$ can be occupied for only one realized platoon, constrained by the lane capacity, that is:
\begin{align}
    & \sum_{\hat{q}\in\mathcal{Q}}\Theta_{q,\hat{q}}^{ij}\leq 1 & \forall ij\in\mathcal{A},q\in\mathcal{Q}\\
    & \sum_{q\in\mathcal{Q}}\Theta_{q,\hat{q}}^{ij}\leq 1 & \forall ij\in\mathcal{A},\hat{q}\in\mathcal{Q}    
\end{align}\par

\vspace{0.5em}
\noindent\underline{\textit{Travel time computation}}\par
Denote the travel time of the realized VP on link $(i,j)$ of $q\rightarrow\hat{q}$ as $t_{q,\hat{q}}^{ij}$. The travel time can be calculated as:
\begin{align}
    \begin{split}
        & t_{q,\hat{q}}^{ij}=\tau_j-\tau_i+(\hat{q}-q)T+\beta_{q,\hat{q}}^{ij}H\\
        &\beta_{q,\hat{q}}=\min\{\beta\in\mathbb{N}|\tau_j-\tau_i+(\hat{q}-q)T+\beta H\geq t_a^{ij}\}
    \end{split} & \forall ij\in\mathcal{A}, q,\hat{q}\in\mathcal{Q}
\end{align}
where $\tau_i$ is the relative arrival time of the VPs at node $i$ as shown in Fig.\ref{Fig_realized_virtualpla}, and $t_a^{ij}$ is the travel time of the background VPs on link $(i,j)$. Therefore, the travel time of realized VPs either has several additional cycles of $T$ or is identical to that of the background platoon. Note that the travel times of the realized VPs are not variables and are already computed prior to solving the problem.\par

\vspace{0.5em}
\noindent\underline{\textit{FIFO constraints}}\par
Under the RC scheme, as the vehicles are organized to travel in a relatively tight manner to promote the network throughput, overtaking is not allowed in the mainline roads to reduce traffic risks. Thus, in the realization of the VP, the first-in-first-out (FIFO) constraint should be satisfied. Considering two VPs of $q_m\rightarrow\hat{q}_m$ and $q_n\rightarrow\hat{q}_n$, the FIFO constraint can be expressed as:
\begin{align}
    & \Theta_{q_m,\hat{q}_m}^{ij}+\Theta_{q_n,\hat{q}_n}^{ij}\leq \zeta_{q_m,\hat{q}_m,q_n,\hat{q}_n}^{ij} & \forall ij\in\mathcal{A}, q_m,\hat{q}_m,q_n,\hat{q}_n\in\mathcal{Q}
\end{align}
where $\zeta_{q_m,\hat{q}_m,q_n,\hat{q}_n}^{ij}$ denotes whether or not the FIFO constraint between the two realized platoons is violated; if such a violation occurs, $\zeta_{q_m,\hat{q}_m,q_n,\hat{q}_n}^{ij}=1$; otherwise, $\zeta_{q_m,\hat{q}_m,q_n,\hat{q}_n}^{ij}=2$. 
For determining the value of $\zeta_{q_m,\hat{q}_m,q_n,\hat{q}_n}^{ij}$, without loss of generality, we assume $q_m<q_n$ to obtain two cases differentiated by $\hat{q}_m>\hat{q}_n$ or $\hat{q}_m<\hat{q}_n$.\par

\begin{figure}[!ht]
	\centering
	\subfloat[][]{\includegraphics[width=0.45\textwidth]{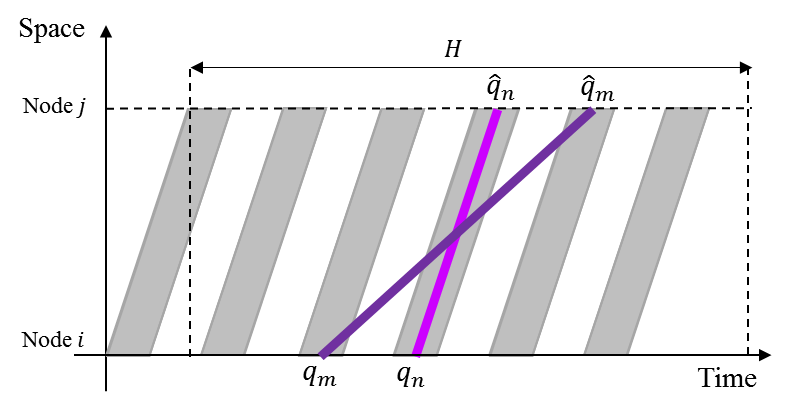}}\hspace{1em}
	\subfloat[][]{\includegraphics[width=0.45\textwidth]{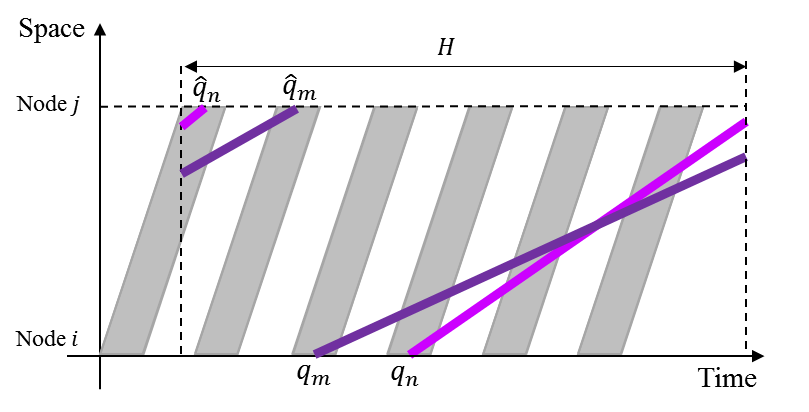}}\\
	\subfloat[][]{\includegraphics[width=0.45\textwidth]{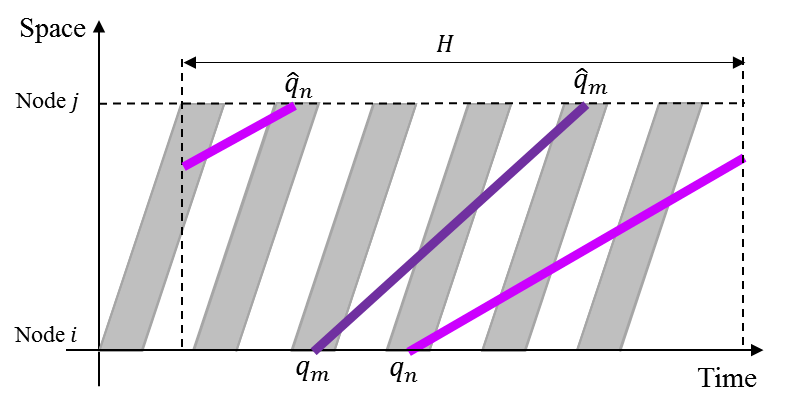}}\\
	\caption[]{Case I of FIFO constraint}
	\label{Fig_FIFO_case1}
\end{figure}

Case I: $q_m<q_n$ and $\hat{q}_m>\hat{q}_n$. Note that it is possible that $\hat{q}_m$ and $q_m$ are in different cycles, and it can be evaluated by comparing $\hat{q}_m$ and $q_m+\alpha_{ij}$, where $\alpha_{ij}$ denotes the number of cycle $T$ required by the background platoons to pass through link $(i,j)$. Therefore, three possible relations between the platoons of $q_m\rightarrow\hat{q}_m$ and $q_n\rightarrow\hat{q}_n$ exist, as shown in Fig.\ref{Fig_FIFO_case1}. It is observed that the FIFO constraint is violated in the first two relations, i.e., $\hat{q}_m\geq q_m+\alpha_{ij}, \hat{q}_n\geq q_n+\alpha_{ij}$ as shown in Fig.\ref{Fig_FIFO_case1}(a) and
$\hat{q}_m<q_m+\alpha_{ij}, \hat{q}_n<q_n+\alpha_{ij}$ as shown in Fig.\ref{Fig_FIFO_case1}(b), where $\zeta_{q_m,\hat{q}_m,q_n,\hat{q}_n}^{ij}$ should be set to $1$.\par
\begin{figure}[!ht]
	\centering
	\subfloat[][]{\includegraphics[width=0.45\textwidth]{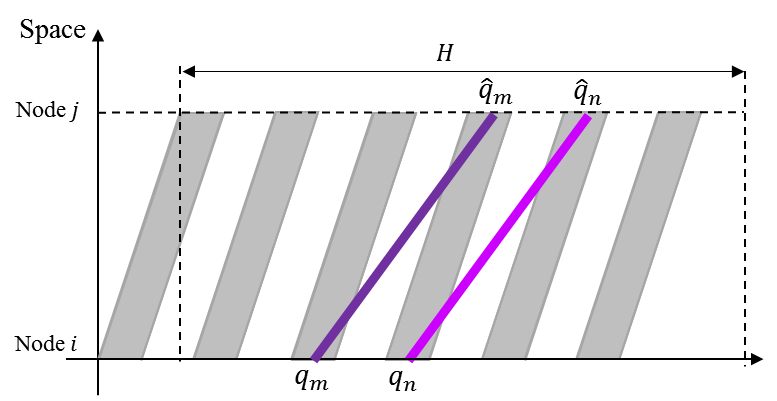}}\hspace{1em}
	\subfloat[][]{\includegraphics[width=0.45\textwidth]{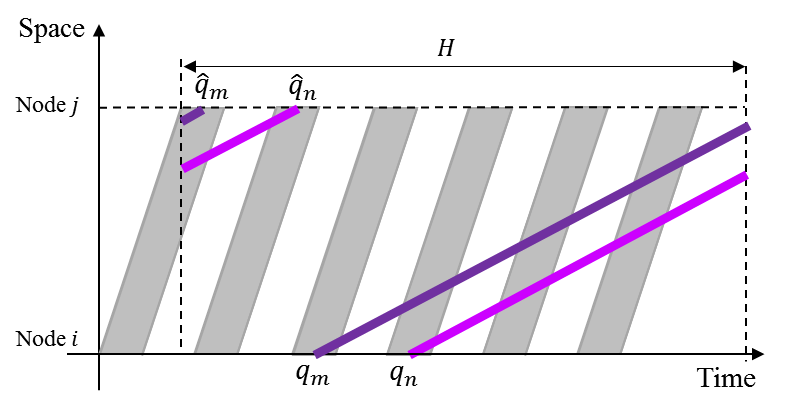}}\\
	\subfloat[][]{\includegraphics[width=0.45\textwidth]{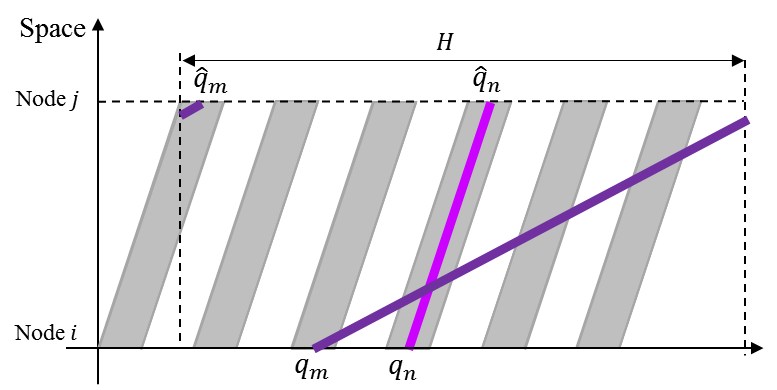}}\hspace{1em}
	\subfloat[][]{\includegraphics[width=0.45\textwidth]{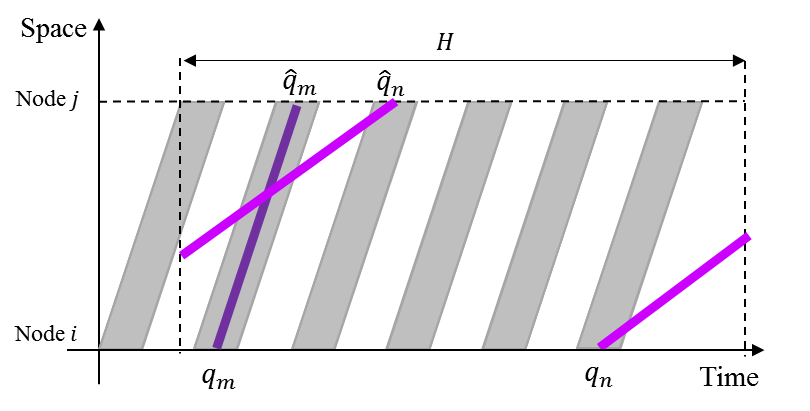}}\\
	\caption[]{Case II of FIFO constraint}
	\label{Fig_FIFO_case2}
\end{figure}

Case II: $q_m<q_n$ and $\hat{q}_m<\hat{q}_n$. For this case, four possible relations between the platoons of $q_m\rightarrow\hat{q}_m$ and $q_n\rightarrow\hat{q}_n$ exist, as shown in Fig.\ref{Fig_FIFO_case2}. The FIFO constraint is violated in the last two relations, i.e., $\hat{q}_m\geq q_m+\alpha_{ij}, \hat{q}_n< q_n+\alpha_{ij}$ as shown in Fig.\ref{Fig_FIFO_case2}(c) and $\hat{q}_m<q_m+\alpha_{ij}, \hat{q}_n\geq q_n+\alpha_{ij}$ as shown in Fig.\ref{Fig_FIFO_case2}(d).

\subsection{Bus itinerary planning} \label{bus_itine_subsec}
As buses have distinct travel characteristics, such as fixed schedule, preset routes and dwelling processes, the traveling of buses should be modeled separately from that of private cars. Moreover, to provide exclusive ROWs for buses along their trips, dedicated VPs are introduced for buses to follow up, differentiating from that for private cars, i.e., regular VPs.\par 

\vspace{0.5em}
\noindent\underline{\textit{Dedicated virtual platoon}}\par
\begin{figure}[!ht]
	\centering
	{\includegraphics[width=0.7\textwidth]{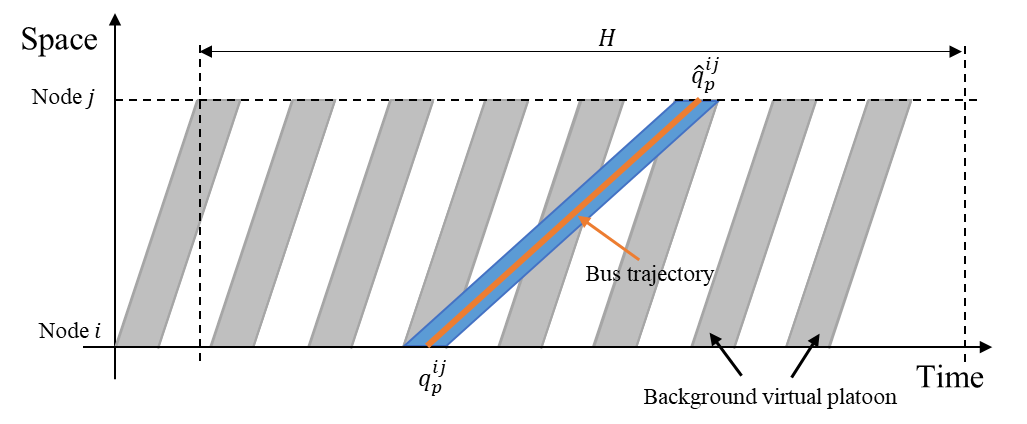}}
	\caption[]{Illustration of a dedicated virtual platoon}
	\label{Fig_dedicated_virtualpla}
\end{figure}

Assume a bus traveling on link $(i,j)$ as shown in Fig.\ref{Fig_dedicated_virtualpla}. We introduce another binary variable $\hat{\Theta}_{q,\hat{q}}^{p,ij}$ to denote the dedicated VP for bus $p$, and $p\in\mathcal{P}_{ij}$ and $\mathcal{P}_{ij}$ denotes the bus set traveling on link $(i,j)$. $\hat{\Theta}_{q,\hat{q}}^{p,ij}=1$ denotes that the VP of $q\rightarrow\hat{q}$ is dedicated for bus $p$ as shown by the blue band in Fig.\ref{Fig_dedicated_virtualpla}; $\hat{\Theta}_{q,\hat{q}}^{p,ij}=0$ otherwise. For each bus, one and only one dedicated VP is required to load the bus on each link, that is: 
\begin{align}
    & \sum_{\hat{q}\in\mathcal{Q}}\sum_{q\in\mathcal{Q}}\hat{\Theta}_{q,\hat{q}}^{p,ij}=1, & \forall ij\in\mathcal{A}, p\in\mathcal{P}_{ij}
\end{align}\par 

To establish the relationship between the dedicated VP and the realized VP, we introduce $\hat{\Theta}_{q,\hat{q}}^{ij}$ to denote whether or not platoon $q\rightarrow\hat{q}$ is dedicated for buses. Note that one dedicated VP can load multiple buses, as long as the platoon capacity constraint is satisfied. Thus, as long as one bus occupies the platoon $q\rightarrow\hat{q}$, $\hat{\Theta}_{q,\hat{q}}^{ij}=1$; otherwise, $\hat{\Theta}_{q,\hat{q}}^{ij}=0$. Then, we have:
\begin{align}
    & M\times\hat{\Theta}_{q,\hat{q}}^{ij}\geq \sum_{p\in\mathcal{P}_{ij}}\hat{\Theta}_{q,\hat{q}}^{p,ij}, & \forall ij\in\mathcal{A}, q,\hat{q}\in\mathcal{Q} \label{bus_size1}\\
    & \hat{\Theta}_{q,\hat{q}}^{ij}\leq \sum_{p\in\mathcal{P}_{ij}}\hat{\Theta}_{q,\hat{q}}^{p,ij}, & \forall ij\in\mathcal{A}, q,\hat{q}\in\mathcal{Q}
\end{align}
 and the constraint of platoon capacity can be expressed as:
\begin{align}
    & s_b\sum_{p\in\mathcal{P}_{ij}}\hat{\Theta}_{q,\hat{q}}^{p,ij}\leq s_a^{ij}, & \forall ij\in\mathcal{A}, q,\hat{q}\in\mathcal{Q}\label{bus_size2} 
\end{align}
where $M$ is a large number, and $s_b$ and $s_a^{ij}$ denote the size of buses and VPs on link $(i,j)$, respectively. Incorporating (\ref{bus_size1}) and (\ref{bus_size2}), we further yield:
\begin{align}
    & s_a^{ij}\times\hat{\Theta}_{q,\hat{q}}^{ij}\geq s_b \sum_{p\in\mathcal{P}_{ij}}\hat{\Theta}_{q,\hat{q}}^{p,ij}, & \forall ij\in\mathcal{A}, q,\hat{q}\in\mathcal{Q}
\end{align}\par
Moreover, the dedicated VPs are feasible only under the premise that it has been realized, which can be expressed as:
\begin{align}
    & \hat{\Theta}_{q,\hat{q}}^{ij}\leq\Theta_{q,\hat{q}}^{ij}, & \forall ij\in\mathcal{A}, q,\hat{q}\in\mathcal{Q}
\end{align}\par

\begin{figure}[!ht]
	\centering
	\subfloat[][Station in the mainline]{\includegraphics[width=0.7\textwidth]{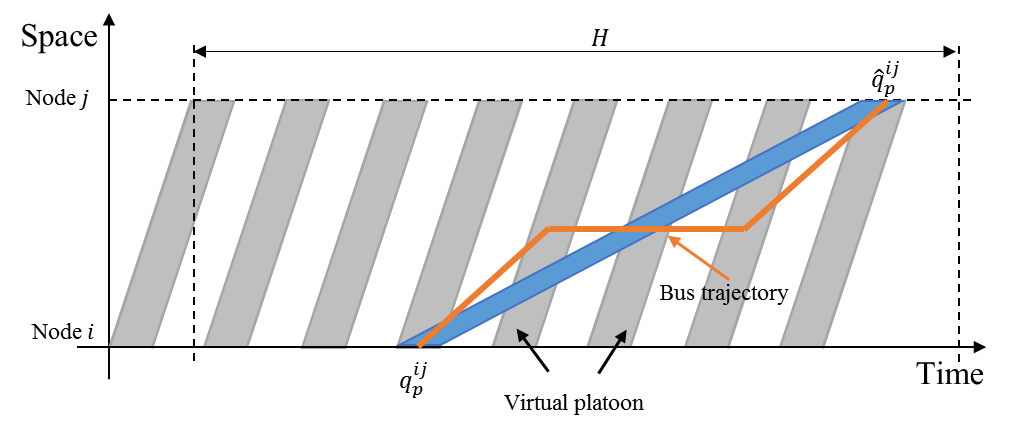}}\\
	\subfloat[][Station at side-platform]{\includegraphics[width=0.7\textwidth]{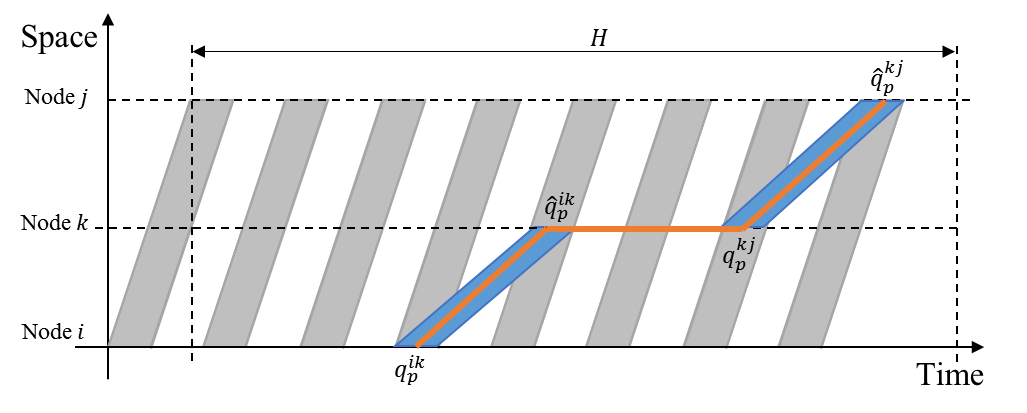}}\\
	\caption[]{Illustration of dedicated virtual platoons at stations}
	\label{Fig_dedicated_dwell}
\end{figure}

For describing the dwelling process of buses at stations, two primary types of bus stations are considered; the first is placed in the mainline road, where the vehicles following the bus will be blocked in the dwelling of the bus; the second is placed as a side-platform, where the main road will be free for the following vehicles when the bus dwells at the station, that is, overtaking is allowed during the time of bus dwelling. For the first type of station, the dwelling time of buses can be embodied by adding to the traveling time, as shown in Fig.\ref{Fig_dedicated_dwell}(a); this is equivalent to a bus traveling at a much slower speed, as the impact on the traffic remains the same. For the second type of station, as the bus in the dwelling seems essentially “disappeared” and will not block the traffic on the main road, the VPs planned to pass through within the dwelling time are available. Thus, to describe the movements on the lane with a side-platform station, we introduce an additional node at the position of the station that divides the original link into two connecting links, as shown in Fig.\ref{Fig_dedicated_dwell}(b); then, the traveling of the bus will be modeled separately in the two links, and the VPs passing within the dwelling time can be released for loading other vehicles. \par

\begin{figure}[!ht]
	\centering
	{\includegraphics[width=0.7\textwidth]{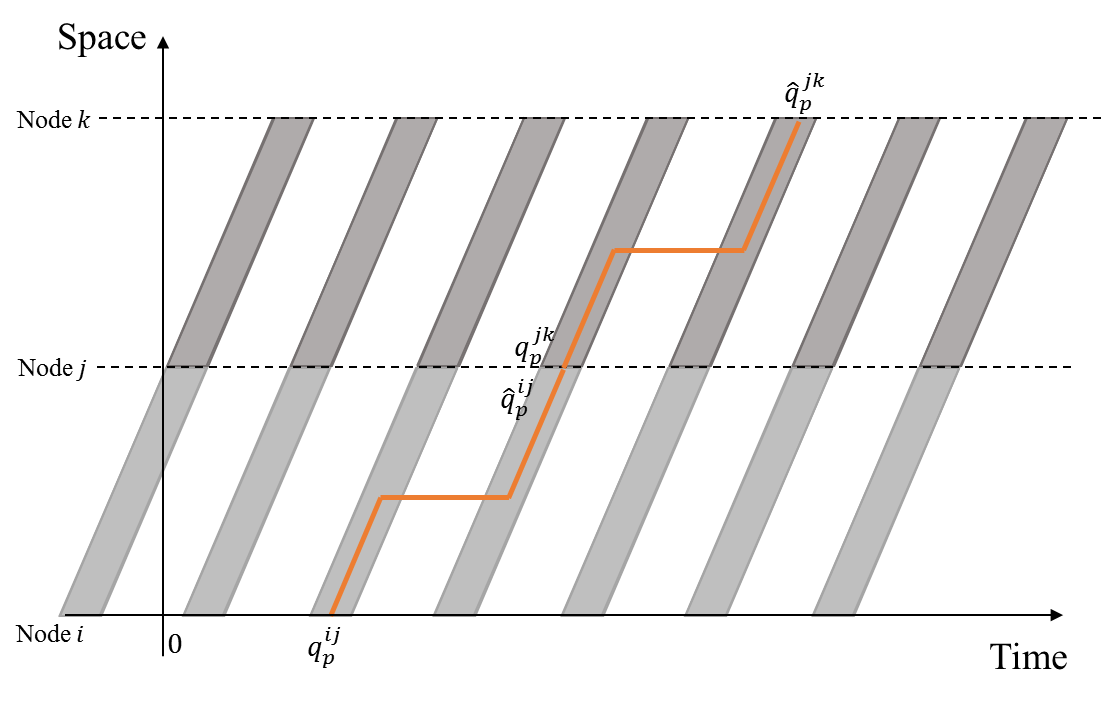}}\\
	\caption[]{Illustration for link connectivity}
	\label{Fig_dedicated_connect}
\end{figure}

\vspace{0.5em}
\noindent\underline{\textit{Link connectivity}}\par
Under the logic of the RC-H scheme, buses will pass through the network continuously without stops (except to dwell at stations), as shown in Fig.\ref{Fig_dedicated_connect}; that is, when the bus departs from link $(i,j)$, it will enter the next link in its route, i.e., link $(j,k)$, immediately. The number of the arrival platoon and the departure platoon of bus $p$ can be obtained as: 
\begin{align}
    & q_p^{ij}=\sum_{\hat{q}\in\mathcal{Q}}\sum_{q\in\mathcal{Q}}q\times\hat{\Theta}_{q,\hat{q}}^{p,ij}, & \forall ij\in\mathcal{A}, p\in\mathcal{P}_{ij}\\
    & \hat{q}_p^{ij}=\sum_{\hat{q}\in\mathcal{Q}}\sum_{q\in\mathcal{Q}}\hat{q}\times\hat{\Theta}_{q,\hat{q}}^{p,ij}, & \forall ij\in\mathcal{A}, p\in\mathcal{P}_{ij}
\end{align} 
and the link connectivity constraint can be expressed as:
\begin{align}
    & q_p^{jk}-\hat{q}_p^{ij}=0, & \forall j\in\mathcal{V}_{in}, p\in\mathcal{P}_{ij}\cup\mathcal{P}_{jk}
\end{align}
where $\mathcal{V}_{in}$ denotes the intersection node set. \par

For a station node, as shown in Fig.\ref{Fig_dedicated_dwell}(c), the gap between the departure platoon $\hat{q}_p^{ij}$ and the arrival platoon $q_p^{jk}$ should reflect the dwelling time of the bus, that is:
\begin{align}
    & \hat{t}_{p}^j=(q_p^{jk}-\hat{q}_p^{ij}+\hat{\phi}_{p}^jQ)\times T, & \forall j\in\mathcal{V}_s, p\in\mathcal{P}_{ij}\cup\mathcal{P}_{jk}\\
    & \underline{\hat{t}}_{p}^j\leq \hat{t}_{p}^j\leq H-\epsilon, & \forall j\in\mathcal{V}_s, p\in\mathcal{P}_{ij}\cup\mathcal{P}_{jk}
\end{align}
where $\mathcal{V}_s$ denotes the set of station nodes and $\epsilon$ is a small number used to eliminate the strict inequality of $\hat{t}_{p}^j<H$ in the formulation; $\hat{t}_{p}^j$ denotes the dwelling time of bus $p$ at station $j$, and $\underline{\hat{t}}_{p}^j$ denotes the lower bound of the dwelling time; $\hat{\phi}_{p}^j\in\mathbb{N}$.

\vspace{0.5em}
\noindent\underline{\textit{Travel time computation}}\par
Similar to the computation for the travel time of realized VPs, the bus travel time can be calculated as:
\begin{align}
    & t_p^{ij}=\tau_j-\tau_i+(\hat{q}_p^{ij}-q_p^{ij})T+\phi_p^{ij}H, & \forall ij\in\mathcal{A}, p\in\mathcal{P}_{ij}\\
    & \underline{t}_p^{ij}\leq t_p^{ij}\leq H+t_{a}^{ij}-\epsilon, & \forall ij\in\mathcal{A}, p\in\mathcal{P}_{ij}
\end{align}
where $t_p^{ij}$ denotes the travel time of bus $p$ on link $(i,j)$, which is a variable as the arrival platoon and departure platoon are decision variables, and $\underline{t}_p^{ij}$ denotes the lower bound of the travel time; $\phi_p^{ij}\in\mathbb{N}$. The inequality $t_p^{ij}<H+t_{a}^{ij}$ implies that the travel time difference between buses and private cars on the link is below $H$ in order to maintain the FIFO constraint. This is reasonable as a dedicated VP with a too large travel time will not be utilized due to the objective of minimizing the travel delay. Furthermore, for the link with a long length where the travel time difference is expected to exceed $H$, we can double the system cycle $H$ or divide the link into two connecting shorter links to handle this issue.\par

\subsection{Traffic assignment of private cars} \label{traffic_assign_subsec}
\begin{figure}[!ht]
	\centering
	{\includegraphics[width=0.7\textwidth]{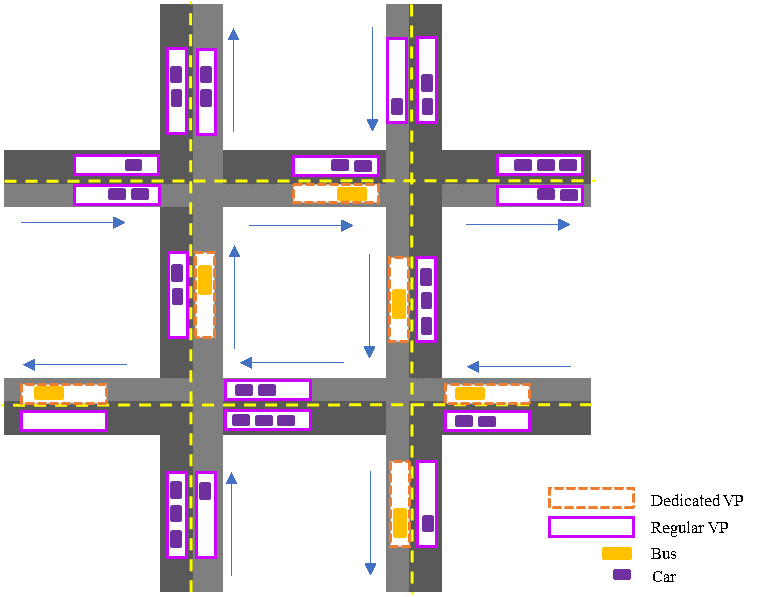}}\\
	\caption[]{Illustration of traffic assignment of private cars}
	\label{Fig_traffic assignment}
\end{figure}

Unlike for buses, all of the lanes are available for private cars. In regular lanes, private cars can follow the preset (background) VPs with no extra delay, while in mixed-traffic lanes, they are required to follow the regular VPs that is realized but not dedicated for the buses. Therefore, the traffic assignment will determine which path and which lane on each link will be used, along with which regular VP on each lane to follow for the traffic of private cars, as Fig.\ref{Fig_traffic assignment} shows.\par

\vspace{0.5em}
\noindent\underline{\textit{Feasible path}}\par
First, the traffic of private cars of different O-D pairs should be assigned to the corresponding feasible paths. Assume that $d^w$ is the stationary demand for O-D pair $w$, and $f^r$ is the traffic flow on path $r$, we have:
\begin{align}
    & \sum_{r\in\mathcal{R}^w}f^r=d^w, & \forall w\in\mathcal{W}
\end{align}
where $\mathcal{R}^w$ is the path set for O-D pair $w$. Denote the set of all feasible paths by $\mathcal{R}$, i.e., $\mathcal{R}=\cup_{w\in\mathcal{W}}\mathcal{R}^w$.\par

\vspace{0.5em}
\noindent\underline{\textit{Traffic assignment}}\par
For the traffic of private cars in the mixed-traffic lane, we introduce variables $\pi_{q,\hat{q}}^{r,ij}$, $\lambda_{q}^{r,ij}$ and $\mu_{\hat{q}}^{r,ij}$ to denote the traffic flow following the VP $q\rightarrow\hat{q}$, the traffic arriving with platoon $q$ and the traffic leaving with platoon $\hat{q}$ on link $(i,j)$ belonging to path $r$, respectively. Thus, we have:
\begin{align}
    & \sum_{q\in\mathcal{Q}}\pi_{q,\hat{q}}^{r,ij}=\mu_{\hat{q}}^{r,ij}, & \forall r\in\mathcal{R}, ij\in\mathcal{A}, \hat{q}\in\mathcal{Q}\\
    & \sum_{\hat{q}\in\mathcal{Q}}\pi_{q,\hat{q}}^{r,ij}=\lambda_{q}^{r,ij}, & \forall r\in\mathcal{R}, ij\in\mathcal{A}, q\in\mathcal{Q}    
\end{align}\par
Note that in the mixed-traffic lane, the traffic can be only assigned to the regular VPs, which have been realized but not dedicated for buses, that is:
\begin{align*}
    \sum_{r\in\mathcal{R}}{\pi_{q,\hat{q}}^{r,ij}}\leq M\times (\Theta_{q,\hat{q}}^{ij}-\hat{\Theta}_{q,\hat{q}}^{ij})
\end{align*}
implying that only when $\Theta_{q,\hat{q}}^{ij}=1$ and $\hat{\Theta}_{q,\hat{q}}^{ij}=0$, there is a chance for $\sum_{r\in\mathcal{R}}{\pi_{q,\hat{q}}^{r,ij}}>0$. Additionally, as the traffic should be constrained by the platoon size, the constraint can be further expressed as:
\begin{align}
    & \sum_{r\in\mathcal{R}}{\pi_{q,\hat{q}}^{r,ij}}\leq s_a^{ij}\times (\Theta_{q,\hat{q}}^{ij}-\hat{\Theta}_{q,\hat{q}}^{ij}), & \forall ij\in\mathcal{A}, q, \hat{q}\in\mathcal{Q}  
\end{align}\par 
For the regular lanes, where private cars will follow the background VPs, we introduce $\tilde{\Theta}_{q,\hat{q}}^{ij}$ to denote the background VPs; if platoon $q\rightarrow\hat{q}$ is on the rhythm of the background RC, $\tilde{\Theta}_{q,\hat{q}}^{ij}=1$; otherwise, $\tilde{\Theta}_{q,\hat{q}}^{ij}=0$. Note that $\tilde{\Theta}_{q,\hat{q}}^{ij}$ is not a variable and is determined by the background RC scheme. In the same way, we introduce variables $\tilde\pi_{q,\hat{q}}^{r,ij}$, $\tilde{\lambda}_{q}^{r,ij}$ and $\tilde{\mu}_{\hat{q}}^{r,ij}$ with the similar definition for the traffic assigned to the mixed-traffic lanes. Thus, we have:
\begin{align}
    & \sum_{q\in\mathcal{Q}}\tilde{\pi}_{q,\hat{q}}^{r,ij}=\tilde{\mu}_{\hat{q}}^{r,ij}, & \forall r\in\mathcal{R}, ij\in\mathcal{A}, \hat{q}\in\mathcal{Q}\\
    & \sum_{\hat{q}\in\mathcal{Q}}\tilde{\pi}_{q,\hat{q}}^{r,ij}=\tilde{\lambda}_{q}^{r,ij}, & \forall r\in\mathcal{R}, ij\in\mathcal{A}, q\in\mathcal{Q}\\
    & \sum_{r\in\mathcal{R}}{\tilde{\pi}_{q,\hat{q}}^{r,ij}}\leq s_a^{ij}\times(l_{ij}-1)\times\tilde\Theta_{q,\hat{q}}^{ij}, &\forall ij\in\mathcal{A}, q, \hat{q}\in\mathcal{Q} 
\end{align}
where $l_{ij}$ is the lane number on link $(i,j)$.\par

\vspace{0.5em}
\noindent\underline{\textit{Node conservation}}\par
For the node that is neither an entrance nor an exit, the conservation constraint is:
\begin{align}
    & (\lambda_{q}^{r,jk}+\tilde\lambda_{q}^{r,jk})-(\mu_{q}^{r,ij}+\tilde\mu_{q}^{r,ij})=0, & \forall r\in\mathcal{R}, j\in\mathcal{V}_{in}\cup\mathcal{V}_{s}, q\in\mathcal{Q}\label{nodeconse_equ}
\end{align}
where $\mathcal{V}_{in}$ and $\mathcal{V}_{s}$ denote the node sets of intersections and stations, respectively. Formula (\ref{nodeconse_equ}) also implies that private cars are not permitted to dwell within the network, similar to the link connectivity constraint of buses.\par
For the node conservation at entrances and exits, we have:
\begin{align}
    & \lambda_{q}^{r,jk}+\tilde\lambda_{q}^{r,jk}=Tf^r, & \forall r\in\mathcal{R}, j\in\mathcal{V}_o, q\in\mathcal{Q}\\
    & \sum_{\hat{q}\in\mathcal{Q}}{(\mu_{\hat{q}}^{r,ij}+\tilde\mu_{\hat{q}}^{r,ij})}=Hf^r, &\forall r\in\mathcal{R}, j\in\mathcal{V}_d
\end{align}
where $\mathcal{V}_o$ and $\mathcal{V}_d$ denote the node sets of the origins and the destinations.

\subsection{Model formulation} \label{Model_formu_subsec}
To enhance the efficiency of the whole network, the objective of the RC-H design is to minimize weighted travel costs of both buses and private cars. The travel time of private cars is composed of that in mixed-traffic lanes and that in regular lanes, given by:
\begin{align}
    & O_{a}=\sum_{ij\in\mathcal{A}}\sum_{q\in\mathcal{Q}}\sum_{\hat{q}\in\mathcal{Q}}\sum_{r\in\mathcal{R}}(\pi_{q,\hat{q}}^{r,ij}\times t_{q,\hat{q}}^{ij}+\tilde\pi_{q,\hat{q}}^{r,ij}\times t_{a}^{ij})
\end{align}\par
The travel time of buses is composed of that moving on mixed-traffic lanes and dwelling at stations, as given by:
\begin{align}
    & O_{b}=\sum_{p\in\mathcal{P}_{ij}}\gamma_p\times(\sum_{ij\in\mathcal{A}} t_{p}^{ij}+\sum_{j\in\mathcal{V}_s}\hat{t}_{p}^{j})
\end{align}
where $\gamma_p$ is the average number of passengers on bus $p$. To combine the travel time of buses with private cars, the travel time of buses are normalized into passengers’ travel time (Christofa et al.,2016). Moreover, a parameter denoting the efficiency priority between the buses and private cars is also introduced, i.e., $\omega_p\in[0,1]$. Then, the objective function can be expressed as:
\begin{align}
    & O=(1-\omega_p)\times O_{a}+\omega_p\times O_{b}
\end{align}\par
Along with the constraints proposed in the models of realized VP design, bus itinerary planning and traffic assignment of private cars, a mixed integer linear program (MILP) can be formulated, i.e., \textbf{MILP-O}, as follows.
\begin{align}
 \mathrm{\textbf{MILP-O}} \nonumber \\
& \min_{\boldsymbol{\Theta},\boldsymbol{\hat\Theta},\boldsymbol{q},\boldsymbol{\hat{q}},\boldsymbol{\phi},\boldsymbol{f},\boldsymbol{\pi},\boldsymbol{\tilde{\pi}},\boldsymbol{\lambda},\boldsymbol{\tilde{\lambda}},\boldsymbol{\mu},\boldsymbol{\tilde{\mu}}}(1-\omega_p)\times O_{a}+\omega_p\times O_{b} 
\end{align}
\begin{align}
&\mathrm{s.t.} \nonumber\\
& \sum_{\hat{q}\in\mathcal{Q}}\Theta_{q,\hat{q}}^{ij}\leq 1 & \forall ij\in\mathcal{A},q\in\mathcal{Q}\\
& \sum_{q\in\mathcal{Q}}\Theta_{q,\hat{q}}^{ij}\leq 1 & \forall ij\in\mathcal{A},\hat{q}\in\mathcal{Q}\\
& \Theta_{q_m,\hat{q}_m}^{ij}+\Theta_{q_n,\hat{q}_n}^{ij}\leq \zeta_{q_m,\hat{q}_m,q_n,\hat{q}_n}^{ij} & \forall ij\in\mathcal{A}, q_m,\hat{q}_m,q_n,\hat{q}_n\in\mathcal{Q}\\
& \Theta_{q,\hat{q}}^{ij}\in\{0,1\} & \forall ij\in\mathcal{A}, q,\hat{q}\in\mathcal{Q}\\
& \sum_{\hat{q}\in\mathcal{Q}}\sum_{q\in\mathcal{Q}}\hat{\Theta}_{q,\hat{q}}^{p,ij}=1 & \forall ij\in\mathcal{A}, p\in\mathcal{P}_{ij}\\
\begin{split}
& s_a^{ij}\times\hat{\Theta}_{q,\hat{q}}^{ij}\geq s_b \sum_{p\in\mathcal{P}_{ij}}\hat{\Theta}_{q,\hat{q}}^{p,ij}\\
& \hat{\Theta}_{q,\hat{q}}^{ij}\leq \sum_{p\in\mathcal{P}_{ij}}\hat{\Theta}_{q,\hat{q}}^{p,ij}    
\end{split} & \forall ij\in\mathcal{A}, q,\hat{q}\in\mathcal{Q}\\
& \Theta_{q,\hat{q}}^{ij}\geq \hat{\Theta}_{q,\hat{q}}^{ij} & \forall ij\in\mathcal{A}, q,\hat{q}\in\mathcal{Q} \\
& \hat{\Theta}_{q,\hat{q}}^{ij}, \hat{\Theta}_{q,\hat{q}}^{p,ij}\in\{0,1\} & \forall ij\in\mathcal{A}, p\in\mathcal{P}_{ij}, q,\hat{q}\in\mathcal{Q}\\
\begin{split}
& q_p^{ij}=\sum_{\hat{q}\in\mathcal{Q}}\sum_{q\in\mathcal{Q}}q\times\hat{\Theta}_{q,\hat{q}}^{p,ij}\\
& \hat{q}_p^{ij}=\sum_{\hat{q}\in\mathcal{Q}}\sum_{q\in\mathcal{Q}}\hat{q}\times\hat{\Theta}_{q,\hat{q}}^{p,ij}\\
& q_p^{ij},\hat{q}_p^{ij}\in\mathcal{Q}
\end{split} & \forall ij\in\mathcal{A}, p\in\mathcal{P}_{ij}\\
& q_p^{jk}-\hat{q}_p^{ij}=0 & \forall j\in\mathcal{V}_{in}, p\in\mathcal{P}_{ij}\cup\mathcal{P}_{jk}\\
\begin{split}
& \hat{t}_{p}^j=(q_p^{jk}-\hat{q}_p^{ij}+\hat{\phi}_{p}^jQ)\times T \\
& \underline{\hat{t}}_{p}^j\leq \hat{t}_{p}^j\leq H-\epsilon\\
& \hat{\phi}_{p}^j\in\mathbb{N}
\end{split} & \forall j\in\mathcal{V}_s, p\in\mathcal{P}_{ij}\cup\mathcal{P}_{jk}\\
\begin{split}
& t_p^{ij}=\tau_j-\tau_i+(\hat{q}_p^{ij}-q_p^{ij})T+\phi_p^{ij}H\\
& \underline{t}_p^{ij}\leq t_p^{ij}\leq H+t_a^{ij}-\epsilon\\
& \phi_{p}^{ij}\in\mathbb{N}
\end{split} & \forall ij\in\mathcal{A}, p\in\mathcal{P}_{ij}\\
& \sum_{r\in\mathcal{R}^w}f^r=d^w & \forall w\in\mathcal{W}\\
\begin{split}
& \sum_{q\in\mathcal{Q}}\pi_{q,\hat{q}}^{r,ij}=\mu_{\hat{q}}^{r,ij}\\
& \sum_{\hat{q}\in\mathcal{Q}}\pi_{q,\hat{q}}^{r,ij}=\lambda_{q}^{r,ij}
\end{split} & \forall r\in\mathcal{R}, ij\in\mathcal{A}, q,\hat{q}\in\mathcal{Q}\\
& \sum_{r\in\mathcal{R}}{\pi_{q,\hat{q}}^{r,ij}}\leq s_a^{ij}\times (\Theta_{q,\hat{q}}^{ij}-\hat{\Theta}_{q,\hat{q}}^{ij}) & \forall ij\in\mathcal{A}, q, \hat{q}\in\mathcal{Q} \\
\begin{split}
& \sum_{q\in\mathcal{Q}}\tilde{\pi}_{q,\hat{q}}^{r,ij}=\tilde{\mu}_{\hat{q}}^{r,ij}\\
& \sum_{\hat{q}\in\mathcal{Q}}\tilde{\pi}_{q,\hat{q}}^{r,ij}=\tilde{\lambda}_{q}^{r,ij}
\end{split} & \forall r\in\mathcal{R}, ij\in\mathcal{A}, q,\hat{q}\in\mathcal{Q}\\
& \sum_{r\in\mathcal{R}}{\tilde{\pi}_{q,\hat{q}}^{r,ij}}\leq s_a^{ij}\times(l_{ij}-1)\times\tilde\Theta_{q,\hat{q}}^{ij} &\forall ij\in\mathcal{A}, q, \hat{q}\in\mathcal{Q} \\
& (\lambda_{q}^{r,jk}+\tilde\lambda_{q}^{r,jk})-(\mu_{q}^{r,ij}+\tilde\mu_{q}^{r,ij})=0 & \forall r\in\mathcal{R}, j\in\mathcal{V}_{in}\cup\mathcal{V}_{s}, q\in\mathcal{Q}\\
& \lambda_{q}^{r,jk}+\tilde\lambda_{q}^{r,jk}=Tf^r & \forall r\in\mathcal{R}, j\in\mathcal{V}_o, q\in\mathcal{Q}\\
& \sum_{\hat{q}\in\mathcal{Q}}{(\mu_{\hat{q}}^{r,ij}+\tilde\mu_{\hat{q}}^{r,ij})}=Hf^r &\forall r\in\mathcal{R}, j\in\mathcal{V}_d
\end{align}

\section{Bilevel solution method} \label{Two_level_sec}
In \textbf{MILP-O}, the numbers of both variables and constraints are enormous and the relations among the variables are complicated; thus, it is burdensome to obtain the optimal solution by directly solving the MILP. In this section, we propose a bilevel solution method for \textbf{MILP-O} to relieve the computational cost, that can provide sub-optimal solutions with reasonable computational efforts and acceptable solution quality.

\subsection{Problem decomposition} \label{Decomp_subsec}
\begin{figure}
    \centering
    \includegraphics[width=0.55\textwidth]{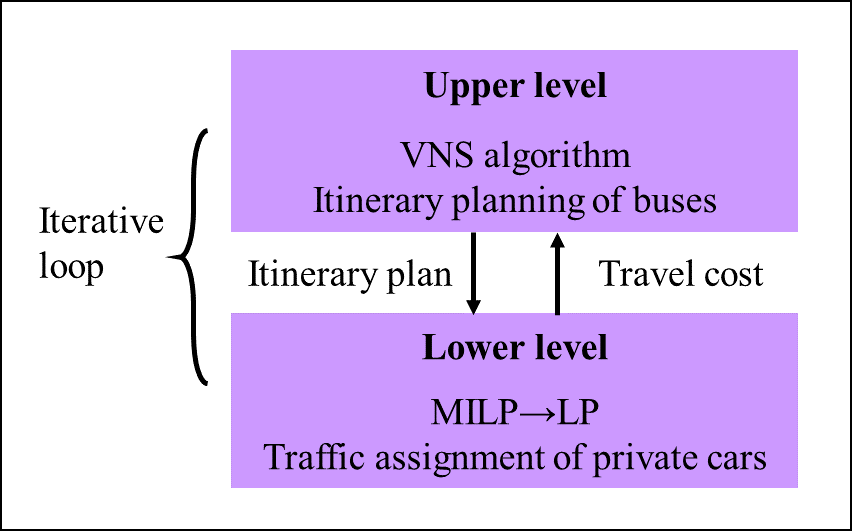}
    \caption{Illustration for the bilevel solution structure}
    \label{Fig_twolevel}
\end{figure}

Considering the scale of the problem, we decompose the original model into two levels, i.e., the upper-level model and the lower-level model. The upper-level model is focused on designing the itinerary of buses along bus route, while the lower-level model solves the problem of traffic assignment of private cars with the given bus itinerary plan. An iterative loop between the two levels is built up as shown in Fig.\ref{Fig_twolevel}. Beginning from the initial solution, the upper level provides an itinerary plan of buses for the lower level, based on which the lower level then solves the traffic assignment problem; after obtaining the solution from the lower level, the performance of the bus planning can be evaluated and coupled back to the upper level; in turn, the upper level improves the planning according to the feedback from the lower level, and then triggers the next iteration. The loop will not be terminated until a preset maximum iteration number is reached, when both the bus itinerary planning and traffic assignment of private cars can be acquired concurrently.

\subsection{Variable neighborhood search for upper level} \label{upper_subsec}
Variable neighborhood search (VNS) algorithm is applied for the upper level; it is a metaheuristic method capable of solving problems with a large scale or a complicated structure. VNS has the advantages of both descent to local minima and escape from the valley valley  (Hansen  et  al., 2006). The idea of VNS algorithm is straight-forward; it improves the solution iteratively by systematic changes of the neighborhood with a local search algorithm (Mladenovic and Hansen, 1997). Specifically, in each iteration, taking the current incumbent solution as the search origin, a solution in the neighborhood will be generated as a new solution through a local search routine; then the performance of the new solution will be measured to determine whether to accept it and instruct the search neighborhood for the next iteration. If the solution has not been improved for a certain number of iterations, the size of the neighborhood scale will be adjusted.\par

For applying the VNS algorithm to the upper level, a plan of the bus itinerary will be generated in every iteration. The variables for the planning problem include the arrival time, the travel time on links and the dwelling time at stations for each bus. From the results of numerical experiments presented in the next section, it can be observed that the travel time of buses on links is always optimized to the minimum value, that is, $t_p^{ij}=\underline{t}_p^{ij}$; it is reasonable as the smaller travel time of buses is, the less delay of buses will incur, and concurrently the less delay of private cars will also be induced due to the reduced interference from buses. Therefore, the travel time of buses is fixed as the minimum value in the upper level, reducing the number of variables remained to be determined.\par
\textbf{Algorithm 1} describes the local search routine of a neighbor solution of the upper level. $\bm{x}^{*}$ denotes the current incumbent solution, which are updated continuously toward the direction of performance improvement. $\bm{x}^{k}$ denotes the newly generated solution in iteration $k$ with the searching origin of $\bm{x}^{*}$ and the searching radius $\Delta_k$. $\Delta_k$ embodies the scale of the searching neighborhood and will be updated according to $\Delta_k=\Delta_k+\epsilon_\Delta$ when the incumbent solution $\bm{x}^{*}$ has no update within a preset iteration times; $p_l^k$ and $p_u^k$ are two “instruction parameters” that determine the variation direction of the variables. Note that $p_l^k$ and $p_u^k$ are updated with iteration number $k$ to promote the quality of the generated solution; specifically, in the early searching process, more freedom is provided for the solution variation to create more possibilities for the solution set and avoid the solution from falling into local optima too early; as the iteration number increases, the variation tendency becomes restricted to focus more on the current solution as most variables are already at their optimal values. Although the proposed searching method cannot guarantee the optimality of solution, it does indeed provide sub-optimal solutions with acceptable quality, as will be validated in the numerical experiments.\par

\noindent\rule{\textwidth}{0.2pt}
\noindent \textbf{Algorithm 1. Local search algorithm} \par
\hangafter 1
\hangindent 3.5em
\noindent \textbf{Step 1}. Generate a random number for each variable $x_i^k$, $p_i^k\in[0,1]$, and calculate the values of the instruction parameters in iteration $k$, $p_u^k$ and $p_l^k$; \par
\hangafter 1
\hangindent 3.5em
\noindent \textbf{Step 2}. Compare the generated number with the instruction parameters; if $p_i^k<p_l^k$, set $x_i^k=x_i^k-\Delta_k$; if $p_i^k>p_u^k$, set $x_i^k=x_i^k+\Delta_k$; otherwise, $x_i^k$ remains unchanged; \par
\hangafter 1
\hangindent 3.5em
\noindent \textbf{Step 3}. Output new solution $\bm{x}^k$.\par
\noindent\rule{\textwidth}{0.2pt}

\subsection{Linear programming approximation for lower level} \label{lower_subsec}
In the lower-level problem, with the given bus itinerary plan, the traffic assignment of private cars is optimized; the optimization is formulated as an MILP problem, as shown for \textbf{MILP-L}:
\begin{align}
 \mathrm{\textbf{MILP-L}} \nonumber \\
& O_L=\min_{\boldsymbol{\Theta},\boldsymbol{f},\boldsymbol{\pi},\boldsymbol{\tilde{\pi}},\boldsymbol{\lambda},\boldsymbol{\tilde{\lambda}},\boldsymbol{\mu},\boldsymbol{\tilde{\mu}}}(1-\omega_p)\times O_a
\end{align}
\begin{align}
&\mathrm{s.t.} \nonumber\\
& \sum_{\hat{q}\in\mathcal{Q}}\Theta_{q,\hat{q}}^{ij}\leq 1 & \forall ij\in\mathcal{A},q\in\mathcal{Q}\\
& \sum_{q\in\mathcal{Q}}\Theta_{q,\hat{q}}^{ij}\leq 1 & \forall ij\in\mathcal{A},\hat{q}\in\mathcal{Q}\\
& \Theta_{q_m,\hat{q}_m}^{ij}+\Theta_{q_n,\hat{q}_n}^{ij}\leq \zeta_{q_m,\hat{q}_m,q_n,\hat{q}_n}^{ij} & \forall ij\in\mathcal{A}, q_m,\hat{q}_m,q_n,\hat{q}_n\in\mathcal{Q}\label{FIFO_eq}\\
& \Theta_{q,\hat{q}}^{ij}\in\{0,1\} & \forall ij\in\mathcal{A}, q,\hat{q}\in\mathcal{Q}\\
& \Theta_{q,\hat{q}}^{ij}\geq \hat{\Theta}_{q,\hat{q}}^{ij} & \forall ij\in\mathcal{A}, q,\hat{q}\in\mathcal{Q} \\
& \sum_{r\in\mathcal{R}^w}f^r=d^w & \forall w\in\mathcal{W}\\
\begin{split}
& \sum_{q\in\mathcal{Q}}\pi_{q,\hat{q}}^{r,ij}=\mu_{\hat{q}}^{r,ij}\\
& \sum_{\hat{q}\in\mathcal{Q}}\pi_{q,\hat{q}}^{r,ij}=\lambda_{q}^{r,ij}
\end{split} & \forall r\in\mathcal{R}, ij\in\mathcal{A}, q,\hat{q}\in\mathcal{Q}\\
& \sum_{r\in\mathcal{R}}{\pi_{q,\hat{q}}^{r,ij}}\leq s_a^{ij}\times (\Theta_{q,\hat{q}}^{ij}-\hat{\Theta}_{q,\hat{q}}^{ij}) & \forall ij\in\mathcal{A}, q, \hat{q}\in\mathcal{Q} \\
\begin{split}
& \sum_{q\in\mathcal{Q}}\tilde{\pi}_{q,\hat{q}}^{r,ij}=\tilde{\mu}_{\hat{q}}^{r,ij}\\
& \sum_{\hat{q}\in\mathcal{Q}}\tilde{\pi}_{q,\hat{q}}^{r,ij}=\tilde{\lambda}_{q}^{r,ij}
\end{split} & \forall r\in\mathcal{R}, ij\in\mathcal{A}, q,\hat{q}\in\mathcal{Q}\\
& \sum_{r\in\mathcal{R}}{\tilde{\pi}_{q,\hat{q}}^{r,ij}}\leq s_a^{ij}\times(l_{ij}-1)\times\tilde\Theta_{q,\hat{q}}^{ij} &\forall ij\in\mathcal{A}, q, \hat{q}\in\mathcal{Q} \\
& (\lambda_{q}^{r,jk}+\tilde\lambda_{q}^{r,jk})-(\mu_{q}^{r,ij}+\tilde\mu_{q}^{r,ij})=0 & \forall r\in\mathcal{R}, j\in\mathcal{V}_{in}\cup\mathcal{V}_{s}, q\in\mathcal{Q}\\
& \lambda_{q}^{r,jk}+\tilde\lambda_{q}^{r,jk}=Tf^r & \forall r\in\mathcal{R}, j\in\mathcal{V}_o, q\in\mathcal{Q}\\
& \sum_{\hat{q}\in\mathcal{Q}}{(\mu_{\hat{q}}^{r,ij}+\tilde\mu_{\hat{q}}^{r,ij})}=Hf^r &\forall r\in\mathcal{R}, j\in\mathcal{V}_d
\end{align}
In the \textbf{MILP-L}, $\hat{\Theta}_{q,\hat{q}}^{ij}$ is the input from the upper level, which denotes the VPs dedicated for buses. Under the premise of a given $\hat{\Theta}_{q,\hat{q}}^{ij}$, the traffic assignment of private cars can be solved, which is much more tractable compared to \textbf{MILP-O} and is solvable for a large-scale problem.\par
Moreover, to accelerate the solution process of the lower level, several approximation techniques are applied to further reduce the solution space and simplify the problem structure in the loop iteration, as described below.\par
\begin{enumerate}
    \item The VPs slower than the bus traveling on the same link are abandoned. As the objective of the optimization is to minimize the total travel cost, the likelihood for the VPs with low speed to be realized for private cars is limited due to the large delay. Then, the candidate set of VPs can be shrunk as shown in Fig.\ref{Fig_abadonslow}. \par
\begin{figure}[!ht]
	\centering
	\subfloat[][Before shrink]{\includegraphics[width=0.45\textwidth]{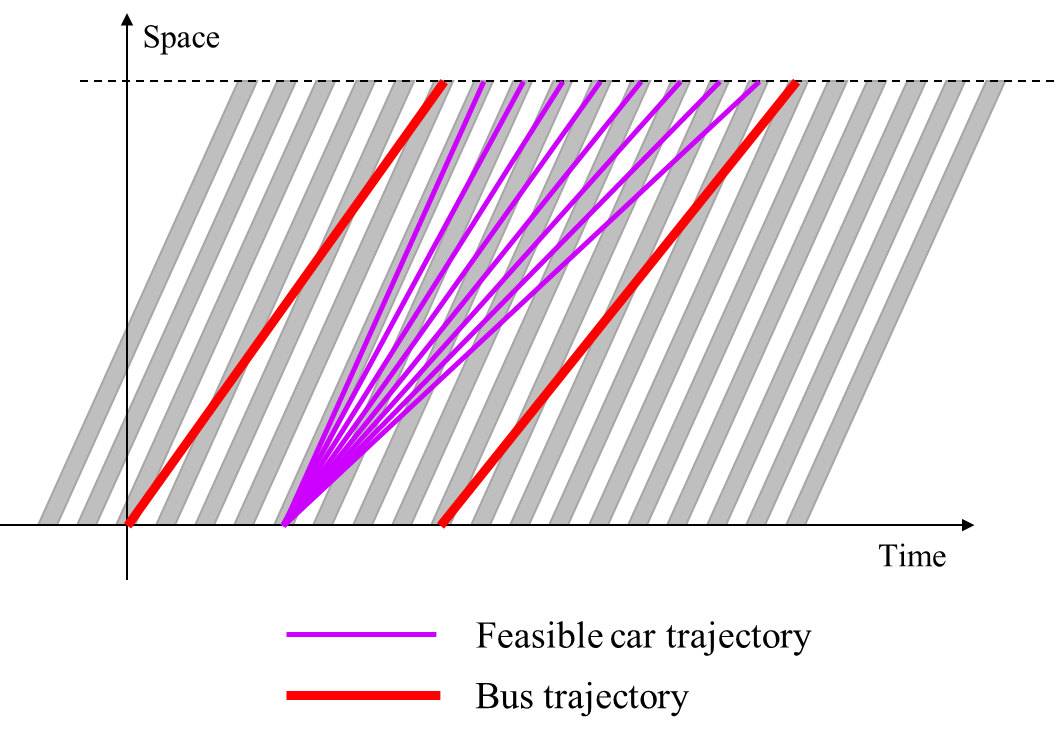}}\hspace{1em}
	\subfloat[][After shrink]{\includegraphics[width=0.45\textwidth]{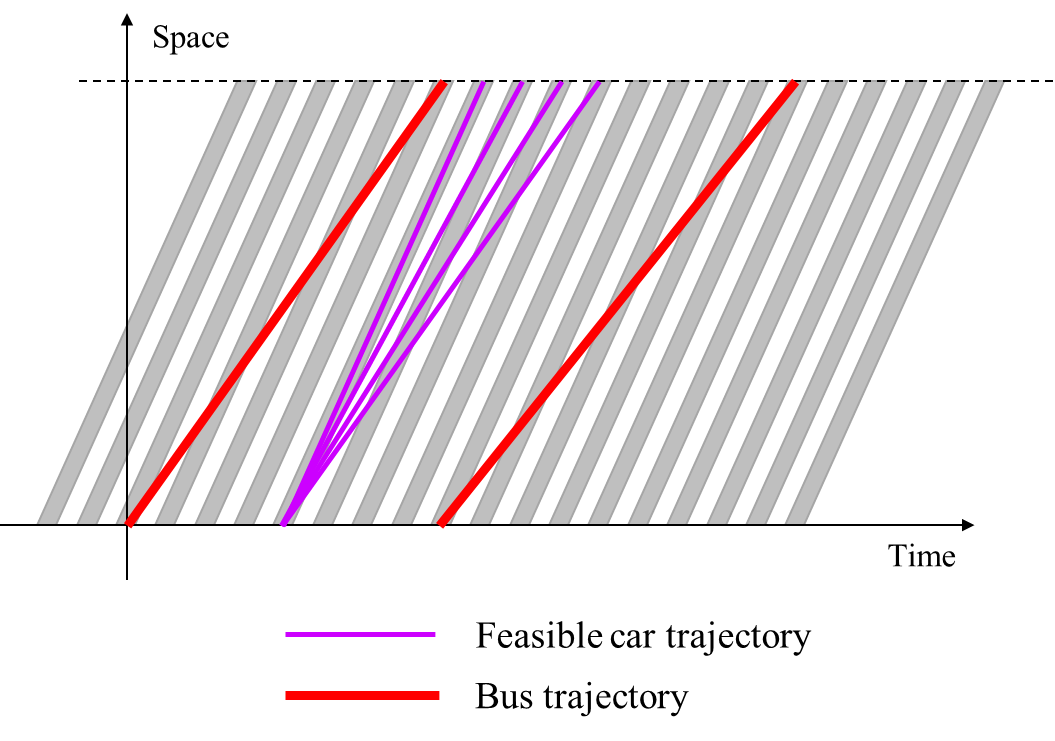}}\\
	\caption[]{ Illustration of abandoning slow virtual platoons}
	\label{Fig_abadonslow}
\end{figure}

\item The FIFO constraint (\ref{FIFO_eq}) is removed in solving the lower-level problem. As the FIFO constraint is imposed on any two different VPs for each link, it generates a large number of constraints for the problem and slows down the solution process to a great extent. From the numerical experiments in the next section, it can be validated that the removal of the FIFO constraint has limited impacts on the results, as few VPs violating the FIFO constraint are realized. Moreover, with some adjustment methods, the violation of the FIFO constraint can be eliminated without affecting the total cost. Fig.\ref{Fig_adjustment} illustrates an example of the adjustment method, where Nos.1 and 2 are the  trajectories violating FIFO constraint and Nos.3 and 4 are the trajectories after re-assignment. Note that this adjustment method is only applicable for the vehicles with the same destination.\par
\begin{figure}[!ht]
	\centering
	\subfloat[][Before adjustment]{\includegraphics[width=0.45\textwidth]{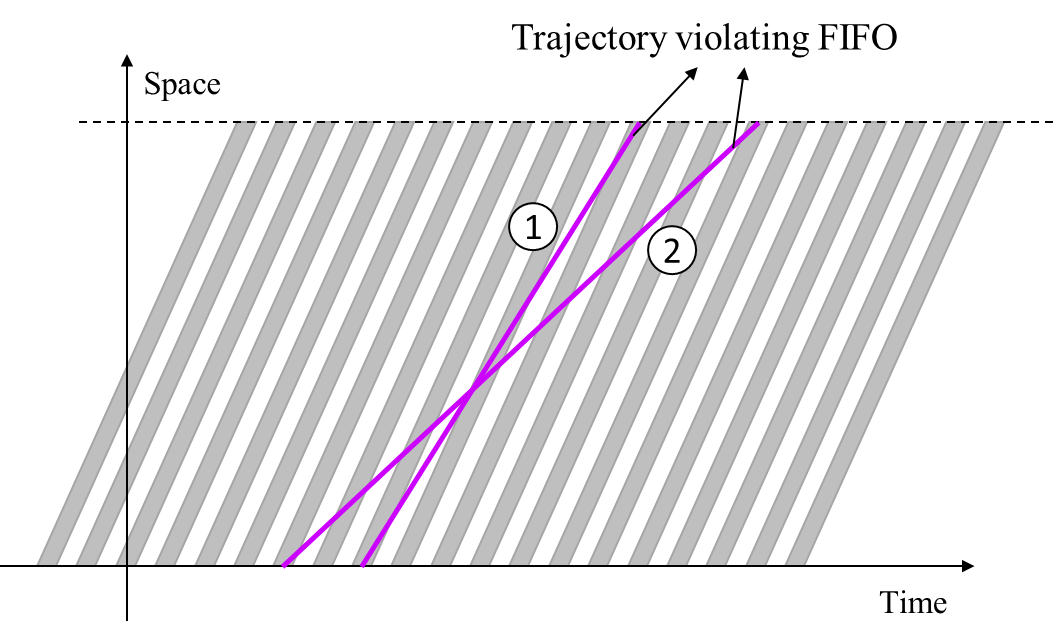}}\hspace{1em}
	\subfloat[][After adjustment]{\includegraphics[width=0.45\textwidth]{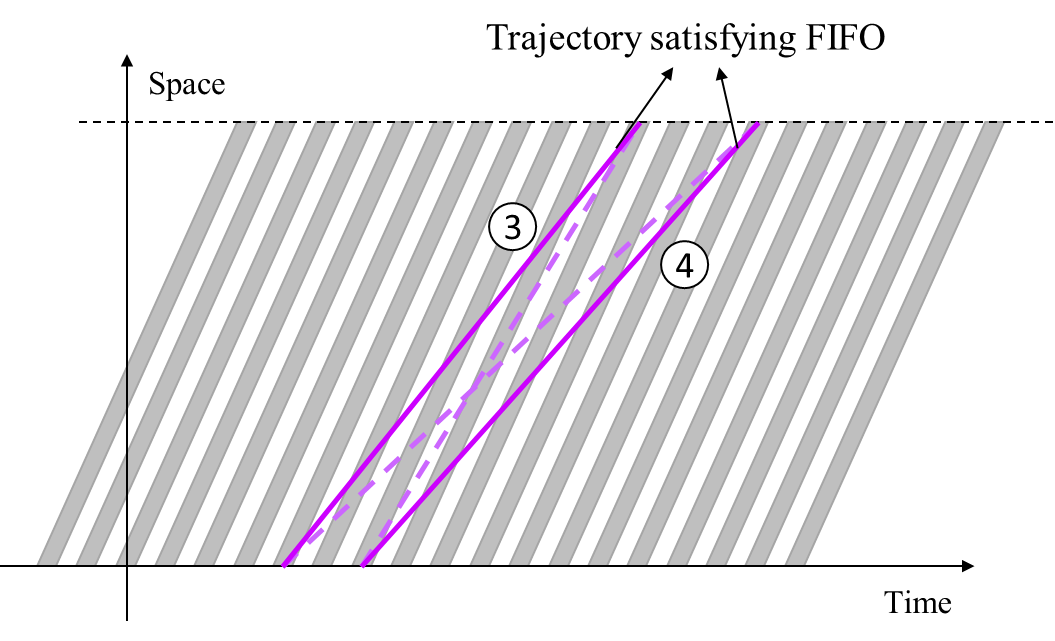}}\\
	\caption[]{Illustration of an adjustment method}
	\label{Fig_adjustment}
\end{figure}
\item The VPs are allowed to be split for different traffic flows. This technique aims to remove the binary variable of $\Theta_{q,\hat{q}}^{ij}$, and then the problem can be re-formulated as a linear program (LP), which can be solved with a high computational efficiency. Similar to technique (2), the impact of the LP approximation is validated in the numerical experiments, where it shows that the approximation has limited impacts as most realized VPs remain non-split.\par
\end{enumerate}

Then, using these three approximation techniques, the lower-level problem can be formulated as an LP with reduced solution space, as shown by \textbf{LP-L}:
\begin{align}
 \mathrm{\textbf{LP-L}} \nonumber \\
& O_L'=\min_{\boldsymbol{f},\boldsymbol{\pi},\boldsymbol{\tilde{\pi}},\boldsymbol{\lambda},\boldsymbol{\tilde{\lambda}},\boldsymbol{\mu},\boldsymbol{\tilde{\mu}}}(1-\omega_p)\times O_a
\end{align}
\begin{align}
&\mathrm{s.t.} \nonumber\\
& \sum_{r\in\mathcal{R}^w}f^r=d^w & \forall w\in\mathcal{W}\\
\begin{split}
& \sum_{q\in\mathcal{Q}}\pi_{q,\hat{q}}^{r,ij}=\mu_{\hat{q}}^{r,ij}\\
& \sum_{\hat{q}\in\mathcal{Q}}\pi_{q,\hat{q}}^{r,ij}=\lambda_{q}^{r,ij}
\end{split} & \forall r\in\mathcal{R}, ij\in\mathcal{A}, q,\hat{q}\in\mathcal{Q}\\
& \sum_{r\in\mathcal{R}}{\pi_{q,\hat{q}}^{r,ij}}\leq s_a^{ij}\times \dot{\Theta}_{q,\hat{q}}^{ij} & \forall ij\in\mathcal{A}, q, \hat{q}\in\mathcal{Q} \\
\begin{split}
& \sum_{q\in\mathcal{Q}}\tilde{\pi}_{q,\hat{q}}^{r,ij}=\tilde{\mu}_{\hat{q}}^{r,ij}\\
& \sum_{\hat{q}\in\mathcal{Q}}\tilde{\pi}_{q,\hat{q}}^{r,ij}=\tilde{\lambda}_{q}^{r,ij}
\end{split} & \forall r\in\mathcal{R}, ij\in\mathcal{A}, q,\hat{q}\in\mathcal{Q}\\
& \sum_{r\in\mathcal{R}}{\tilde{\pi}_{q,\hat{q}}^{r,ij}}\leq s_a^{ij}\times(l_{ij}-1)\times\tilde\Theta_{q,\hat{q}}^{ij} &\forall ij\in\mathcal{A}, q, \hat{q}\in\mathcal{Q} \\
& (\lambda_{q}^{r,jk}+\tilde\lambda_{q}^{r,jk})-(\mu_{q}^{r,ij}+\tilde\mu_{q}^{r,ij})=0 & \forall r\in\mathcal{R}, j\in\mathcal{V}_{in}\cup\mathcal{V}_{s}, q\in\mathcal{Q}\\
& \lambda_{q}^{r,jk}+\tilde\lambda_{q}^{r,jk}=Tf^r & \forall r\in\mathcal{R}, j\in\mathcal{V}_o, q\in\mathcal{Q}\\
& \sum_{\hat{q}\in\mathcal{Q}}{(\mu_{\hat{q}}^{r,ij}+\tilde\mu_{\hat{q}}^{r,ij})}=Hf^r &\forall r\in\mathcal{R}, j\in\mathcal{V}_d
\end{align}
In the formulation of \textbf{LP-L}, $\dot{\Theta}_{q,\hat{q}}^{ij}$ denotes whether the corresponding VP can be realized, which is determined by the results of bus itinerary planning, i.e., $\hat{\Theta}_{q,\hat{q}}^{ij}$, as shown in Fig.\ref{Fig_abadonslow}; if it can be realized, $\dot{\Theta}_{q,\hat{q}}^{ij}=1$; otherwise, $\dot{\Theta}_{q,\hat{q}}^{ij}=0$. As the bus itinerary planning is determined by the upper bound, $\dot{\Theta}_{q,\hat{q}}^{ij}$ has been preset when solving the lower-level model; thus, all of the variables are continuous as expected.\par
Note that the \textbf{LP-L} is only applied in the iteration loop to accelerate the solving process, and \textbf{MILP-L} will be solved after the termination of the iteration, as the exact solution provided by the bilevel solution procedure.

\subsection{Iterative algorithm} \label{iterative_subsec}
\textbf{Algorithm 2} describes the detailed iterative process of the bilevel solution method. $S^{*}$ denotes the current incumbent solution of the bilevel problem and $S^{k}$ denotes the newly generated solution in iteration $k$, obtained by solving the \textbf{LP-L} with the input of $\bm{x}^k$. If the new solution $S^{k}$ is better than $S^{*}$, $S^{*}$ will be updated by $S^{k}$, and $\bm{x}^{*}$ will be updated by $\bm{x}^k$ as the new searching origin for \textbf{Algorithm 1}; otherwise, $S^{k}$ will be abandoned and the next iteration will be triggered.\par
\noindent\rule{\textwidth}{0.2pt}
\noindent \textbf{Algorithm 2. Iterative algorithm} \par
\hangafter 1
\hangindent 3.5em
\noindent {\textbf{Step 1}}. Set the initial values of step length $\Delta_0$ and instruction parameters $p_u^0$, $p_l^0$, $0\leq p_l^0<p_u^0\leq1$; obtain the initial solution of upper level, $\bm{x}^0$, and input it to the lower level to obtain the initial solution, $S_0$. Set the iteration number $k=1$ and let $\Delta_k=\Delta_0, p_l^k=p_l^0,p_u^k=p_u^0,\bm{x}^{*}=\bm{x}^0,S^{*}=S_0$; \par
\hangafter 1
\hangindent 3.5em
\noindent {\textbf{Step 2}}. Generate a new solution $\bm{x}^{k}$ in the neighborhood of $\bm{x}^{*}$ by \textbf{Algorithm 1} with the parameters of $\Delta_k$, $p_l^k$ and $p_u^k$; \par
\hangafter 1
\hangindent 3.5em
\noindent {\textbf{Step 3}}. Input solution $\bm{x}^{k}$ to the lower level and solve the \textbf{LP-L} to obtain a new solution of the bilevel problem, $S_k$; \par
\hangafter 1
\hangindent 3.5em
\noindent {\textbf{Step 4}}. Compare the new solution $S_k$ with current incumbent solution $S^{*}$; if the new solution is better, update the current incumbent solution by $S^{*}=S_k$, $\bm{x}^{*}=\bm{x}^{k}$; otherwise, abandon the new solution;\par
\hangafter 1
\hangindent 3.5em
\noindent {\textbf{Step 5}}. Update the iteration number as $k=k+1$. If $k$ is beyond the maximum iteration time, stop the iteration and go to Step 6; otherwise, update the value of $\Delta_k$, $p_l^k$ and $p_u^k$, and go back to Step 2;\par
\hangafter 1
\hangindent 3.5em
\noindent {\textbf{Step 6}}. Input solution $\bm{x}^{*}$ to the lower level and solve the \textbf{MILP-L} to obtain the exact solution of the bilevel problem.\par
\noindent\rule{\textwidth}{0.2pt}

\section{Numerical experiments} \label{Numerical_sec}
In the numerical experiments, two scenarios are tested, i.e., a toy example and a real-world network implementation. In the toy example, we test the proposed RC-H scheme in a small-scale network. As the optimal solution of \textbf{MILP-O} can be obtained for the small-scale problem, the performance of the RC-H scheme will be measured and the performances characteristics of the bilevel solution method will be validated. In the real-world network scenario, the applicability of the proposed methods on a large-scale problem in practice will be verified, and the advantages compared with the traffic signal control strategies will also be demonstrated by simulation experiments.
\subsection{Toy example} \label{toy_exam_subsec}

\subsubsection{Experimental setups} \label{setup_toy_subsubsec}

\noindent\underline{\textit{Road network settings}}\par
\begin{figure}[!ht]
	\centering
	{\includegraphics[width=0.95\textwidth]{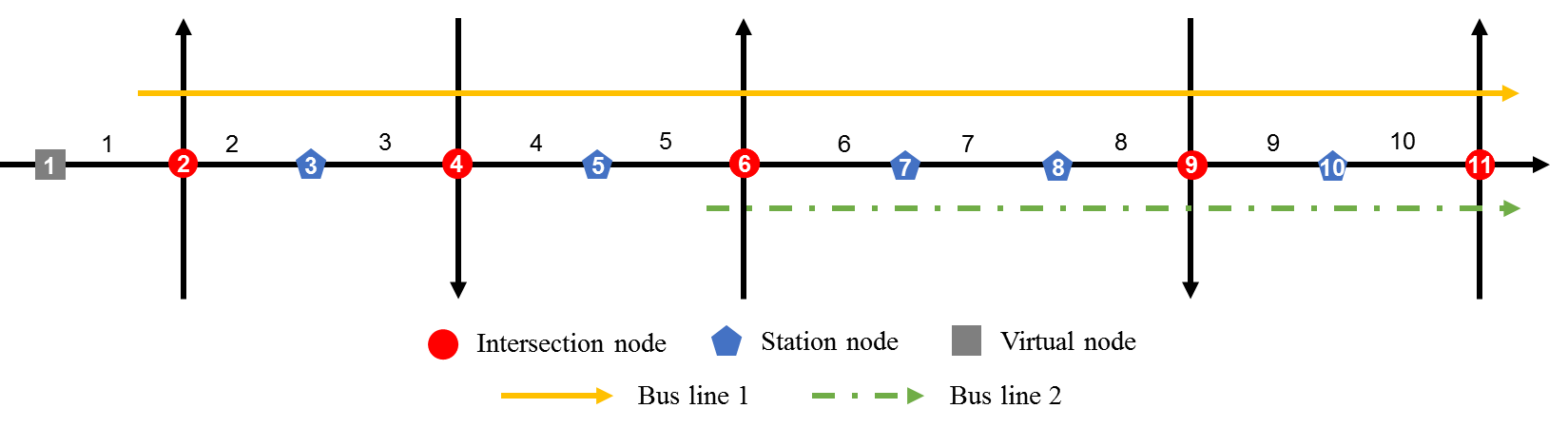}}
	\caption[]{Illustration for toy example scenario}
	\label{Fig_toy_network}
\end{figure}

We consider a simple scenario of one corridor as shown in Fig.\ref{Fig_toy_network}. In this scenario, there are 11 nodes composed of 5 intersection nodes, 5 station nodes and 1 virtual node, and 10 links with one lane on each link. Note that node 1 is designed as a virtual node and link 1 is a virtual link for the traffic organization in the entrance with no actual length, which can be recognized as a waiting zone for vehicles prior to their entry into the network.
\par

\vspace{0.5em}
\noindent\underline{\textit{Traffic demand settings}}\par
Two bus lines are considered: bus line 1 from node 2 to node 11, and bus line 2 from node 6 to node 11; note that there are shared links for the two bus lines as shown in Fig.\ref{Fig_toy_network}. One O-D pair of private car traffic is also set from node 2 to node 11 with various demand levels. 

\vspace{0.5em}
\noindent\underline{\textit{RC-H scheme settings}}\par
The cycles of background RC scheme and bus system are set as $T=10s$ and $H=120s$, respectively, so that the number of background VPs for one RC-H cycle is $Q=12$. The free-flow travel time of private cars and buses on each link, i.e., $\underline{t}_a^{ij}$ and $\underline{t}_p^{ij}$, are listed in Table \ref{toy_travtim_table}, and the minimum dwelling time at stations of buses is set as $\underline{\hat{t}}_{p}^{j}=40s$. The VP size and bus size are set as $s_a=4$ and $s_b=2$, respectively, implying that one VP can accommodate at most two buses. The number of passengers on buses is set as $\gamma_p=20$.
\begin{table}[]
    \centering
    \caption{Free-flow travel time of vehicles}
    \label{toy_travtim_table}
    \begin{tabular}{c|c c c c c c c c c c}
    \hline
    Link & 1 & 2 & 3 & 4 & 5 & 6 & 7 & 8 & 9 & 10\\
    \hline
    $\underline{t}_a^{ij}(s)$ & 0 & 10 & 20 & 10 & 20 & 10 & 10 & 20 & 10 & 20\\
    $\underline{t}_p^{ij}(s)$ & 0 & 20 & 30 & 20 & 30 & 20 & 20 & 30 & 20 & 30\\
    \hline
    \end{tabular}
\end{table}

\vspace{0.5em}
\noindent\underline{\textit{Heuristic algorithm settings}}\par
The parameters of the upper-level VNS algorithm are set as follows: the initial step length is the minimum unit of variables, i.e, $\Delta_0=T$, and the variation of the step length is also set as $\epsilon_\Delta=T$, and the upper bound of the step length is set as $H$. For the initial solution of the upper-level model, to guarantee the feasibility of solving the lower-level problem, the buses with shared routes collaborate to travel together, reducing the number of dedicated VPs and then reserve more capacity for private cars. The maximum iteration number for the iterative algorithm is set to $2000$.

\subsubsection{Results of \textbf{MILP-O}} \label{MILP_toy_subsubsec}
First, we investigate the maximum traffic volume of private cars that can be accommodated by the RC-H scheme, defined as “maximum admissible traffic”. For one RC-H cycle, there are $Q$ background VPs for each lane, where at least two of them will serve for the dedicated VP of a single bus when the dwelling process exists, as shown in Fig.\ref{Fig_dedicated_dwell}(b); therefore, the maximum admissible traffic of private cars in a single cycle can be calculated as $(Q-2)\times s_a$. Taking the maximum admissible traffic as the upper bound of the traffic demand of private cars, i.e., $\beta_1=1$, three other demands are set as low level ($\beta_1=0.2$), medium level ($\beta_1=0.5$) and high level ($\beta_1=0.8$), respectively. Moreover, two priorities are considered as bus priority ($\omega_p=0.9$) and private car priority ($\omega_p=0.1$) in the tests. Therefore, a total of $8$ cases are tested as shown in Table \ref{toy_MILP_table)}. \par
\begin{table}[]
    \centering
    \caption{Results of \textbf{MILP-O}}
    \label{toy_MILP_table)}
    \begin{tabular}{c|c c c c c c c c}
    \hline
    Case & 1 & 2 & 3 & 4 & 5 & 6 & 7 & 8\\
    \hline
    $\omega_p$ & 0.9 & 0.9 & 0.9 & 0.9 & 0.1 & 0.1 & 0.1 & 0.1 \\
    $\beta_1$ & 0.2 & 0.5 & 0.8 & 1 & 0.2 & 0.5 & 0.8 & 1\\
    $O_a^{m}$ & 1040 & 2600 & 4160 & 5200 & 1040 & 2600 & 4160 & 5200\\
    $O_a^{opt}$ & 1153 & 2960 & 6160 & 9600 & 1153 & 2860 & 5040 & 6800\\
    $\Delta O_a$ & 10.90\% & 13.85\% & 48.08\% & 84.62\% & 10.90\% & 10.00\% & 21.15\% & 30.77\% \\
    $O_b^{m}$ & 13200 & 13200 & 13200 & 13200 & 13200 & 13200 & 13200 & 13200\\
    $O_b^{opt}$ & 13200 & 13200 & 13200 & 13200 & 13200 & 13600 & 16400 & 24400\\    
    $\Delta O_b$ & 0.00\% &0.00\% &0.00\% &0.00\% &0.00\% & 3.03\% & 24.24\% & 84.85\% \\
    $O^{m}$ &14240 & 15800 & 17360 & 18400 & 14240 & 15800 & 17360 &	18400\\
    $O^{opt}$ & 14353 & 16160 & 19360 & 22800 & 14353 & 16460 & 21440 & 31200\\
    $\Delta O$ & 0.80\% & 2.28\% & 11.52\% & 23.91\% & 0.80\% & 4.18\% & 23.50\% & 69.57\% \\
    \hline
    \end{tabular}
\end{table}
The optimization problem is implemented by YALMIP and solved by CPLEX 12.6 solver using MATLAB 2016a on a desktop with a 3.6-GHz Intel Core and 16 GB of RAM. Table \ref{toy_MILP_table)} shows the results of \textbf{MILP-O}, where $O_a^m$ and $O_a^{opt}$ denote the free-flow total travel cost and the cost under the MILP solution of private cars, respectively, and $\Delta O_a$ denotes the gap between $O_a^m$ and $O_a^{opt}$, which reflects the traffic efficiency; $O_b^m$, $O_b^{opt}$ and $\Delta O_b$ have similar meanings but for buses, and $O^m$, $O^{opt}$ and $\Delta O$ are for the total cost of private cars and buses. As the travel time of buses is calculated with the unit of person, i.e., multiplied by the passenger number, $O^m$, $O^{opt}$ and $\Delta O$ can be regarded as the summation of personal travel cost. \par
From Table \ref{toy_MILP_table)}, the results are affected by both traffic demand level and efficiency priority. The main observations are concluded as follows:
\begin{itemize}
    \item For private cars, the travel cost increases with traffic demand, and the rate of increase under bus priority is higher than that under private car priority as the VPs are designed preferentially for dedicated VPs under bus priority.
    \item For buses, the travel cost also increases with traffic demand under private car priority as the relationship of buses and private cars is undoubtedly competitive for the limited VPs. However, under bus priority, the travel cost of buses holds the minimum value and has no variation with the traffic demand increase; it implies that by setting bus priority, the performance of bus operation can be guaranteed under all demand cases. 
    \item From the view of personal cost, the travel cost under bus priority is lower than that under private car priority, especially under high demand, indicating that giving priority for buses can reduce the personal delay and is consistent to the encouragement of the transit application. Hence, bus priority is recommended for the personal efficiency and is set for the following tests.
\end{itemize}\par
\begin{figure}[!ht]
	\centering
	\subfloat[][Under bus priority]{\includegraphics[width=0.45\textwidth]{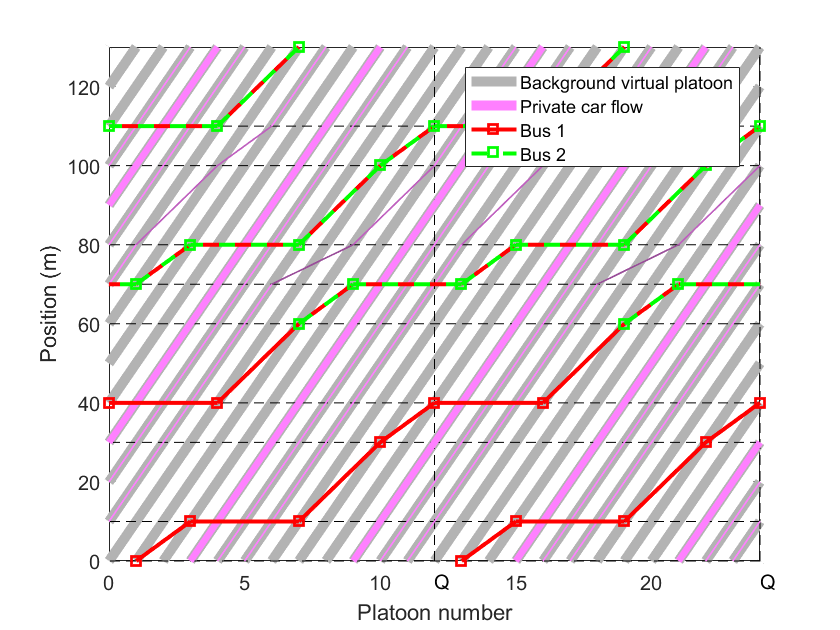}}\hspace{1em}
	\subfloat[][Under private car priority]{\includegraphics[width=0.45\textwidth]{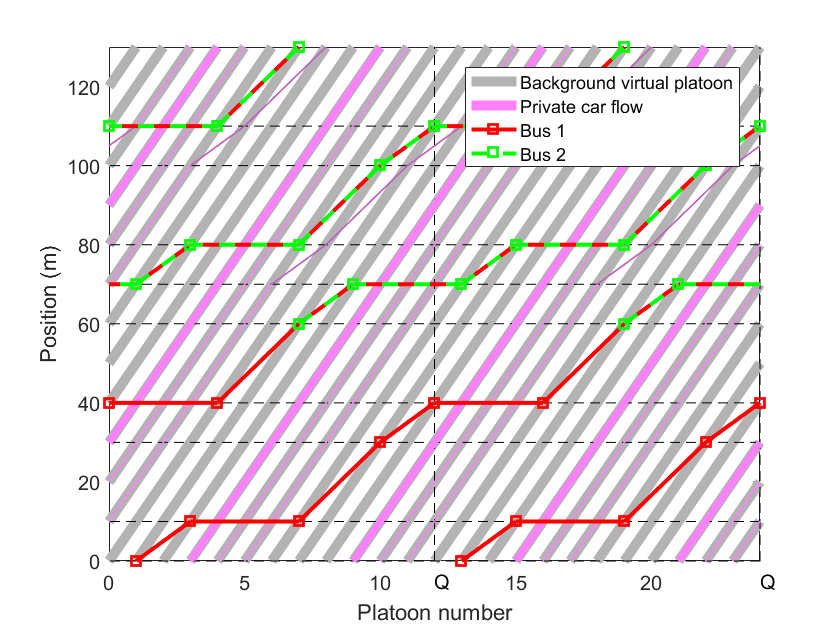}}\\
	\caption[]{Results of \textbf{MILP-O} under low demand}
	\label{Fig_toy_resu_MILP1}
\end{figure}
\begin{figure}[!ht]
	\centering
	\subfloat[][Under bus priority]{\includegraphics[width=0.45\textwidth]{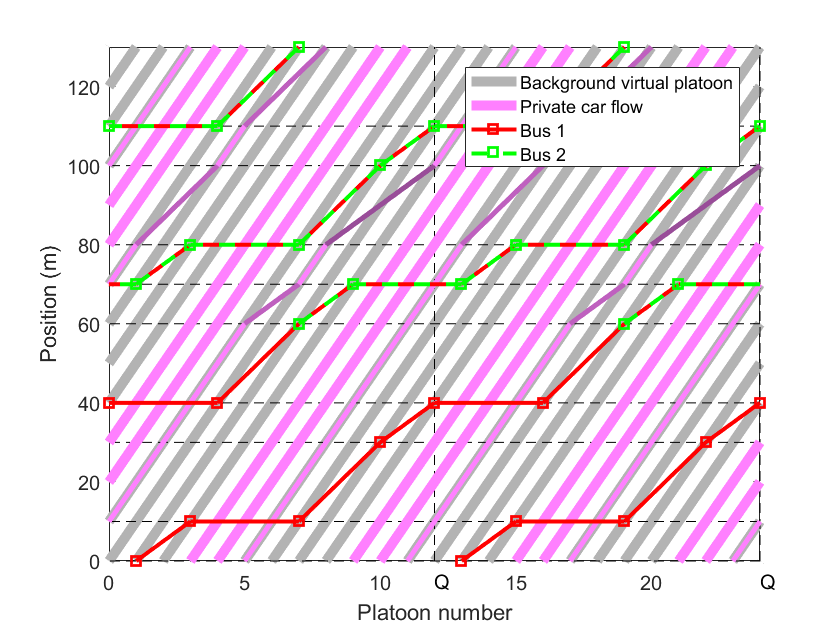}}\hspace{1em}
	\subfloat[][Under private car priority]{\includegraphics[width=0.45\textwidth]{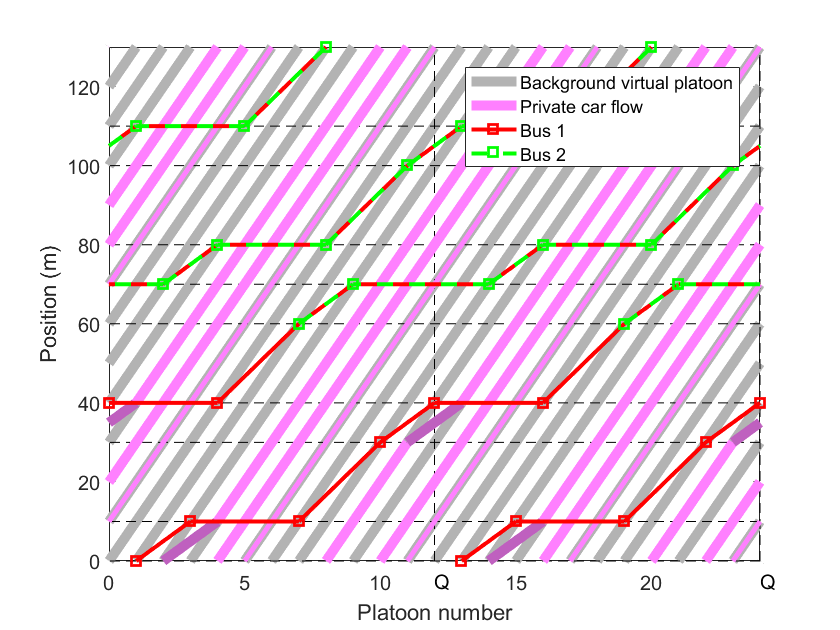}}\\
	\caption[]{Results of \textbf{MILP-O} under medium demand}
	\label{Fig_toy_resu_MILP2}
\end{figure}
\begin{figure}[!ht]
	\centering
	\subfloat[][Under bus priority]{\includegraphics[width=0.45\textwidth]{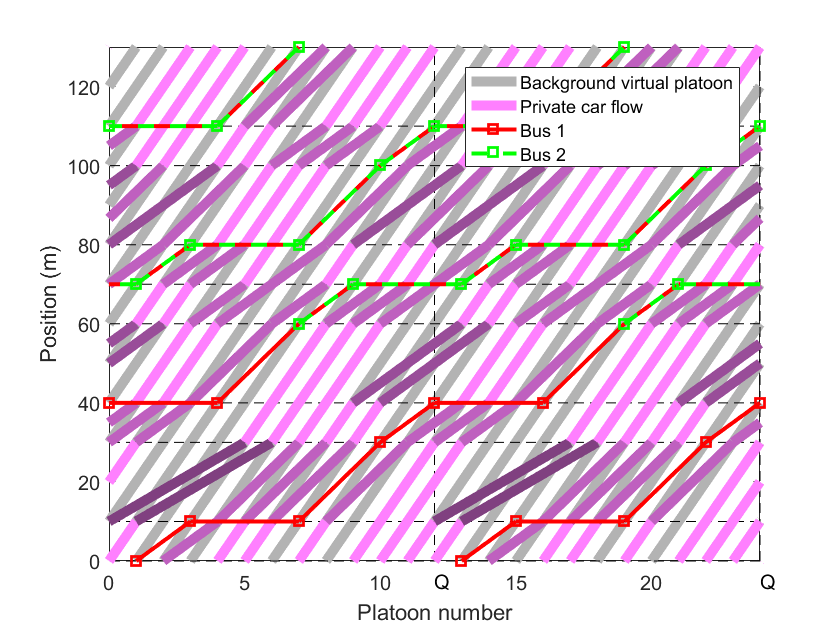}}\hspace{1em}
	\subfloat[][Under private car priority]{\includegraphics[width=0.45\textwidth]{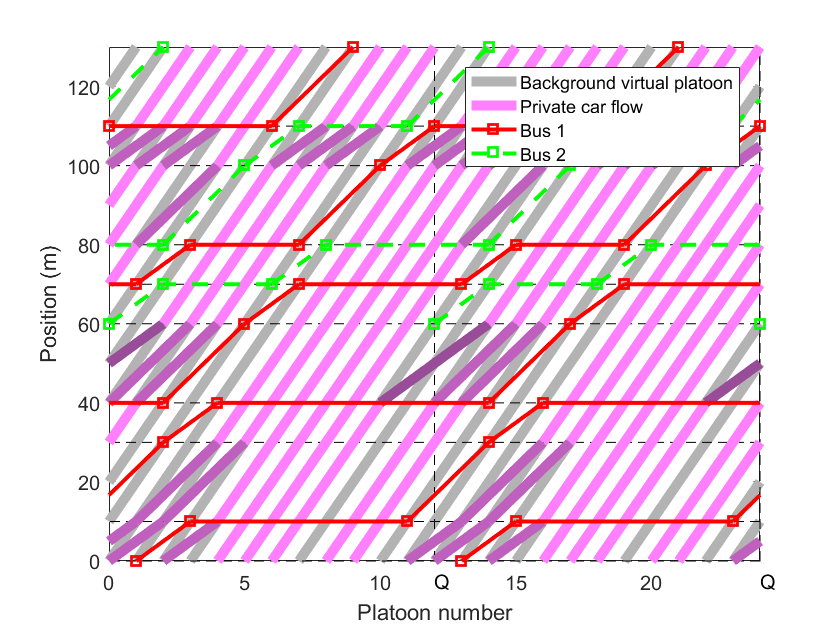}}\\
	\caption[]{Results of \textbf{MILP-O} under high demand}
	\label{Fig_toy_resu_MILP3}
\end{figure}
\begin{figure}[!ht]
	\centering
	\subfloat[][Under bus priority]{\includegraphics[width=0.45\textwidth]{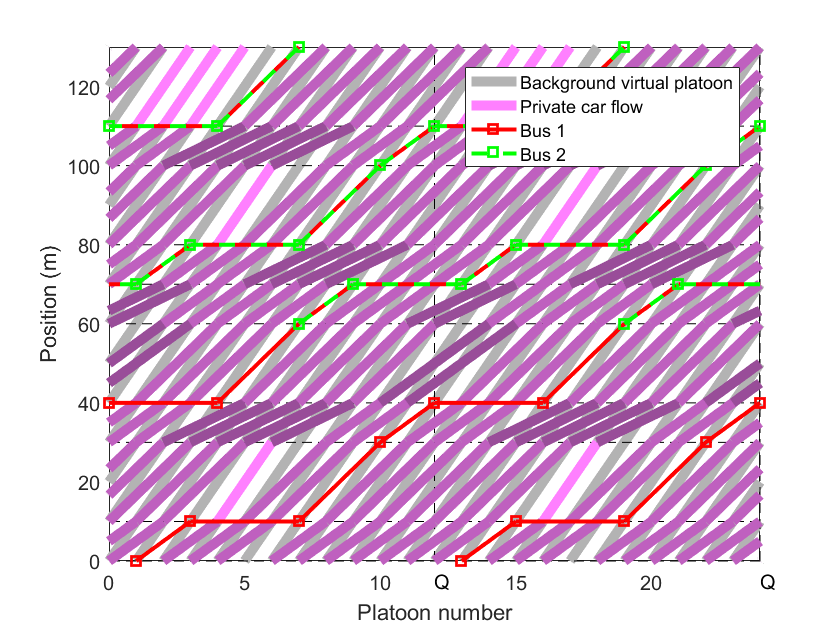}}\hspace{1em}
	\subfloat[][Under private car priority]{\includegraphics[width=0.45\textwidth]{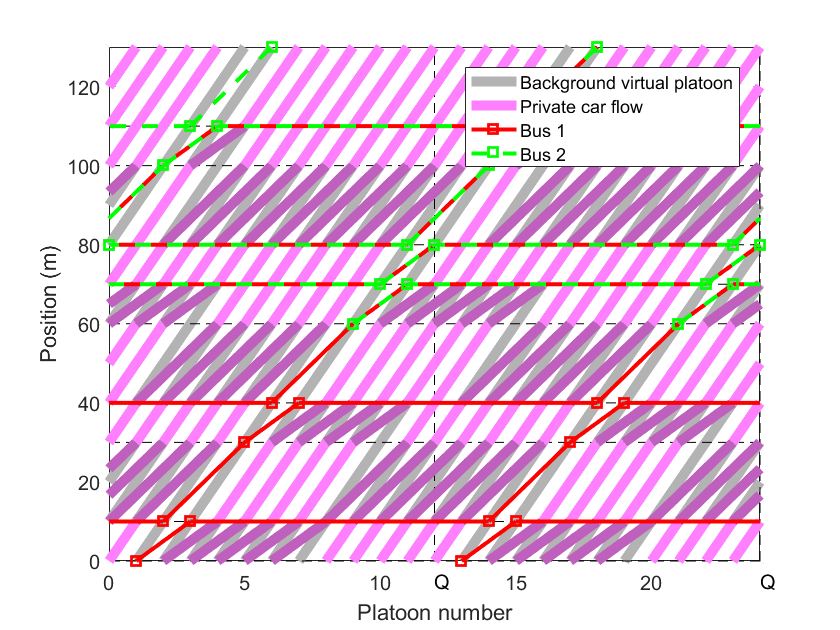}}\\
	\caption[]{Results of \textbf{MILP-O} under maximum admissible demand}
	\label{Fig_toy_resu_MILP4}
\end{figure}

Figs.\ref{Fig_toy_resu_MILP1}-\ref{Fig_toy_resu_MILP4} draw the trajectories of buses and private cars from the results of \textbf{MILP-O} under different demand levels, where (a) and (b) correspond to bus priority and private car priority, respectively. The pink lines denote the flows of private cars, where a darker color implies more delay and a thicker line implies more traffic volume; the red and green lines denote the trajectories of buses 1 and 2, respectively. From Figs.\ref{Fig_toy_resu_MILP1}-\ref{Fig_toy_resu_MILP4}, several additional observations for the MILP solution are as follows:
\begin{itemize}
    \item The travel time of buses on all links is equal to the minimum travel time. That means that the extra travel cost of buses is due solely to the extra dwelling time at stations, supporting the setting of the bus travel time on links as the minimum value when solving the upper-level model, as discussed in Section \ref{upper_subsec}.
    \item Different buses collaborate to travel together for most of the shared links. It is expected as the collaboration could reduce the number of dedicated VPs and then reserve more capacity for private cars, in accordance with the idea of the initial solution generation of the upper-level model.
    \item Under bus priority, more delayed VPs are realized for private cars, especially when the demand is high, while buses are persistently granted to fulfill their trips with the minimum travel time. By the contrary, under private car priority where the efficiency consideration of buses gives way to that of private cars, extra dwelling times of buses are generated. In particular, when the demand is sufficiently high ($\beta_1=1$), the dwelling time of buses  approaches the maximum value, consequently generating an effect of “green wave bands” for private cars as shown in Fig.\ref{Fig_toy_resu_MILP4}(b), remarkably promoting the efficiency of private cars. 
    \item In Fig.\ref{Fig_toy_resu_MILP4} where $\beta_1=1$, there are some links in which all background VPs are either designed for dedicated VPs or regular VPs with full traffic volume of private cars, as implied by the definition of maximum admissible traffic. 
\end{itemize}

\subsubsection{Results of heuristic algorithm} \label{heu_toy_subsubsec}
To test the effectiveness of the bilevel solving method and the proposed heuristic algorithm, comparisons are conducted between the results of the optimal solution and the heuristic solution. Table \ref{toy_heu_table} shows the results under different demand levels, where $O^{opt}$, $O_{ap}^{heu}$, and $O_{ex}^{heu}$ denote the total cost of buses and private cars under the optimal solution of \textbf{MILP-O}, heuristic solution with LP approximations and heuristic solution without LP approximations, respectively; $\Delta O^{heu}$ is the gap between $O^{opt}$ and $O_{ex}^{heu}$. For the heuristic algorithms, the cost of the solution with approximation techniques is identical to the exact solution, verifying the reasonability of the approximation techniques. Moreover, the gap between the optimal solution and the heuristic solution is also tiny, with the maximum value of \textbf{$1.46\%$}; it validates that the heuristic method can provide a sub-optimal solution for \textbf{MILP-O} with a high-quality performance.\par 
\begin{table}[]
    \centering
    \caption{Results obtained by the heuristic method}
    \label{toy_heu_table}
    \begin{tabular}{c|c c c}
    \hline
    Case & 1 & 2 & 3\\
    \hline
    $\beta_2$ & 0.2 & 0.5 & 0.8 \\
    $O^{opt}$ & 10791 & 11152 & 11792\\
    $O_{ex}^{heu}$ & 10948 & 11152 & 11792\\  
    $O_{ap}^{heu}$ & 10948 & 11152 & 11792\\ 
    $\Delta O^{heu}$ & 1.46\% & 0.00\% & 0.00\% \\
    \hline
    \end{tabular}
\end{table}

\begin{figure}[!ht]
	\centering
	\subfloat[][Approximate solution]{\includegraphics[width=0.45\textwidth]{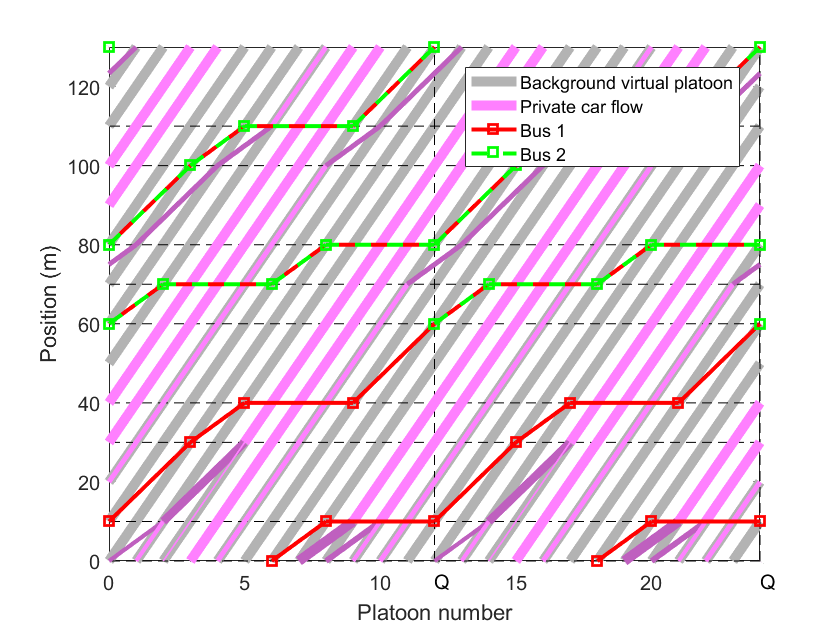}}\hspace{1em}
	\subfloat[][Exact solution]{\includegraphics[width=0.45\textwidth]{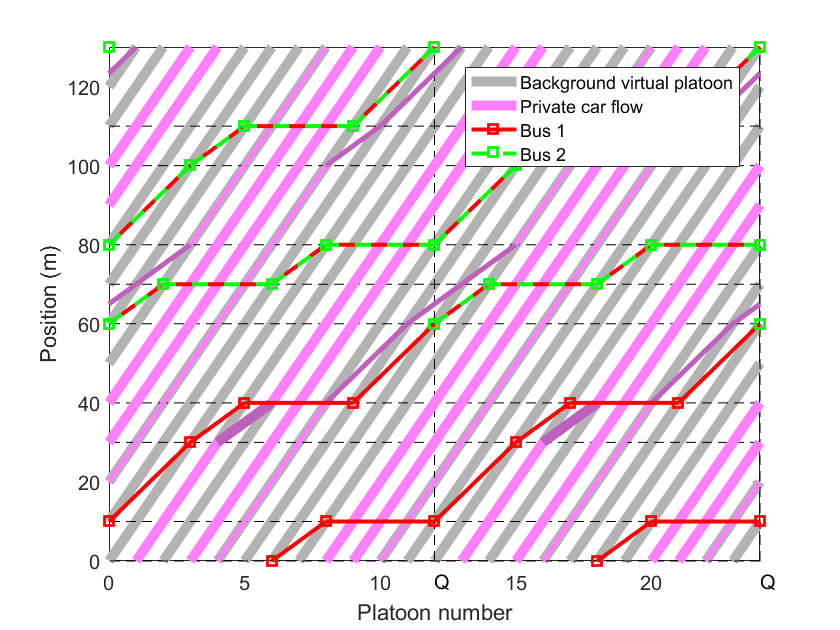}}\\
	\subfloat[][Iteration curve]{\includegraphics[width=0.45\textwidth]{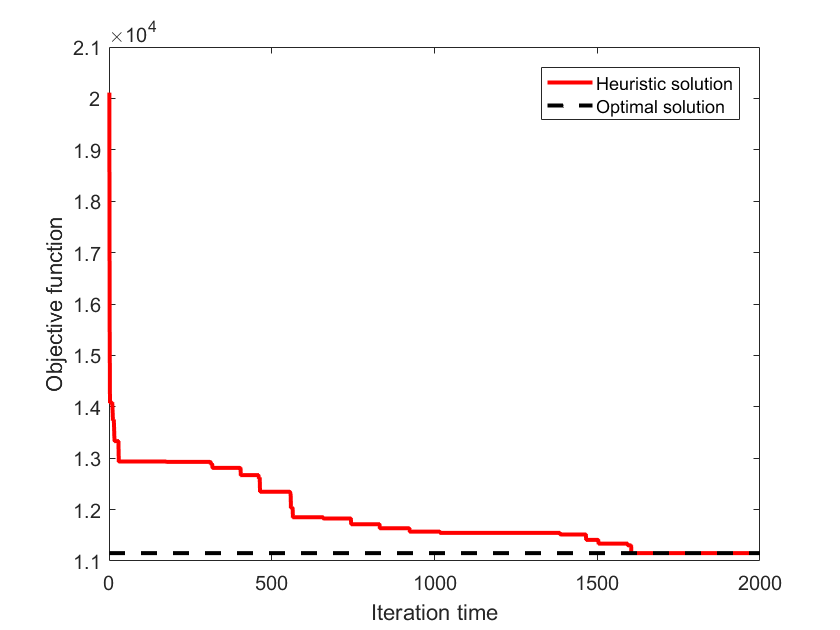}}
	\caption[]{Results of heuristic algorithm under medium demand}
	\label{Fig_toy_resu_appro}
\end{figure}

\begin{figure}[!ht]
	\centering
	\subfloat[][Under low demand]{\includegraphics[width=0.45\textwidth]{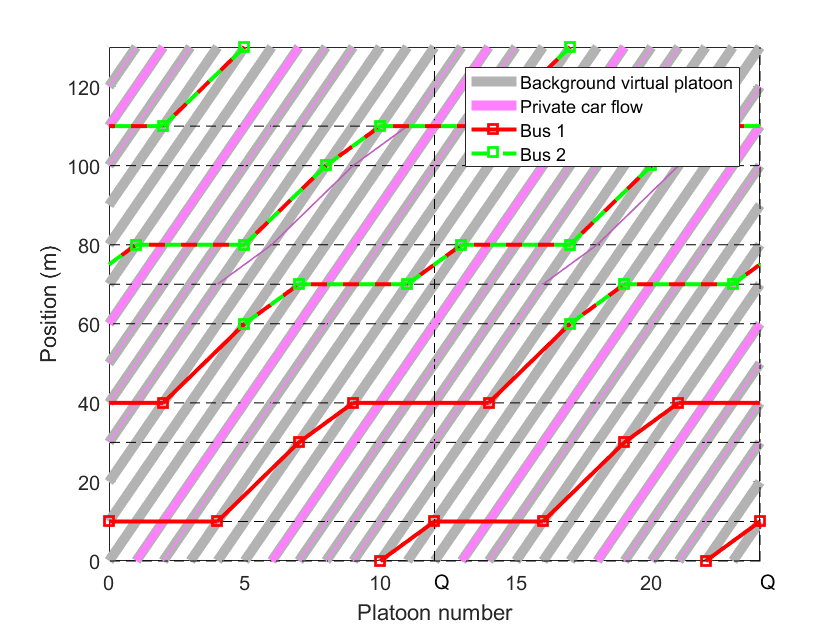}}\hspace{1em}
	\subfloat[][Under high demand]{\includegraphics[width=0.45\textwidth]{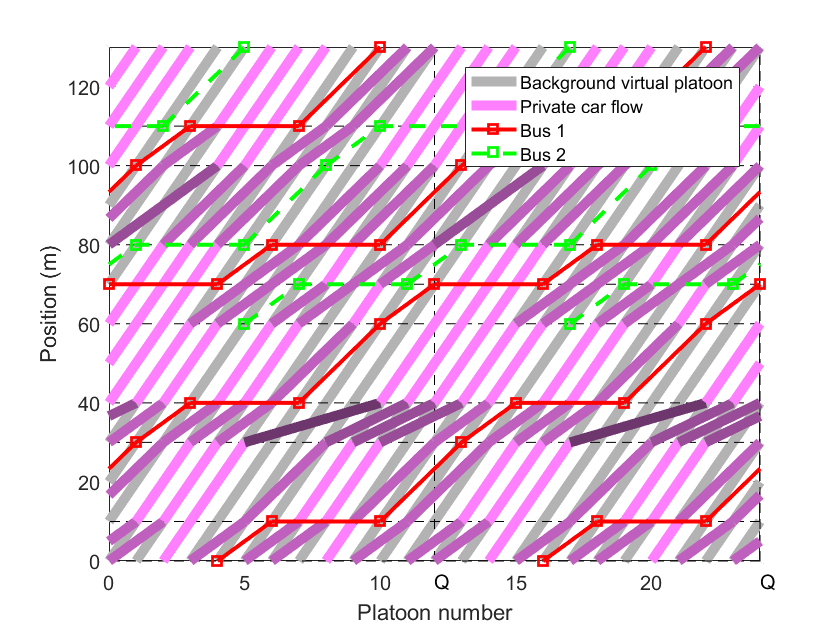}}\\
	\caption[]{Results of heuristic algorithm under different demands}
	\label{Fig_toy_resu_heu}
\end{figure}

Taking the case of medium demand as an example, Fig.\ref{Fig_toy_resu_appro} shows the results of heuristic algorithms where (a) and (b) draw the solutions obtained with LP approximations and without approximations, respectively. Comparing Figs.\ref{Fig_toy_resu_appro} (a) and (b), we observe that only limited regular VPs violating the constraints of FIFO and the indivisibility of VPs have been realized, so their impacts on the results is slight, consistent with the results listed in Table \ref{toy_heu_table}. Fig.\ref{Fig_toy_resu_appro}(c) draws the curve of the objective value in the iteration process. It illustrates that the solution has been improved continuously by the bilevel model, and the optimal solution is obtained after approximately 1600 iterations.\par  
Fig.\ref{Fig_toy_resu_heu} shows the results of the heuristic algorithm under low and high demand cases. It is worthy of noting that in some cases (e.g., high demand), although the solution provided by the heuristic algorithm achieves the minimum objective value as the optimal solution of \textbf{MILP-O} does, different itinerary plans of buses and private cars are observed, implying that the optimal solution is non-unique.\par

\subsection{Real-world network implementation} \label{real_net_subsec}

\subsubsection{Experimental setups} \label{setup_real_subsubsec}

\noindent\underline{\textit{Road network settings}}\par

\begin{figure}
    \centering
    \includegraphics[width=0.5\textwidth]{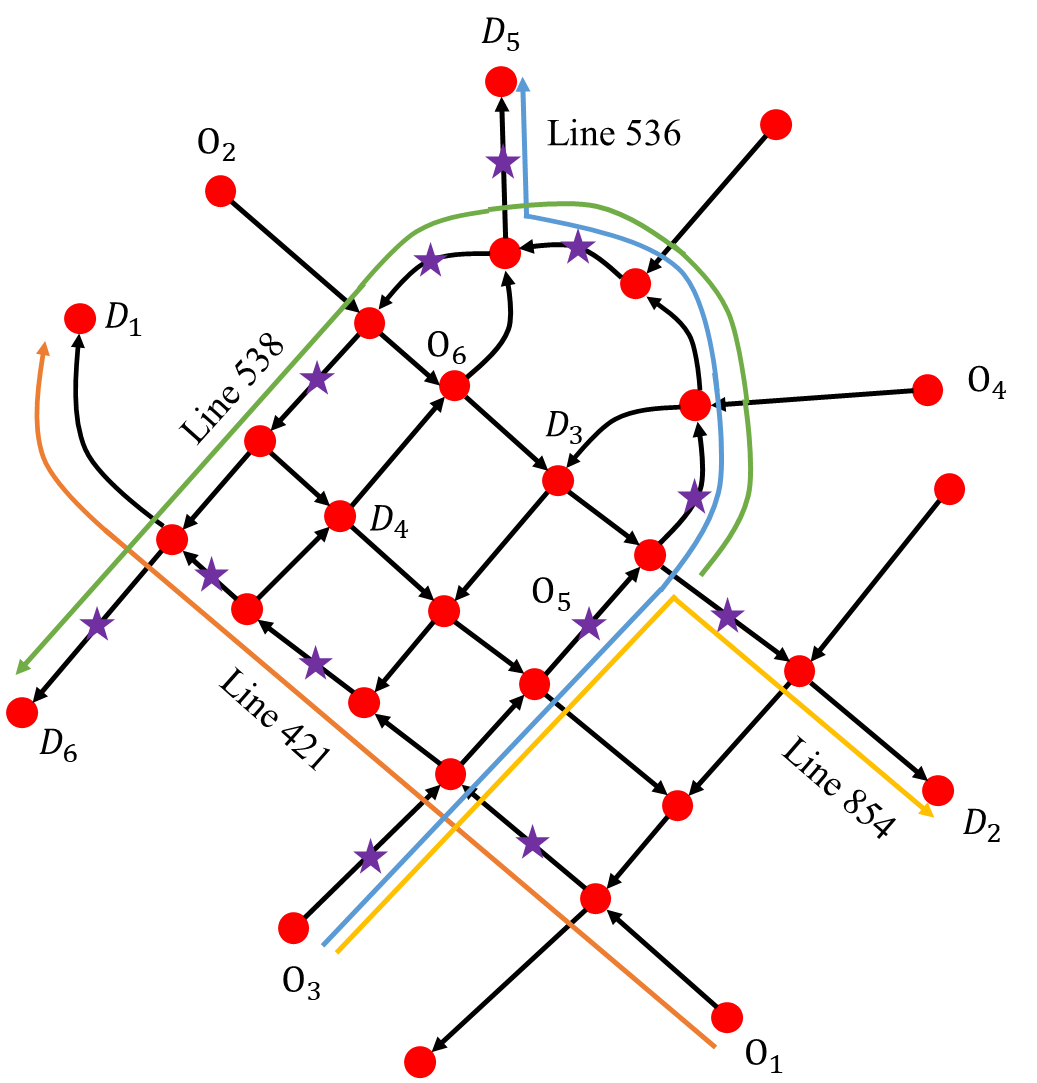}
    \caption{Road network of Wangjing district}
    \label{Fig_real_network}
\end{figure}

Considering a road network from a part of Wangjing district, Beijing, China, for the implementation of the the RC-H scheme, the network is re-organized with one-way roads, as shown in Fig.\ref{Fig_real_network}. The network contains a total of 18 intersections, 5 entrances and 6 exits. The number of lanes is set to $2$ of one way for each link; thus, for the road with buses traveling, there is one mixed-traffic lane and one regular lane. Based on the real-world situation, 4 bus lines and 12 bus stations are set as shown in Fig.\ref{Fig_real_network}, where the routes of buses are denoted by colorful lines and bus stations are denoted by purple stars. The detailed settings of the bus lines can be found in Table \ref{real_bus_table}.\par

\begin{table}[]
    \centering
    \caption{Settings for bus lines}
    \label{real_bus_table}
    \begin{tabular}{c|c c c}
    \hline
    Number & Name & Route & Station\\
    \hline
    1 & Line 421 & $1\rightarrow15$ & 4;9;11 \\
    2 & Line 854 & $7\rightarrow30$ & 6;21;28 \\
    3 & Line 536 & $7\rightarrow40$ & 6;21;32;37;39 \\
    4 & Line 538 & $27\rightarrow14$ & 32;37;41;22;13 \\
    \hline
    \end{tabular}
\end{table}

\begin{figure}
    \centering
    \includegraphics[width=0.4\textwidth]{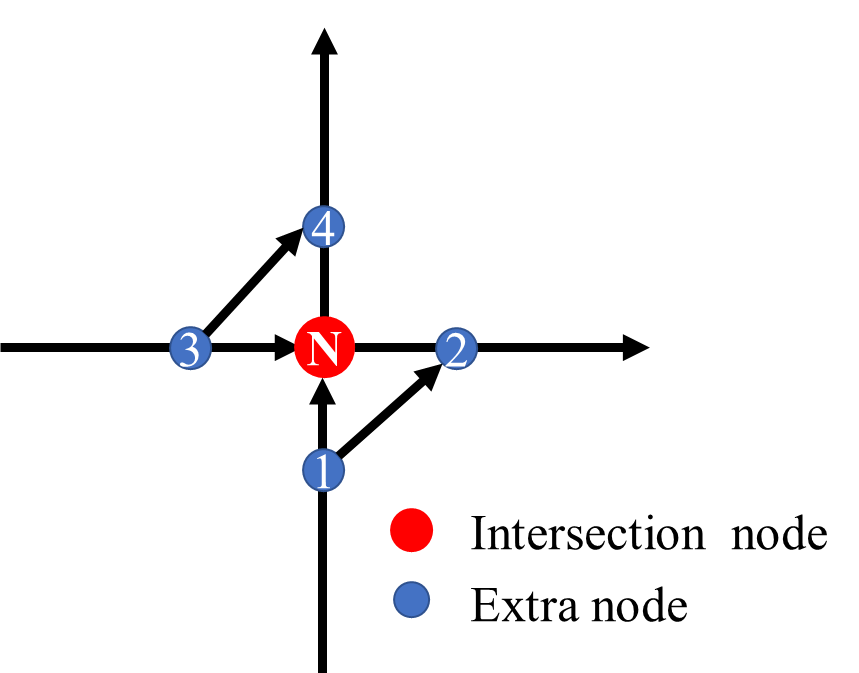}
    \caption{Extra nodes and links within intersections}
    \label{Fig_extranode_inters}
\end{figure}

\begin{figure}
    \centering
    \includegraphics[width=0.6\textwidth]{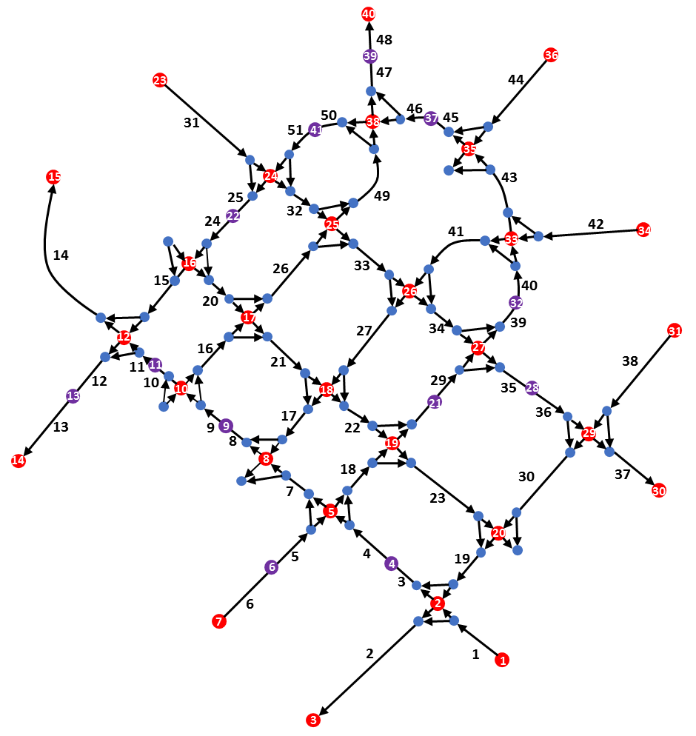}
    \caption{Illustration for the overall network topology}
    \label{Fig_wholenetwork}
\end{figure}

According to the road design required by the RC scheme, the additional nodes and links at the intersections are introduced for the traffic of right-turning and left-turning movements. For each intersection, 4 addtional nodes are placed  and then 4 links are formed for different movements as shown in  Fig.\ref{Fig_extranode_inters}. Consequently, along with the virtual nodes and links for entrances, the whole network is composed of 101 nodes and 129 links, as shown in Fig.\ref{Fig_wholenetwork}.\par 

\vspace{0.5em}
\noindent\underline{\textit{Traffic demand settings}}\par
As the real traffic demand is limited and the network is transferred to one-way roads, we design various demand patterns in the tests, which can fully verify the performance of the RC-H scheme. Three types of O-D pairs are considered: (1) from an entrance to an exit; (2) from an entrance to a junction; (3) from a junction to an exit; the junction is defined as the node inside the network with traffic generation or absorption, e.g., parking facilities and buildings. Two different O-D pairs are set for each type and then a total of 6 O-D pairs for the demand of private cars is obtained as shown in Fig.\ref{Fig_real_network} and Table \ref{real_OD_table}. Moreover, to improve the computational efficiency, for each O-D pair, only the three shortest paths are candidates for the traffic assignment. In fact, as observed in the following test results, even the third shortest path is not utilized, validating that the shrinkage of the path space is reasonable and its impact on the traffic performance is negligible.\par
\begin{table}[]
    \centering
    \caption{Settings for the O-D pairs of private cars}
    \label{real_OD_table}
    \begin{tabular}{c|c c c}
    \hline
    Number & Type & Origin & Destination  \\
    \hline
    1 & 1 & 1 & 15\\
    2 & 1 & 23 & 30\\
    3 & 2 & 7 & 26\\
    4 & 2 & 34 & 17\\
    5 & 3 & 21 & 40\\
    6 & 3 & 25 & 14\\
    \hline
    \end{tabular}
\end{table}

\vspace{0.5em}
\noindent\underline{\textit{RC-H scheme settings}}\par
Prior to the design of the RC-H scheme, an RC scheme should be generated that provides the background rhythm and VPs for the RC-H scheme. The design method of the RC scheme is borrowed from  Lin et al. (2021), except that the objective function is the total personal delay of both buses and private cars. Therefore, an MILP problem is formulated to optimize the design of the RC scheme as follows:
\begin{align}
& \min_{\boldsymbol{\tau},\boldsymbol{\alpha},\boldsymbol{r},\boldsymbol{\theta}}\sum_{ij\in\mathcal{A}}(Q\times s_a\times r_a^{ij}+\sum_{p\in\mathcal{P}}\gamma_p\times r_p^{ij})
\end{align}
\begin{align}
&\mathrm{s.t.} \nonumber\\
\begin{split}
    & \tau_j-\tau_i+\alpha_a^{ij}T\geq\underline{t}_a^{ij}\\
    & \tau_j-\tau_i+\alpha_a^{ij}T<\underline{t}_a^{ij}+T\\
    & r_a^{ij}=\tau_j-\tau_i+\alpha_a^{ij}T-\underline{t}_a^{ij}\\
    & \tau_i\in[0,T)\\
    & \alpha_a^{ij}\in\mathbb{Z}
\end{split} & \forall ij\in\mathcal{A}\\
\begin{split}
    & \tau_j-\tau_i+\alpha_p^{ij}T\geq\underline{t}_p^{ij}\\
    & \tau_j-\tau_i+\alpha_p^{ij}T<\underline{t}_p^{ij}+T\\
    & r_p^{ij}=\tau_j-\tau_i+\alpha_p^{ij}T-\underline{t}_p^{ij}\\
    & \alpha_p^{ij}\in\mathbb{Z}
\end{split} & \forall ij\in\mathcal{A}, p\in\mathcal{P}_{ij}\\
\begin{split}
    & \tau_{i_1}-\tau_{i_2}+T\theta_{i_1,i_2}=\frac{T}{2}\\
    & \theta_{i_1,i_2}\in \{0,1\}
\end{split} & \forall (i_1,i_2)\in\mathcal{C}
\end{align}
where $r_a^{ij}$ and $r_p^{ij}$ denote the delay of VPs and bus $p$ on link $(i,j)$, respectively, $\underline{t}_a^{ij}$ and $\underline{t}_p^{ij}$ denote the free-flow travel time of VPs and bus $p$ on link $(i,j)$, respectively, and $\mathcal{C}$ denotes the set of conflict points. By solving the MILP problem, a background RC scheme can be acquired. Note that the RC scheme is a preliminary design with some simplifications, such as a preset cycle length and an equal platoon size on all links. As the focus of this study is the control scheme for handling heterogeneous traffic, the elaborate design of RC scheme is not investigated in depth, and readers can refer to  Lin et al. (2021) for the more details. \par

The setting of the parameters related to RC-H scheme is identical to the toy example, including $T=10s$, $H=120s$, $Q=12$, $s_a=4$, $s_b=2$, $\gamma_p=20$. Additionally, the minimum travel time is calculated by the free-flow speed of vehicles, which is set as $15m/s$ and $10m/s$ for private cars and buses, respectively.

\subsubsection{Results of RC-H scheme} \label{results_real_subsubsec}
First, we calculate the maximum admissible traffic for the network. For a single RC-H cycle, the maximum traffic of a regular lane can be calculated by $Q\times s_a$, while the maximum traffic of a mixed-traffic lane is identical to the toy example, i.e., $(Q-2)\times s_a$. Thus, the total capacity of a road with traveling buses is $(2Q-2)\times s_a$, and that without buses is $2Q\times s_a$. The maximum admissible traffic is defined as the upper bound of the traffic demand, i.e., $\beta_2=1$.\par
To test the impact of traffic demand, three different demand levels are considered as $\beta_2=0.2$ (low), $\beta_2=0.5$ (medium) and $\beta_2=0.8$ (high) for all O-D pairs. Apart from the balanced demand scenarios, imbalanced demands are also considered to test the performance under different demand patterns. Three cases are tested where the demand of one O-D pair type is double of the other two types. The settings of the traffic demand are shown in Table \ref{real_results_table}.\par
\begin{table}[]
    \centering
    \caption{Results of RC-H scheme}
    \label{real_results_table}
    \begin{tabular}{c|c c c c c c}
    \hline
   Pattern & \multicolumn{3}{c}{Balanced} & \multicolumn{3}{c}{Imbalanced} \\
   \hline
    Case & 1 & 2 & 3 & 4 & 5 & 6 \\
    \hline
    $\beta_2$ (Type 1) & 0.2 & 0.5 & 0.8 & 1.2 & 0.6 & 0.6\\
    $\beta_2$ (Type 2) & 0.2 & 0.5 & 0.8 & 0.6 & 1.2 & 0.6\\
    $\beta_2$ (Type 3) & 0.2 & 0.5 & 0.8 & 0.6 & 0.6 & 1.2\\
    $O_a^m$ & 7700 & 19250 & 30800 & 31126 & 31020 & 30254\\ 
    $O_a^{heu}$ & 7700 & 19250 & 31395 & 31414 & 31960 & 31316\\
    $\Delta O_a$ & 0.00\% & 0.00\% & 1.93\% & 0.92\% & 3.03\% & 3.51\%\\
    $O_b^m$ & 36180 & 36180 & 36180 & 36180 & 36180 & 36180\\
    $O_b^{heu}$ & 37180 & 37380 & 37580 & 37580 & 36980 & 37980\\
    $\Delta O_b$ & 2.76\% & 3.32\% & 3.87\% & 3.87\% & 2.21\% & 4.98\% \\
    $O^m$ & 43880 & 55430 & 66980 & 67306 & 67200 & 66434\\
    $O^{heu}$ & 44880 & 56630 & 68975 & 68994 & 68940 & 69296 \\
    $\Delta O$ & 2.28\%	& 2.16\% & 2.98\% & 2.51\% & 2.59\%	& 4.31\% \\
    \hline
    \end{tabular}
\end{table}
Heuristic method is applied in solving the optimization of the RC-H scheme as the MILP is insolvable for the large-scale problem, and the results are shown in Table \ref{real_results_table}. $O_a^m$ and $O_a^{heu}$ denote the free-flow total travel cost and the cost under the heuristic solution of private cars, respectively, and $\Delta O_a$ denotes the gap between $O_a^m$ and $O_a^{heu}$; $O_b^m$, $O_b^{heu}$ and $\Delta O_b$ have similar meanings but for buses, and $O^m$, $O^{heu}$ and $\Delta O$ are for the personal cost of both private cars and buses. It is observed that when the demand is not high, no extra delay of private cars is generated, implying that all private cars can pass through the network with the background rhythm and the interference from buses can be dodged by the RC-H scheme. Moreover, even when the demand is high (i.e., $\beta_2=0.8$), the gap to the minimum travel cost is only $1.93\%$, which embodies the advantage of RC-H in accommodating high demand traffic. Furthermore, although extra travel costs of buses are observed, the gaps are below $5\%$ in all cases, and the extra personal travel costs are also limited as implied by the results of $\Delta O$.\par
\begin{figure}[!ht]
	\centering
	\subfloat[][Bus 3]{\includegraphics[width=0.45\textwidth]{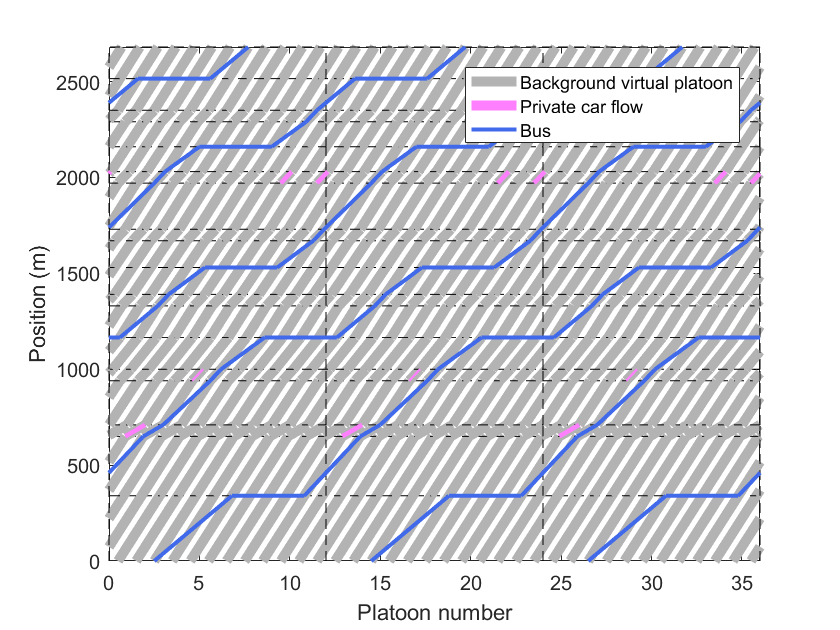}}\hspace{1em}
	\subfloat[][Bus 4]{\includegraphics[width=0.45\textwidth]{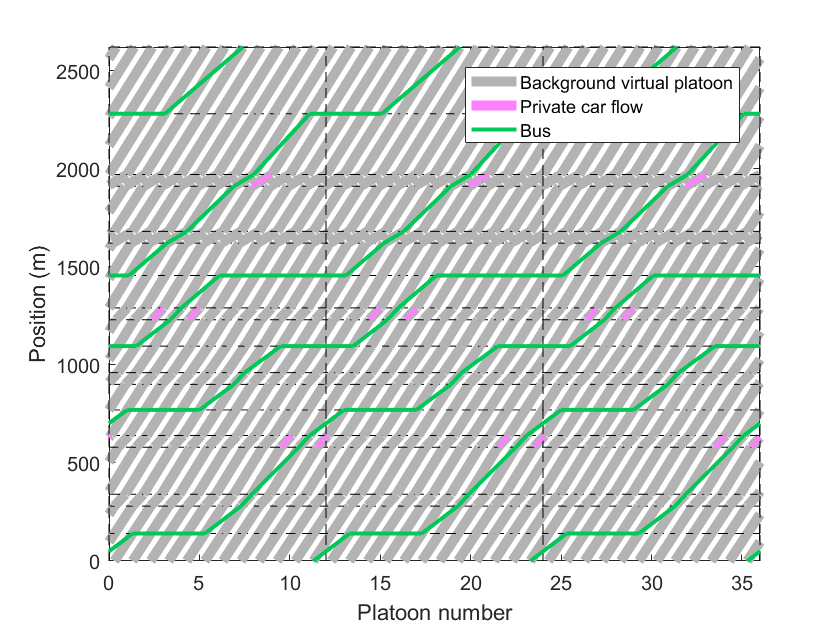}}\\
	\caption[]{Results on the bus itineraries under medium demand}
	\label{Fig_real_resu_case2B}
\end{figure}

\begin{figure}[!ht]
	\centering
	\subfloat[][O-D pair 3 on regular lanes]{\includegraphics[width=0.45\textwidth]{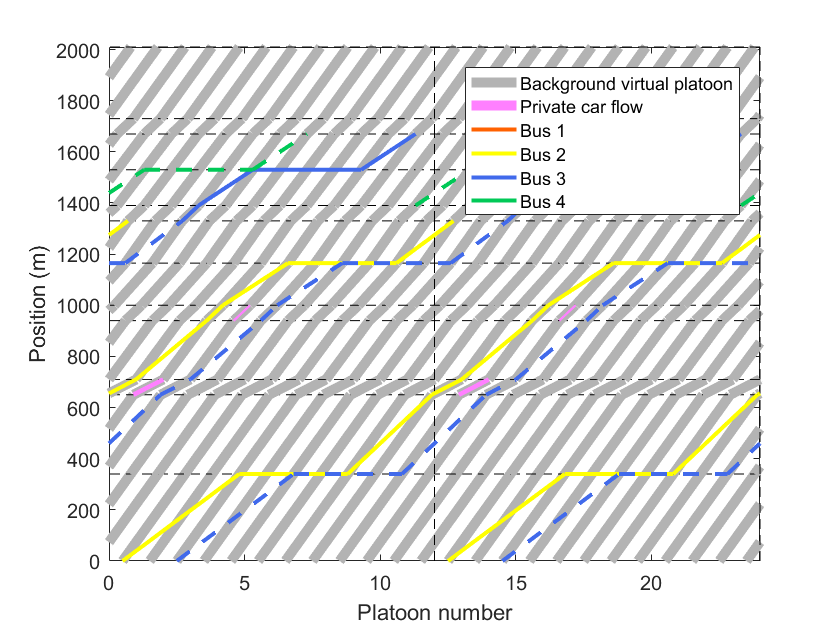}}\hspace{1em}
	\subfloat[][O-D pair 3 on mixed lanes]{\includegraphics[width=0.45\textwidth]{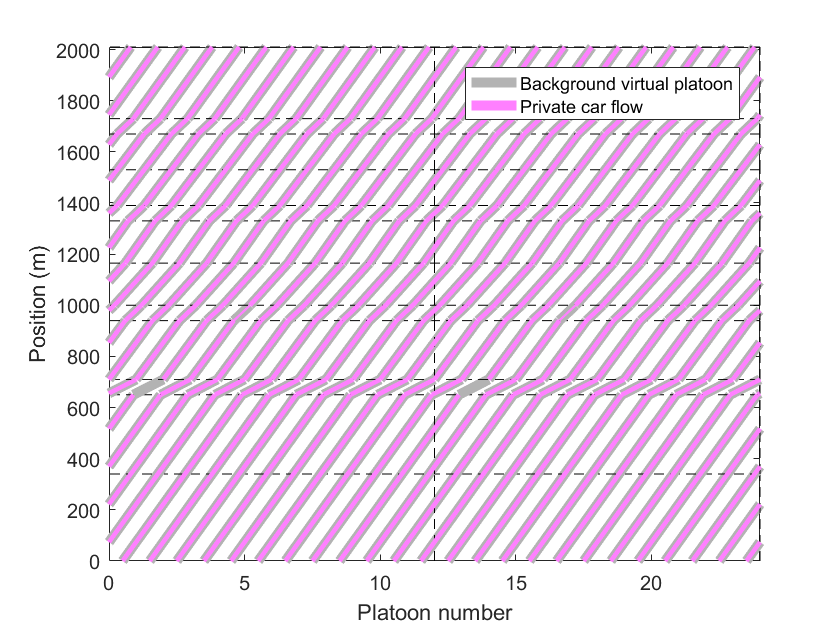}}\\
	\caption[]{Results on private car movements under medium demand}
	\label{Fig_real_resu_case2A}
\end{figure}

 
Figs.\ref{Fig_real_resu_case2B} and \ref{Fig_real_resu_case2A} show the results of the RC-H scheme under the medium demand (i.e., case 2), where the itineraries of buses 3 and bus 4, and the traffic of O-D pair 3 are shown as examples. Figs.\ref{Fig_real_resu_case2B}(a) and (b) show the bus itineraries of buses 3 and 4, respectively, along with the traffic of private cars on the bus paths to illustrate the traffic assignment results in the mixed lanes. Fig.\ref{Fig_real_resu_case2A} shows the traffic of O-D pair 3 on its shortest path, where the left column is the traffic in the mixed-traffic lane along with the bus trajectories on the related links to illustrate the interaction with buses, and the right column is the traffic on the regular lane. According to the results, most of the private cars are assigned to the regular lanes as no extra delay will be generated in the regular lanes; moreover, the traffic assigned to the mixed-traffic lane is either following the VPs with no generated extra delay or yielding to the capacity constraint of the corresponding regular lane. As the traffic demand is not high, the regular lanes accommodate most of the traffic, and the private cars assigned to the mixed-traffic lanes also dodge the interference from buses. Therefore, private cars follow the background VPs in both lanes without extra costs.\par 

\begin{figure}[!ht]
	\centering
	\subfloat[][Bus 3]{\includegraphics[width=0.45\textwidth]{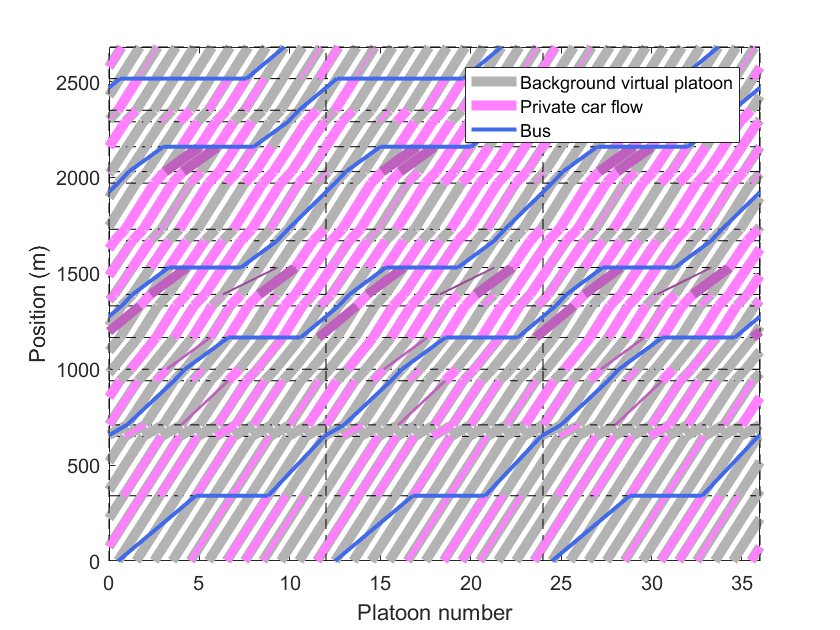}}\hspace{1em}
	\subfloat[][Bus 4]{\includegraphics[width=0.45\textwidth]{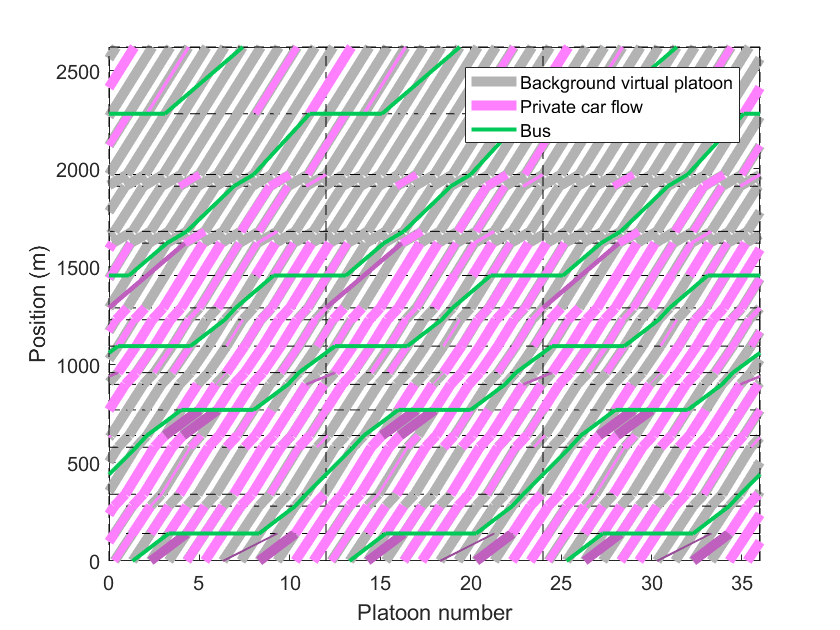}}\\
	\caption[]{Results on the bus itineraries under high demand}
	\label{Fig_real_resu_case3B}
\end{figure}

\begin{figure}[!ht]
	\centering
	\subfloat[][O-D pair 3 on regular lanes]{\includegraphics[width=0.45\textwidth]{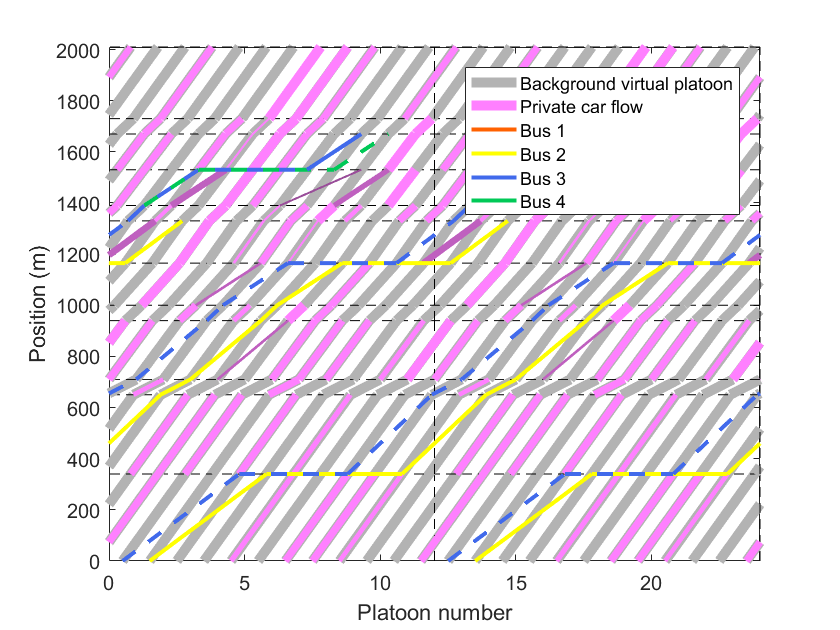}}\hspace{1em}
	\subfloat[][O-D pair 3 on mixed lanes]{\includegraphics[width=0.45\textwidth]{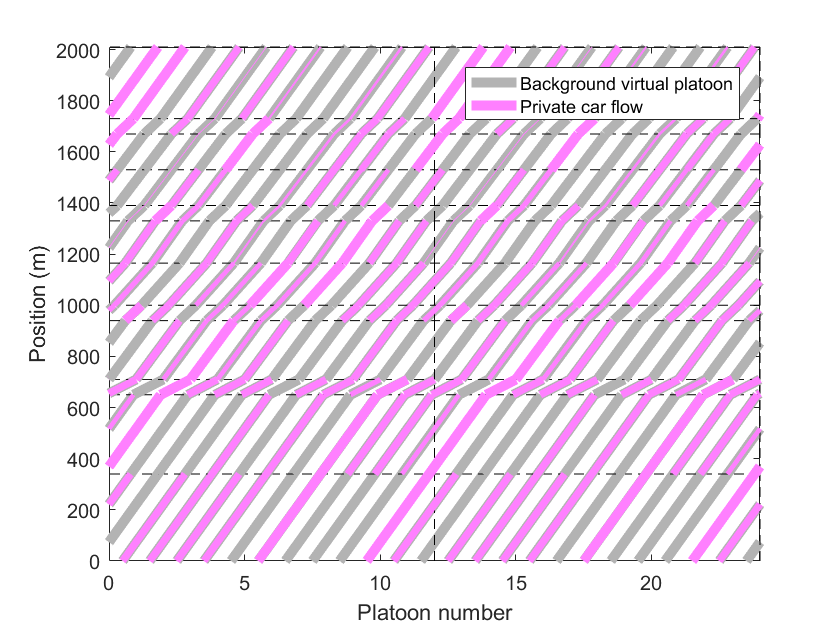}}\\
	\caption[]{Results on private car movements under high demand}
	\label{Fig_real_resu_case3A}
\end{figure}

Figs.\ref{Fig_real_resu_case3B} and \ref{Fig_real_resu_case3A} shows the results under the high demand (i.e., case 3). Compared with case 2, as the demand increases, some delayed VPs are realized for the traffic assignment of private cars, inducing the extra travel cost as shown in Table \ref{real_results_table}. Moreover, according to Fig.\ref{Fig_real_resu_case3A}(a), buses 3 and bus 4 collaborate to travel together on some of the shared links, weakening the interference effect of the buses on the mobility of private cars.\par 

\begin{figure}[!ht]
	\centering
	\subfloat[][Bus 3]{\includegraphics[width=0.45\textwidth]{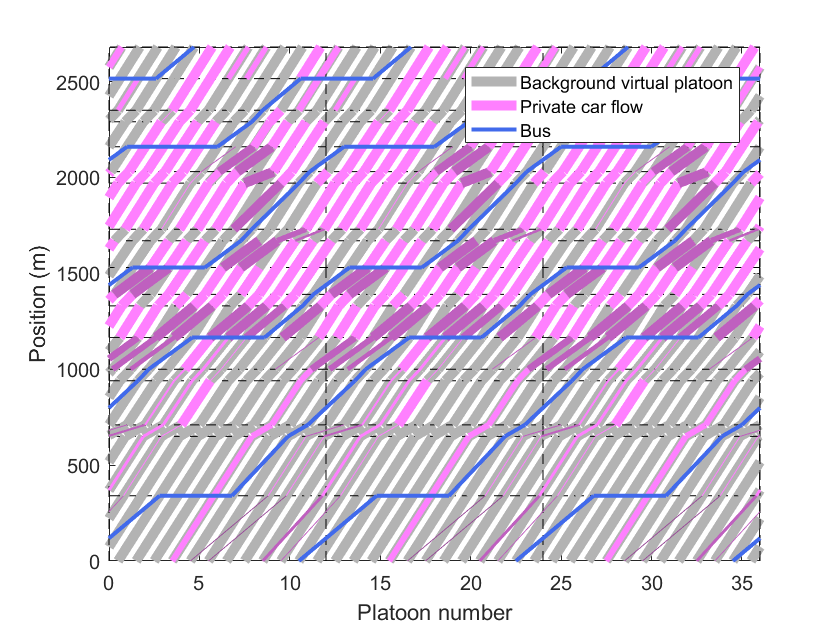}}\hspace{1em}
	\subfloat[][Bus 4]{\includegraphics[width=0.45\textwidth]{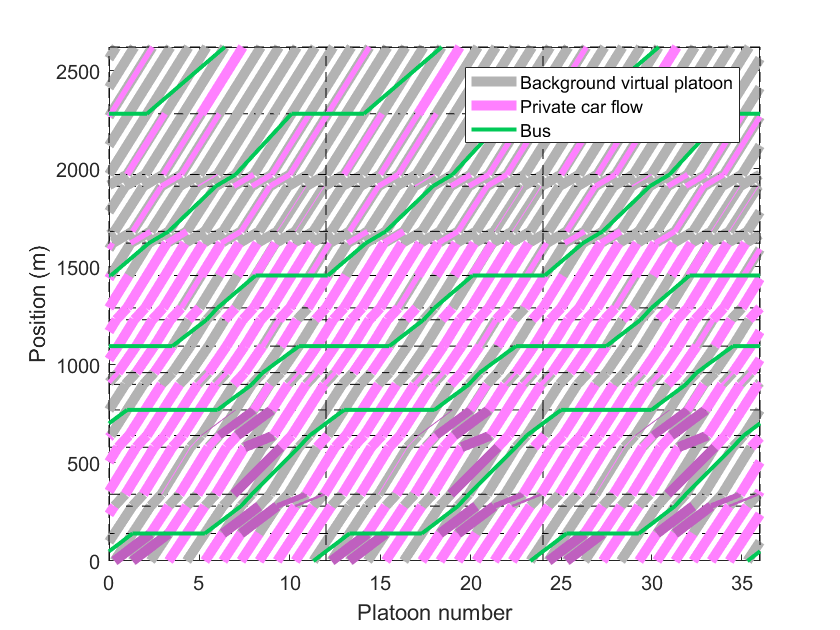}}\\
	\caption[]{Results on the bus itineraries under imbalanced demand}
	\label{Fig_real_resu_case5B}
\end{figure}

\begin{figure}[!ht]
	\centering
	\subfloat[][O-D pair 3 on regular lanes]{\includegraphics[width=0.45\textwidth]{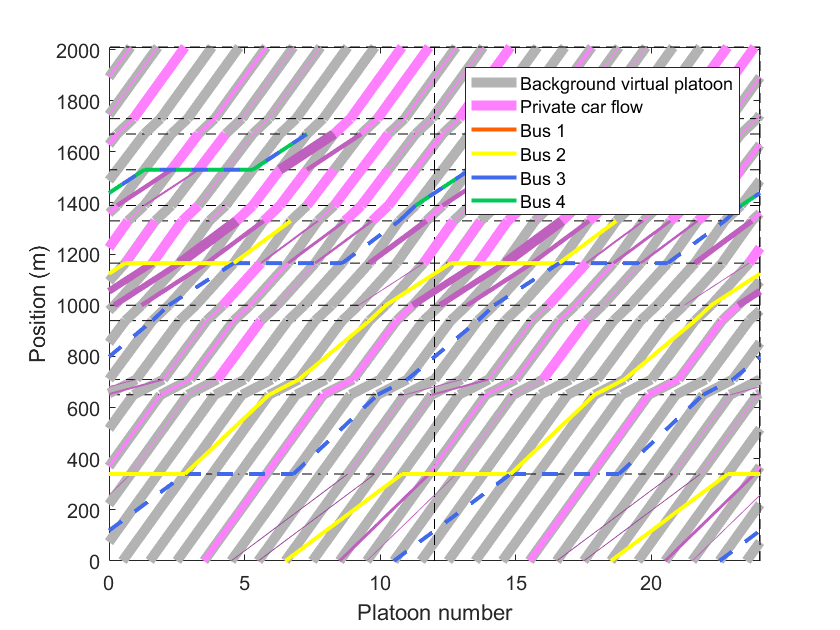}}\hspace{1em}
	\subfloat[][O-D pair 3 on mixed lanes]{\includegraphics[width=0.45\textwidth]{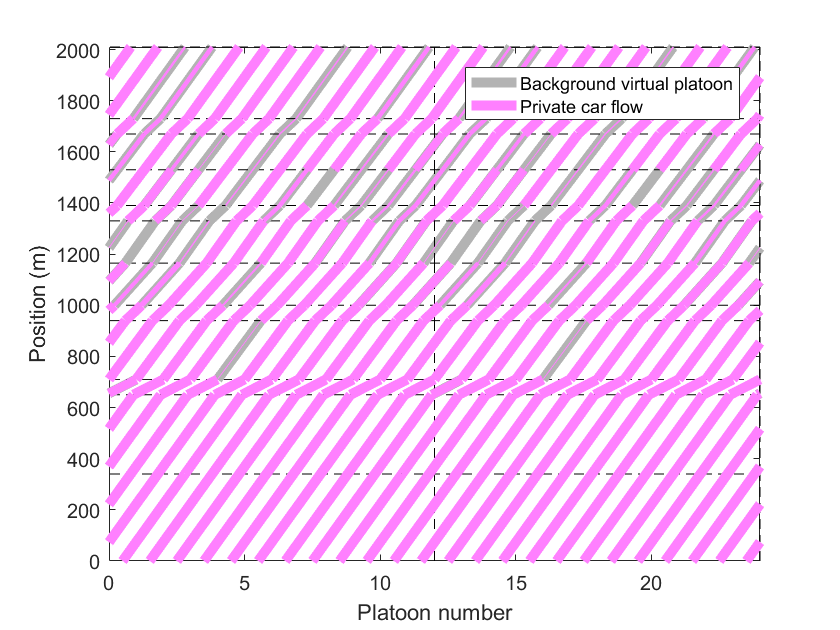}}\\
	\caption[]{Results on private car movements under imbalanced demand}
	\label{Fig_real_resu_case5A}
\end{figure}

In addition, different demand patterns generate quite distinct results as observed from the results presented in Table \ref{real_results_table}. It is observed that the private cars of different O-D pairs are affected by buses to different extents. Specifically, O-D pair 1 interacts with only bus 1 while O-D pairs 3 and 5 interact with buses 2, 3 and 4 in their trips. Therefore, the private cars of O-D pairs 3 and 5 have less admissible demand and more efficiency damage due to the interference from buses. From Table \ref{real_results_table}, it can be obtained that in cases 5 and 6 where O-D pairs 3 and 5 have double demands, greater extra costs are generated compared to case 4. Figs.\ref{Fig_real_resu_case5B} and \ref{Fig_real_resu_case5A} show the results under the imbalanced demand of case 5. In Fig.\ref{Fig_real_resu_case5B}, there are several links that reach the associated lane capacities since the demand is doubled for O-D pair 3, and then more delayed VPs are realized as shown in Fig.\ref{Fig_real_resu_case5A}(a). Moreover, as the demand level of O-D pair 3 is set to 1.2, which is even larger than the upper bound, some private cars are assigned to other paths with detours due to the capacity constraint, further increasing the vehicle costs.

\subsubsection{Simulation comparisons with traffic signal control} \label{simu_real_subsubsec}
As the traffic demand of private cars is assumed to be stationary, the RC-H scheme describes a stable state of the whole network during a single control cycle. Moreover, the discretization of vehicles is ignored in the traffic assignment as the flow of private cars is regarded as continuous variables. Therefore, to further validate the performance under a more general vehicle arrival scenario, simulation tests of online itinerary planning are conducted. The online planning is based on the results of the RC-H scheme, which provides the rhythm and realized VPs for each link. For comparison purpose, traffic signal control (TSC) strategies with the implementation of dedicated bus lanes (DBLs) and without DBLs are also tested in the simulations.\par 
Simulation of Urban Mobility (SUMO) is chosen as the platform for the simulation experiments, which has been widely applied for the tests of traffic operations and management (Behrisch et al., 2011). For TSC strategies, impartially, the time lengths of different phases are set to identical values as the sizes of VPs in the RC-H scheme are the same for all links; moreover, two phase lengths are tested, i.e., $15s$ and $30s$, with the clearance time of $2s$; thus, two cycle lengths are generated as $T=34s$ and $T=64s$, respectively; along with the implementation of DBLs or not, a total of four TSC strategies are tested. The parameter settings of vehicles under TSC strategies are identical to the RC-H scheme, and the vehicle following mode is set as "Cooperative Adaptive Cruise Control (CACC)" in SUMO. Note that the minimum headway of private cars is set as $1s$, which is essentially equivalent to the headway setting within the VP of the RC-H scheme\footnote{In the RC-H scheme, the vehicle headway within the VP is set as $1s$ and the headway between the intersecting VPs is set as $2s$; as the platoon size is set as $4$, the total time for two intersecting VPs to pass through is calculated as $(3+2+3+2)s$, which is identical to the set of the RC cycle length.}. The arrival of private cars is assumed to be Poisson arrival for both TSC strategies and the RC-H scheme, and the traffic demand level varies from $0.1$ to $1.5$ with an interval of $0.1$. For each scenario, the simulation test is run for $1$ hour with $5$ repetitions to test the stable traffic and eliminate randomness. \par 
\begin{figure}[!ht]
	\centering
	\subfloat[][Travel time of private cars]{\includegraphics[width=0.5\textwidth]{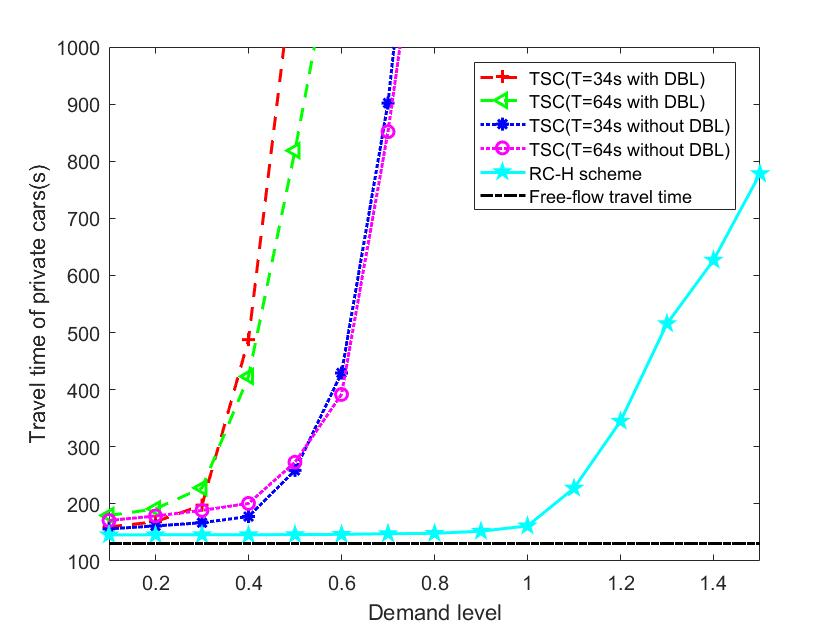}}\\
	\subfloat[][Travel time of buses]{\includegraphics[width=0.5\textwidth]{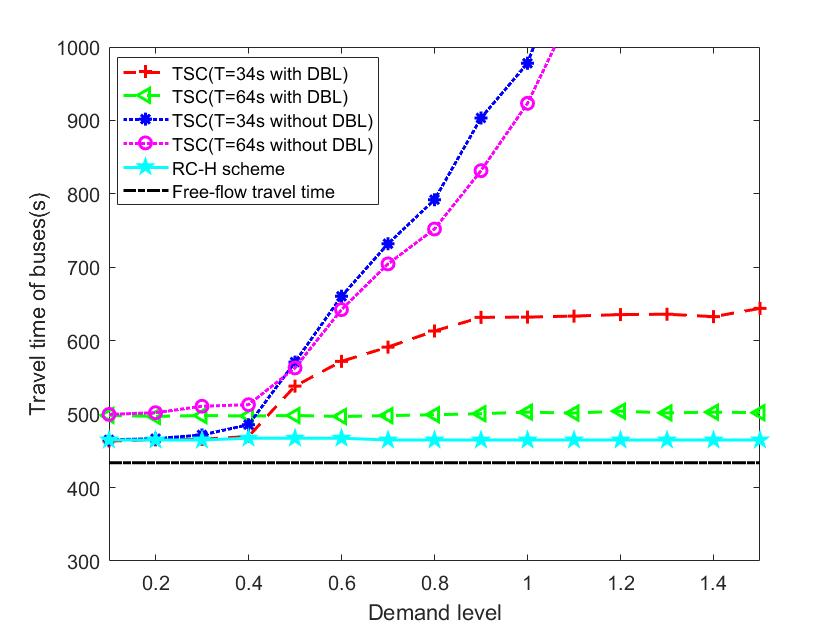}}\\
	\subfloat[][Throughput]{\includegraphics[width=0.5\textwidth]{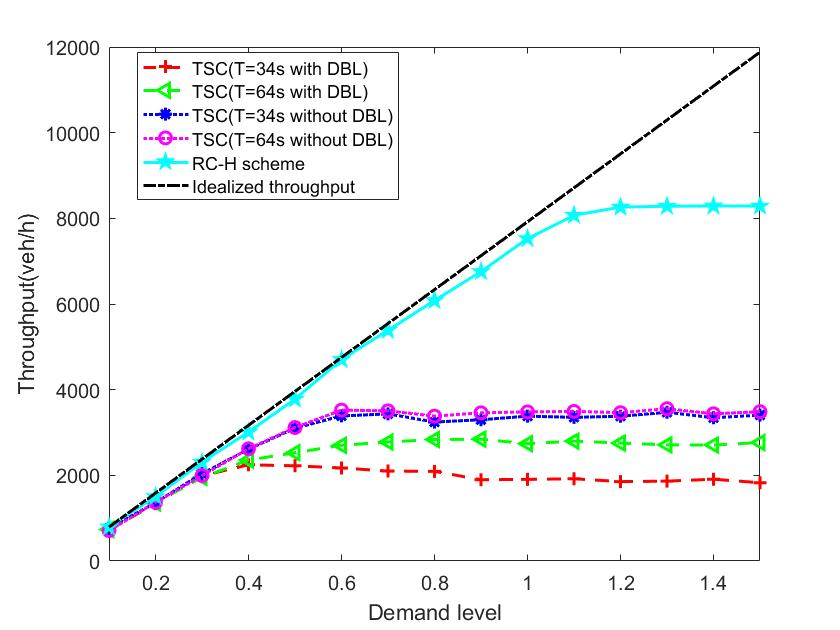}}\\
	\caption[]{Results of simulation experiments}
	\label{Fig_real_simu}
\end{figure}
Fig.\ref{Fig_real_simu} shows the results of simulation experiments, where (a), (b) and (c) draw the average travel time of private cars, the average travel time of buses, and vehicle throughput of the network, respectively. The vehicle throughput is calculated as the number of vehicles that have completed their trips during the simulation time. From Fig.\ref{Fig_real_simu}, several observations are summarized as follows:
\begin{itemize}
    \item For private cars, the travel time holds a value slightly higher than the free-flow travel time when the traffic demand is low, and increases with the traffic demand under all of the tested control strategies. However, the increasing trend under TSC strategies appears to be more aggressive than that for the RC-H scheme as the travel time increases exponentially under TSC strategies when the demand level exceeds $0.4$, while under the RC-H scheme, it begins to increase until the demand level exceeds $0.9$. 
    \item Comparing the four TSC strategies, a longer cycle induces more delays when the demand is low, but the opposite effect appears when the demand is high; that is because a longer cycle of TSC can accommodate more traffic due to the lower time loss, as can also be validated in Fig.\ref{Fig_real_simu}(c). The implementation of DBLs damages the efficiency of private cars to a great extent as the ROW of private cars has been deprived on the lanes with bus traveling. 
    \item By contrast, buses benefit strongly by the implementation of DBLs. Under the TSC without DBLs, the travel time of buses increases with the traffic demand of private cars, due to the interference from the traffic flow in the travel and the vehicle queue at intersections. On the other hand, with DBL implementation, the delay of buses can be alleviated to a bounded value, especially under the TSC with the longer cycle. The result is different for the TSC with the shorter cycle, where the delay of buses still increases with the traffic demand for a certain range; that is because traffic gridlocks are incurred at several intersections, where the buses will be trapped even on the DBLs. Moreover, the delay of buses holds the smallest value under the RC-H scheme among all control strategies, as the priority ROW can be granted for buses along the whole trip, and gridlock will not appear under the centralized control.\par 
    \item Finally, the network throughput first increases with the demand level, and then remains almost non-increasing when the demand level exceeds a certain value. The maximum throughput of RC-H is approximately $8000 veh/h$, while the maximum throughput is below $4000 veh/h$ for all TSC strategies. The result for the throughput further validates the capability of the RC-H scheme to handle a large traffic demand.
\end{itemize}

\section{Concluding remarks} \label{Concluding_sec}
Enhancing the service quality of bus transit is an important topic to promote the attractiveness of buses and leverage the advantages of high-occupancy vehicles in alleviating the urban traffic issues of congestion and pollution. With the assistance of automated driving and V2X communication technologies, vehicles can be controlled accurately and organized efficiently, enabling the whole traffic to be more predictable and reliable. Leveraging the emerging technologies, this study proposes an innovative control scheme for the heterogeneous automated traffic composed of buses and private cars in a network, i.e., RC-H. Inheriting the idea of rhythmic control scheme proposed by  Lin  et  al. (2020), we organize the traffic of buses and private cars in a rhythmic manner to enhance the transit service quality and the network mobility. Specifically, "dedicated virtual platoons" are designed for buses to provide exclusive right-of-ways (ROWs) on their routes and guarantee the schedule adherence through the network, while "regular virtual platoons" are designed for the traffic flow of private cars to improve the traffic efficiency under the interference from buses. By following the paces of virtual platoons, both buses and private cars can pass through the network without any stops and collisions.\par
Along with the design of virtual platoons, the problems of bus itinerary planning and the traffic assignment of private cars are jointly optimized to minimize the total cost of all vehicles, which is formulated as a mixed integer linear program (MILP). By solving the MILP, an optimal design of RC-H scheme can be obtained. A bilevel heuristic solution method is proposed to alleviate the computational burden, where the upper level provides the plan of the bus itinerary by applying the variable neighborhood search (VNS) algorithm, and the lower level solves the traffic assignment problem, which is transformed into a linear program (LP) by approximation techniques. It is validated that the heuristic method provides a sub-optimal solution with a high-quality performance. \par
Numerical experiments are conducted under both a toy example and a real-work network implementation. The results show that the travel cost of buses and private cars are affected by both traffic demand level and efficiency priority. By setting bus priority, no extra delay of buses will be generated under all demand cases and the personal travel cost can be minimized. Furthermore, simulation experiments are tested to verify the performance under a more general scenario and to compare to traffic signal control (TSC) strategies with dedicated bus lane (DBL) implementation. The results confirm the capability of RC-H in alleviating delays of private cars and granting the schedule adherence of buses concurrently. Compared to the TSC strategies, the RC-H scheme has distinctive advantages for handling massive traffic demand, reflected in both vehicle delay reduction and network throughput promotion.\par
Several promising research directions can be suggested. First, the RC-H scheme is optimized based on a given RC scheme obtained from a simplified method in this study. A more elaborate RC scheme or the incorporation of the design of the RC scheme into the RC-H scheme could further improve the performance. Moreover, to solve the optimization efficiently, a bilevel solution model with a heuristic algorithm is proposed, but the iteration time is still restricted due to the computational cost. Thus, design of a more tractable algorithm or application of a high-quality computational tool such as parallel computing is also highly valuable. Furthermore, while buses have some unique features such as fixed routes and preset time schedules, accommodating the general heterogeneous traffic composed of various vehicle types is worthy of further investigation, which could be more complex due to the distinctive travel features of the different vehicles. Finally, a common cycle length is assumed for the overall network to maintain the rhythm harmonization in this study. It is expected that the relaxation of the rhythm union has the potential to handle various demand patterns in the different parts of the network, and to accommodate the heterogeneous traffic of various vehicle types simultaneously.

\section*{Acknowledgements}
The research is supported in part by the Tsinghua-Daimler Joint Research Center for Sustainable Transportation, and the Tsinghua University-Toyota Research Center.
\section*{References}
\noindent Abu-Lebdeh, G. and Benekohal, R. F. (1997). Development of traffic control and queue management procedures for oversaturated arterials. Transportation Research Record, 1603(1):119–127.\par 
\noindent Administration, N. H. T. S. et al. (2013). Preliminary statement of policy concerning automated vehicles. Washington, DC, pages 1–14.\par 
\noindent Baker, R. J., Collura, J., Dale, J. J., Head, L., Hemily, B., Ivanovic, M., Jarzab, J., McCormick, D., Obenberger, J., Smith, L., et al. (2002). An overview of transit signal priority. Technical report.\par 
\noindent Balke, K. N., Dudek, C. L., and Urbanik, T. (2000). Development and evaluation of intelligent bus priority concept. Transportation Research Record, 1727(1):12–19.\par 
\noindent Behrisch, M., Bieker, L., Erdmann, J., and Krajzewicz, D. (2011). Sumo–simulation of urbanmobility: an overview. In Proceedings of SIMUL 2011, The Third International Conference on Advances in System Simulation. ThinkMind.\par 
\noindent Chang, G.-l., Vasudevan, M., and Su, C.-c. (1996). Modelling and evaluation of adaptive bus preemption control with and without automatic vehicle location systems. Transportation Research
Part A: Policy and Practice, 30(4):251–268.\par 
\noindent Chen, X., Li, M., Lin, X., Yin, Y., and He, F. (2020a). Rhythmic control of automated traffic – part i: Concept and properties at isolated intersections. Transportation Science, Submitted, arXiv:2010.04322 [math.OC].\par 
\noindent Chen, X., Lin, X., He, F., and Li, M. (2020b). Modeling and control of automated vehicle access on dedicated bus rapid transit lanes. Transportation Research Part C: Emerging Technologies, 120:102795.\par 
\noindent Christofa, E., Ampountolas, K., and Skabardonis, A. (2016). Arterial traffic signal optimization: A person-based approach. Transportation Research Part C: Emerging Technologies, 66:27–47.\par 
\noindent Currie, G. (2006). Bus rapid transit in australasia: Performance, lessons learned and futures. Journal of Public Transportation, 9(3):1.\par 
\noindent Daganzo, C. F. (2007). Urban gridlock: Macroscopic modeling and mitigation approaches. Transportation Research Part B: Methodological, 41(1):49–62.\par 
\noindent Deng, T. and Nelson, J. D. (2011). Recent developments in bus rapid transit: a review of the
literature. Transport Reviews, 31(1):69–96.\par 
\noindent Dion, F. and Hellinga, B. (2002). A rule-based real-time traffic responsive signal control system with transit priority: application to an isolated intersection. Transportation Research Part B:
Methodological, 36(4):325–343.\par 
\noindent Eichler, M. and Daganzo, C. F. (2006). Bus lanes with intermittent priority: Strategy formulae and an evaluation. Transportation Research Part B: Methodological, 40(9):731–744.\par 
\noindent Garrow, M. and Machemehl, R. (1999). Development and evaluation of transit signal priority strategies. Journal of Public Transportation, 2(2):4.\par 
\noindent Hansen, P., Mladenovic, N., and Urosevic, D. (2006). Variable neighborhood search and local branching. Computers $\&$ Operations Research, 33(10):3034–3045.\par 
\noindent He, Q., Head, K. L., and Ding, J. (2012). Pamscod: Platoon-based arterial multi-modal signal
control with online data. Transportation Research Part C: Emerging Technologies, 20(1):164–184.\par 
\noindent He, Q., Head, K. L., and Ding, J. (2014). Multi-modal traffic signal control with priority, signal actuation and coordination. Transportation research part C: emerging technologies, 46:65–82.\par 
\noindent Hu, J., Park, B., and Parkany, A. E. (2014). Transit signal priority with connected vehicle technology. Transportation research record, 2418(1):20–29.\par 
\noindent Hu, J., Park, B. B., and Lee, Y.-J. (2015). Coordinated transit signal priority supporting transit progression under connected vehicle technology. Transportation Research Part C: Emerging
Technologies, 55:393–408.\par 
\noindent Levinson, H. S., Zimmerman, S., Clinger, J., and Gast, J. (2003). Bus rapid transit: Synthesis of case studies. Transportation Research Record, 1841(1):1–11.\par 
\noindent Li, M., Yin, Y., Zhang, W.-B., Zhou, K., and Nakamura, H. (2011). Modeling and implementation of adaptive transit signal priority on actuated control systems. Computer-Aided Civil and Infrastructure Engineering, 26(4):270–284.\par 
\noindent Liao, C.-F. and Davis, G. A. (2007). Simulation study of bus signal priority strategy: taking advantage of global positioning system, automated vehicle location system, and wireless communications. Transportation research record, 2034(1):82–91.\par 
\noindent Lin, X., Li, M., Yin, Y., Shen, Zuojun, M., and He, F. (2020). Rhythmic control of automated
traffic part ii: Grid network rhythm and online routing. Transportation Science, Forthcoming,
arXiv:2010.05416 [math.OC].\par 
\noindent Lin, X., Li, M., Yin, Y., Shen, Zuojun, M., and He, F. (2021). Rhythmic control of automated traffic on general road networks. Working paper.\par 
\noindent Ma, W., Head, K. L., and Feng, Y. (2014). Integrated optimization of transit priority operation at isolated intersections: A person-capacity-based approach. Transportation Research Part C:
Emerging Technologies, 40:49–62.\par 
\noindent Ma, W., Ni, W., Head, L., and Zhao, J. (2013). Effective coordinated optimization model for transit priority control under arterial progression. Transportation Research Record, 2366(1):71–83.\par 
\noindent Mesbah, M., Sarvi, M., and Currie, G. (2008). New methodology for optimizing transit priority at the network level. Transportation Research Record, 2089(1):93–100.\par 
\noindent Mladenovic, N. and Hansen, P. (1997). Variable neighborhood search. Computers $\&$ operations
research, 24(11):1097–1100.\par 
\noindent Smith, H. R., Hemily, B., and Ivanovic, M. (2005). Transit signal priority (tsp): A planning and implementation handbook.\par 
\noindent Sunkari, S. R., Beasley, P. S., Urbanik, T., and Fambro, D. B. (1995). Model to evaluate the impacts of bus priority on signalized intersections. Transportation Research Record, pages 117–123.\par 
\noindent Urbanik, T., Holder, R., and Fitzgerald, A. (1977). Evaluation of priority techniques for high occupancy vehicles on arterial streets. Technical report, Texas Transportation Institute, Texas
A $\&$ M University.\par 
\noindent Viegas, J. and Lu, B. (1997). Traffic control system with intermittent bus lanes. IFAC Proceedings Volumes, 30(8):865–870.\par 
\noindent Viegas, J. and Lu, B. (2004). The intermittent bus lane signals setting within an area. Transportation Research Part C: Emerging Technologies, 12(6):453–469.\par

\newpage

\setcounter{table}{0}
\renewcommand{\thetable}{A-\arabic{table}}

\section*{Appendix A: Notations}
\begin{table}[!ht]
  \centering
  \label{Nomenclature_table}
  \normalsize
  \begin{tabular}{ l  l }
   \hline
        \multicolumn{2}{l}{\textbf{Sets}} \\
        $\mathcal{A}$ & Link set\\
        $\mathcal{V}$ & Node set\\
        $\mathcal{V}_s$ & Node set of station\\
        $\mathcal{V}_o$ & Node set of origin\\
        $\mathcal{V}_d$ & Node set of destination\\
        $\mathcal{V}_{in}$ & Node set of intersection\\
        $\mathcal{Q}$ & Virtual platoon set\\
        $\mathcal{P}_{ij}$ & Bus set on link $(i,j)$\\
        $\mathcal{W}$ & O-D pair set\\
        $\mathcal{R}$ & Path set\\
   \hline
		\multicolumn{2}{l}{\textbf{Parameters}} \\
		$T$ & Cycle of RC scheme\\
        $H$ & Cycle of RC-H scheme\\
        $Q$ & Cycle number of RC for one RC-H cycle \\
        $\tau_i$ & Relative time of node $i$\\
        $\underline{t}_a^{ij}$ & Minimum travel time of private cars on link $(i,j)$\\
        $t_a^{ij}$ & Travel time of background virtual platoon on link $(i,j)$\\
        $t_{q,\hat{q}}^{ij}$ & Travel time of the realized virtual platoon of $q\rightarrow \hat{q}$ on link $(i,j)$\\
        $\underline{t}_{p}^{ij}$ & Minimum travel time of bus $p$ on link $(i,j)$\\
        $\underline{\hat{t}}_{p}^{j}$ & Minimum dwelling time of bus $p$ at station $j$\\
        $s_a^{ij}$ & Virtual platoon size on link $(i,j)$\\
        $s_b$ & Bus size\\
        $l_{ij}$ & Number of lanes on link $(i,j)$\\
        $\gamma_p$ & Number of passengers on bus $p$\\
        $d^w$ & Demand of O-D pair $w$\\
    \hline
		\multicolumn{2}{l}{\textbf{Decision variables}} \\
        $\Theta_{q,\hat{q}}^{ij}$ & Realized virtual platoon of $q\rightarrow \hat{q}$ on link $(i,j)$ \\
        $\hat{\Theta}_{q,\hat{q}}^{p,ij}$ & Dedicated virtual platoon of $q\rightarrow \hat{q}$ for bus $p$ on link $(i,j)$ \\
        $q_p^{ij}$ & Number of arrival platoon for bus $p$ on link $(i,j)$\\
        $\hat{q}_p^{ij}$ & Number of departure platoon for bus $p$ on link $(i,j)$\\
        $t_p^{ij}$ & Travel time of  bus $p$ on link $(i,j)$\\
        $\hat{t}_{p}^{j}$ & Dwelling time of bus $p$ at station $j$\\
        $f_r$ & Traffic flow of private cars on path $r$\\
        $\pi_{q,\hat{q}}^{r,ij}$ & Traffic of path $r$ of $q\rightarrow \hat{q}$ on mixed-traffic lanes of link $(i,j)$ \\
        $\tilde{\pi}_{q,\hat{q}}^{r,ij}$ & Traffic of path $r$ of $q\rightarrow \hat{q}$ on regular lanes of link $(i,j)$ \\        
        $\lambda_{q}^{r,ij}$ & Traffic of path $r$ of arriving with platoon $q$ on mixed-traffic lanes of link $(i,j)$ \\ 
        $\tilde{\lambda}_{q}^{r,ij}$ & Traffic of path $r$ of arriving with platoon $q$ on regular lanes of link $(i,j)$ \\   
        $\mu_{\hat{q}}^{r,ij}$ & Traffic of path $r$ of departing with platoon $\hat{q}$ on mixed-traffic lanes of link $(i,j)$ \\      
        $\tilde{\mu}_{\hat{q}}^{r,ij}$ & Traffic of path $r$ of departing with platoon $\hat{q}$ on regular lanes of link $(i,j)$ \\        
   \hline
  \end{tabular} \\
\end{table}


\end{document}